\documentclass[11pt,reqno]{amsart}

\usepackage{amsmath}\allowdisplaybreaks  
\usepackage{amssymb}        
\usepackage{amsfonts}       
\usepackage{amsthm}         
\usepackage{mathrsfs}       
\usepackage{cases}          
\usepackage{siunitx}        

\usepackage{array}
\usepackage{booktabs}       
\usepackage{longtable}      
\usepackage{multirow}
\usepackage{makecell}       
\usepackage{tabularx}       
\usepackage{nicematrix}     

\usepackage{tikz}\usetikzlibrary{tikzmark, arrows.meta}   
\usetikzlibrary{positioning}
\usepackage{float} 
\usepackage{xcolor}           
\usepackage{graphicx}         
\usepackage{subcaption}

\usepackage[linesnumbered, ruled, vlined]{algorithm2e}  

\usepackage[colorlinks=true, linkcolor=blue, citecolor=blue, urlcolor=blue]{hyperref} 
\usepackage[nameinlink]{cleveref}	

\usepackage[left=2cm, right=2cm, top=2cm, bottom=2cm]{geometry}

\usepackage{bm} 
\usepackage[normalem]{ulem}           
\usepackage{lineno}         

\makeatletter
\newcommand\figcaption{\def\@captype{figure}\caption} 
\newcommand\tabcaption{\def\@captype{table}\caption}
\makeatother

\newtheorem{theorem}{Theorem}[section]          
\newtheorem{example}[theorem]{Example}          
\newtheorem{remark}[theorem]{Remark}            

\begin{document}
	\title[DataTransfer: Neural network based interpolation across non-nested meshes]{DataTransfer: Neural network based interpolation across non-nested meshes}
	
	\author[J.~Hao, Y.~Huang, N.~Yi]{Jiaxiong Hao$^{\dag}$, Yunqing Huang$^\S$, Nianyu Yi$^\dag$}
	\address{$\dag$ Hunan Key Laboratory for Computation and Simulation in Science and Engineering, School of Mathematics and Computational Science, Xiangtan University, Xiangtan 411105, Hunan, P.R.China} \email{moodbear@qq.com (J. Hao);\  yinianyu@xtu.edu.cn (N. Yi)}
	\address{$\S$ National Center for Applied Mathematics in Hunan, Key Laboratory of Intelligent Computing \& Information Processing of Ministry of Education, Xiangtan University, Xiangtan 411105, Hunan, P.R.China} \email{huangyq@xtu.edu.cn}
	
	\begin{abstract}
	In mesh-based numerical simulations, the interpolation of mesh-defined functions across different meshes is a critical task, and achieving high-precision interpolation is of great significance for improving the computational efficiency and numerical stability of algorithms. 
	This paper proposes neural network based function mapping model across meshes, wherein the interpolation process is reformulated as a data-driven regression problem over scattered function data. 
	Conventional interpolation and projection-based approaches are highly dependent on mesh connectivity and corresponding geometric properties, which renders such methods computationally costly and sensitive to mismatches between source and target meshes. 
	The proposed method constructs a neural network approximator using nodal function values on the source mesh to obtain a global representation of the function, which can then be interpolated onto any other meshes.
To investigate the network architectural impacts on model performance, three representative feedforward network structures are numerically compared in this work: multi-layer perceptrons, extreme learning machines, and network incorporating radial basis function hidden units. 
The results reveal distinct trade-offs among  accuracy, computational efficiency and model robustness, among which the radial basis function-based network achieves the most desirable overall performance balance, enabling fast and precise function calculation. 
Numerical experiments conducted on non-nested meshes validate the efficacy of the proposed model in both function interpolation and cross-mesh data transmission tasks. 
    	The code used in this paper is publicly available on GitHub\footnotemark.
	\end{abstract}

	\footnotetext{GitHub: \url{https://github.com/FEMmaster/DataTransfer}.}
	
	
	
	\maketitle
\section{Introduction}
	Numerical simulations in engineering and applied sciences typically rely on mesh-based discretization of partial differential equations (PDEs) over complex domains. 
	In problems involving multiphysics, mesh adaptivity, or evolving geometries, data transfer between non-nested meshes is a critical prerequisite for ensuring solution accuracy, consistency, and numerical stability. 
	This challenge is particularly pronounced in multiphysics simulations, where different physical fields often demand meshes tailored to distinct geometric and resolution requirements, resulting in frequent data exchange between non-nested meshes. 
	A quintessential example is fluid–structure interaction: fluid domains require meshes capable of resolving boundary layers and accommodating large deformations, while structural domains prioritize geometric fidelity and sharp interface representation. 
	Similarly, in structural–thermal analysis, structural mechanics favors geometric conformity, whereas thermal diffusion simulations necessitate fine resolution in regions with steep temperature gradients \cite{Farhat1998}. 
	These challenges are further exacerbated in time-dependent simulations, where meshes are regenerated to track moving interfaces or evolving solution features. 
	While mesh adaption enhances accuracy and efficiency, it inherently entails repeated projection of solutions between successive non-nested meshes. 
	Beyond these applications, cross-mesh data transfer also arises naturally in numerical analysis, such as in multigrid methods \cite{Huang2006} and domain decomposition techniques \cite{Dolean2015}. 
	Existing techniques are often highly specialized, tailored to specific mesh types or problem settings, and lack a universal, efficient solution. 
	Consequently, accurate, stable, and property-preserving mesh-defined function interpolation across different meshes have become indispensable for reliable simulation workflows \cite{Bangerth2003, Nochetto2009}.
	
	Interpolation methods for transferring function data across different meshes can be broadly categorized into mesh-dependent and mesh-free approaches, depending on whether they explicitly rely on the underlying mesh connectivity. 	
	Representative mesh-dependent techniques include element-wise finite element interpolation, nearest-neighbor transfer, and local weighted averaging schemes. 
	All such methods rely on the geometric and topological relationships between the source and target meshes, which necessitates point-location queries to identify the host element within the source mesh for each target point. As the scale and complexity of the mesh increase, the point-location process becomes computationally prohibitive.
	Even with the adoption of spatial search structures such as kd-trees or octrees, geometric queries typically remain the dominant performance bottleneck. Beyond their substantial computational cost, these mesh-dependent methods are also sensitive to mesh quality, often exhibiting discontinuities or diminished accuracy in regions with large deformations \cite{Maddison2017, Harbecht2025}. A more rigorous representative of this category is the Galerkin projection method, which formulates the interpolation problem as a variational projection onto the target function space. However, its enhanced mathematical consistency comes at the cost of constructing a supermesh and performing computationally expensive numerical integration over intersecting subdomains.

In contrast, mesh-free methods eliminate the requirement for mesh connectivity by directly reconstructing approximation functions from scattered nodal data. Prominent examples include classical radial basis function (RBF) interpolation \cite{Li2002, Wendland2004, Bucelli2024}, moving least squares (MLS) \cite{Lancaster1981}, and global polynomial fitting. A key advantage of these approaches is  their ability to be evaluated at arbitrary target locations, make them particularly well-suited for dynamic, adaptive, or non-nested discretizations. Nevertheless, their practical implementation remains constrained by several inherent limitations.
Classical RBF interpolation centers the basis functions at the given date points, with expansion coefficients determined by enforcing interpolation conditions across all sample points. This gives rise to a dense, globally coupled linear system, whose computational cost and numerical conditioning degrade rapidly as the size of the data set increases--posing a significant challenge for large-scale simulations.
Moving least squares methods offer enhanced locality and flexibility compared to RBF interpolation, but their performance is limited by the need for repeated neighbor searches and local weight adjustments during the evaluation process. 
Global polynomial fitting provides a compact analytical representation but lacks the flexibility to resolve complex field features over irregular point distributions. Increasing the polynomial degree to improve flexibility often introduces spurious oscillations and compromises numerical stability. These challenges motivate the exploration of alternative function representation paradigms that combine the advantages of mesh-free flexibility and computational efficiency.

    A promising direction to address these limitations is to frame cross-mesh data transfer as a learned function approximation problem, leveraging neural networks to model the mapping between spatial coordinates and field values. Neural networks naturally excel at capturing complex nonlinear relationships, making them well-suited for cross-mesh interpolation tasks without relying on mesh connectivity information. 
    Early applications of neural network based interpolation primarily focused on sparse, scattered-data scenarios, with particular prominence in environmental and geoscientific research \cite{Rigol2001, Antoni2001, Qi2018, Yin2022, Zhang2020}.
In the field of scientific computing, feedforward neural networks (FNNs) have been widely employed for function reconstructing and interpolation from scattered samples, with applications spanning PDE surrogate modeling \cite{Chen2021} and physical field reconstruction \cite{Yang2023}. Llanas et al. \cite{Llanas2006} constructively demonstrated that single-hidden-layer feedforward networks can approximate any finite set of distinct sample points, thereby highlighting the inherent interpolation capability of shallow neural network architectures. Cao et al. \cite{Cao2010} further proved that exact interpolation networks can achieve near-optimal approximation performance for the target function.
Subsequently, Chen and Cao \cite{Chen2016} proposed neural network based interpolation and quasi-interpolation operators for scattered data, deriving rigorous error bounds for their approximation performance. Yu and Cao \cite{Yu2022} further developed modified four-layer feedforward neural networks as quasi-interpolation operators for scattered data, establishing direct, converse, and inverse approximation theorems. Additionally, extreme learning machines have been shown to effectively improve interpolation performance and suppress Runge-type oscillations \cite{Auricchio2024}. Collectively, these findings underscore the significant potential of neural network based interpolation as a learned extension of classical mesh-free reconstruction techniques, well-suited for non-nested mesh function transfer and scattered data processing.

A recent study \cite{hao2025nn} employed a standard fully connected neural network as a mesh-free surrogate to facilitate  finite element solutions transfer between successive non-nested meshes within an \(hr\)-adaptive finite element framework. This work not only demonstrates the practical viability of neural network based data transfer in adaptive finite element computations but also highlights its potential to address the limitations of classical methods. However, the selection of network architecture in such transfer models remains largely empirical, and its systematic impact on approximation accuracy, training efficiency, and computational cost has not been thoroughly investigated. 

To address the aforementioned limitations of traditional interpolation methods, we formalize cross-mesh data transfer as a scattered regression problem, eliminating the need for prior mesh connectivity information. Three distinct neural network architectures are systematically evaluated and compared across meshes with varying topologies and resolutions. We identify the Radial Basis Function-augmented Extreme Learning Machine (RBF-ELM) as the optimal solution, as it balances approximation accuracy, computational efficiency, and mesh agnosticism. The RBF-ELM model leverages localized Gaussian kernel feature maps and closed-form least-squares output weight learning, achieving ultra-high approximation precision, low computational cost, and strong robustness to non-uniform sampling and mesh topology variations. Comprehensive numerical experiments confirm that the RBF-ELM outperforms fully connected neural networks, standard extreme learning machines, and traditional piecewise linear interpolation in terms of accuracy, efficiency, and stability.
Our primary contributions are summarized as follows:
    \begin{itemize}
    \item A rigorous, mesh-agnostic mathematical formulation of the cross-mesh data transfer problem as a scattered data regression task, which eliminates the reliance on mesh connectivity information and element-wise geometric queries.

\item A systematic comparative analysis of three neural network architectures for mesh interpolation, quantifying their performance trade-offs in terms of approximation accuracy, computational cost, and robustness to non-uniform sampling and mesh topology variations.

\item The proposal of the RBF-ELM model, which integrates the strengths of radial basis function interpolation and extreme learning machines to mitigate the sampling sensitivity of global approximators, enabling fast, accurate, and robust data transfer across arbitrary unstructured and non-nested meshes.

\item Extensive numerical validation covering smooth oscillatory functions, multi-peak Gaussian functions, and weakly singular functions, as well as iterative mesh transfer tasks. The results demonstrate that the RBF-ELM model preserves field continuity and maintains negligible approximation error even after $200$ repeated cross-mesh transfers.  
    \end{itemize}
    
    The remainder of this paper is organized as follows. 
    \Cref{sec:INTER} is dedicated to neural network based interpolation for mesh-to-mesh data transfer on non-nested  meshes. 
    The problem formulation is first presented in \Cref{sec:PROBL}, followed by a detailed consideration of  three neural network models: \Cref{sec:MLP} describes a standard fully connected neural network trained via gradient descent, \Cref{sec:ELM} introduces an extreme learning machine framework featuring fixed random projections and closed-form output weights, and \Cref{sec:RBF} focuses on the integrating of radial basis functions and neural network architecture within the interpolation framework. 
    Representative numerical experiments for each model are conducted on both quasi-uniform and adaptive meshes. 
    Finally, \Cref{sec:APPLIC} investigates the RBF-ELM model in iterative mesh-to-mesh data transfer tasks, while a concise summary of the work and a discussion of potential future extensions are provided in \Cref{sec:CONCLU}. 
     
\section{Neural network based interpolation}\label{sec:INTER}
Data transfer between unstructured  and non-nested meshes is a prevalent challenge in scientific computing. A typical example arises in time-dependent simulations with adaptive finite element discretizations, where the solution at the previous time step must be transferred to the mesh employed at the current time step to continue the computation. This issue is even more pronounced in multiphysics problems, where different physical fields are discretized with the finite element spaces associated with an independently mesh, and require data exchange to ensure consistent coupling. 
Traditional interpolation approaches address this task primarily through mesh-dependent procedures that demand explicit geometric and topological information about the meshes. This dependence motivates an alternative perspective: by reformulating the transfer problem as a data-driven regression problem over scattered spatial samples, a global surrogate model for the underlying field can be constructed. Such a surrogate not only approximates the field accurately but also supports pointwise evaluation without recourse to mesh connectivity.
In the present work, this interpolation task is formulated within a finite element framework. A practical challenge in this context stems from two aspects: first, the non-uniform distribution of sample points inherited from adaptive meshes, where dense clusters of points appear in locally refined regions while other areas remain sparsely sampled; second, the intricate spatial structures of finite element approximations, which require the interpolation method to capture fine-grained features without introducing spurious oscillations. These factors collectively determine the need for a robust interpolation approach that can adapt to diverse mesh topologies and data distributions. We begin with a precise mathematical formulation of the problem.
	
\subsection{Problem Formulation}\label{sec:PROBL}
	Let $\Omega \subset \mathbb{R}^d$ denote a bounded physical domain. 
	We consider two finite element meshes of $\Omega$: a source mesh $\mathcal{T}_A = (\mathcal{N}_A, \mathcal{C}_A)$ and a target mesh $\mathcal{T}_B = (\mathcal{N}_B, \mathcal{C}_B)$, where $\mathcal{N}_A,\mathcal{N}_B$ denote the nodal sets and $\mathcal{C}_A,\mathcal{C}_B$ the corresponding element sets. 
	On the source mesh $\mathcal{T}_A$, we define the finite element space 
	\[
	V_h := \left\{ v \in C^0(\Omega):\left.v\right|_K \in P_k(K),\; \forall\ K \in \mathcal{T}_A \right\},
	\] 
    which is a subspace of $H^1(\Omega)$ consisting of globally continuous functions that restrict to polynomials of degree at most $k$ on each element.
	The finite element solution $u_h \in V_h$ admits the nodal representation 
	\[
	u_h(\mathbf{x}) = \sum_{i=1}^{N} u_i \phi_i(\mathbf{x}), \quad \mathbf{x} \in \Omega,
	\] 
    where $N := |\mathcal{N}_A|$, $u_i \in \mathbb{R}$ denotes the solution value at node $\mathbf{x}_i \in \mathcal{N}_A$, and the basis functions satisfy $\phi_j(\mathbf{x}_i)=\delta_{ji}$, so that $u_h(\mathbf{x}_i)=u_i$ for all $\mathbf{x}_i \in \mathcal{N}_A$. 
    In addition, $u_h$ is assumed to be available at interior quadrature points of each element in $\mathcal{C}_A$. 
    Values of $u_h$ at nodal and quadrature points are readily available because these locations are fixed with respect to the mesh. 
    At arbitrary spatial locations, however, evaluating $u_h$ requires a nontrivial element search.
	
	To circumvent this geometric search, a globally continuous surrogate $\tilde{u} : \Omega \rightarrow \mathbb{R}$ is introduced to approximate the discrete finite element solution $u_h(\mathbf{x})$ across the domain: 
	\[
	\tilde{u}(\mathbf{x}) \approx u_h(\mathbf{x}), \quad \forall\ \mathbf{x} \in \Omega.
	\]
	This surrogate permits evaluation at arbitrary spatial locations without requiring any geometric correspondence or mesh connectivity between $\mathcal{T}_A$ and $\mathcal{T}_B$.
	Its construction is based on a discrete training dataset extracted from pointwise samples of the finite element solution $u_h$, 
	\[
	\mathcal{D}_A = \{(\mathbf{x}_i,u_h(\mathbf{x}_i)): \mathbf{x}_i \in \mathcal{S}_A\},
	\] 
	where the sample set $\mathcal{S}_A$ is typically chosen as the nodal set $\mathcal{N}_A$, possibly augmented with additional interior points such as quadrature points. 
	Based on this, the objective is to construct a parametric regression model $\tilde{u}_{\theta}(\mathbf{x})$ with parameters $\theta$, such that 
	\[
	\tilde{u}_{\theta}(\mathbf{x})\approx u_h(\mathbf{x}), \quad \forall\ \mathbf{x}\in\Omega.
	\] 
	The model parameters $\theta$ are chosen by fitting $\tilde{u}_{\theta}$ to the sampled data in $\mathcal{D}_A$, for instance through the empirical mean squared error (MSE): 
	\[
	\min_{\theta}\frac{1}{|\mathcal{S}_A|}\sum_{\mathbf{x}_i \in \mathcal{S}_A}\bigl(\tilde{u}_{\theta}(\mathbf{x}_i)-u_h(\mathbf{x}_i)\bigr)^2.
	\]
    Depending on the model class, this fitting step may be carried out either through iterative optimization algorithms or through the direct solution of a linear least-squares problem. 
    In the following sections, we introduce several representative neural-network models for mesh-to-mesh transfer on non-nested meshes and compare their accuracy, efficiency and robustness under such sampling distributions. 
    

	
\subsection{Data Transfer via Multilayer Perceptron}\label{sec:MLP}
\begin{figure}[htbp]
		\includegraphics[width=0.7\textwidth]{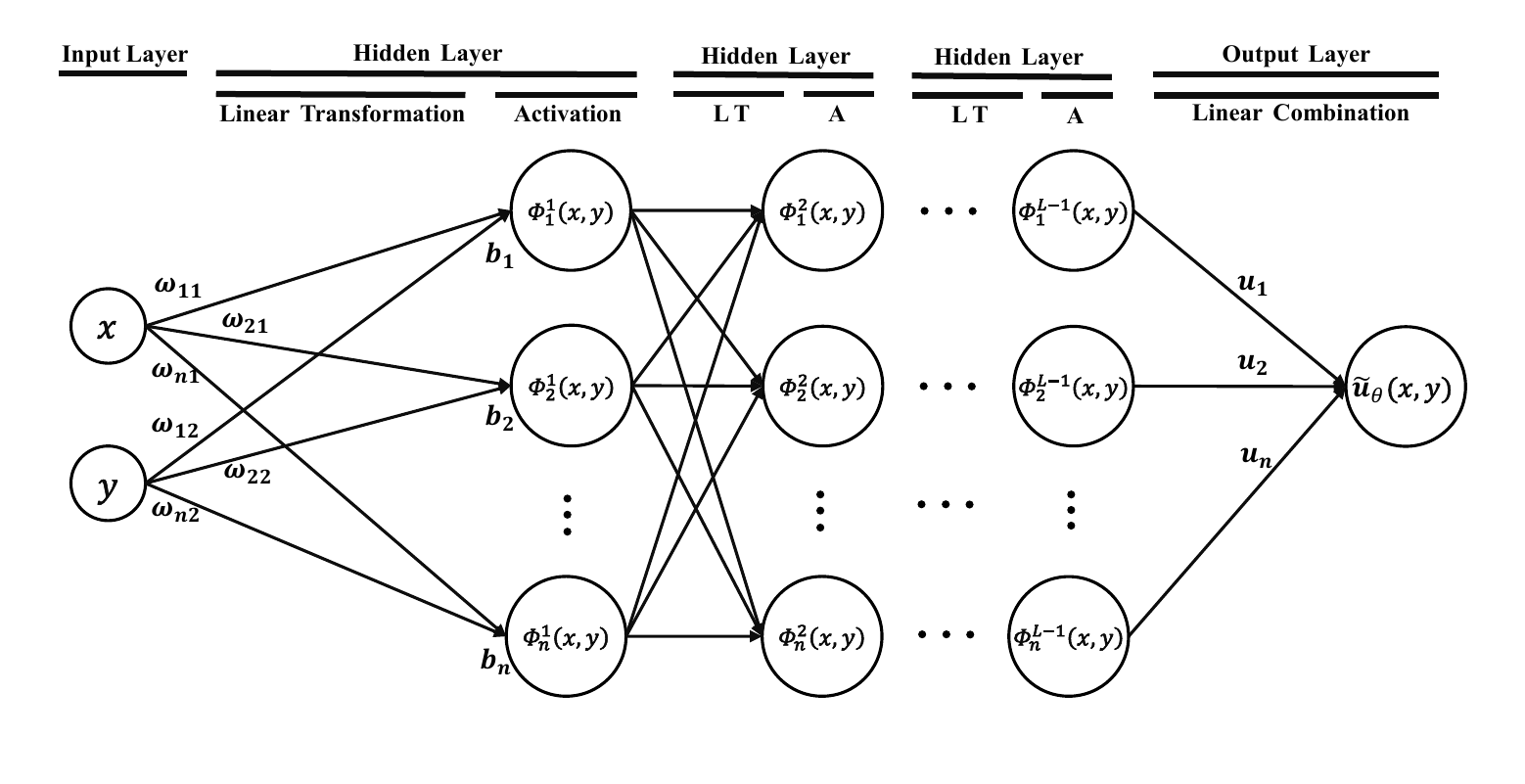}
		\caption{Schematic architecture of a fully connected feedforward neural network.}
		\label{fig:NN}
	\end{figure}

	Motivated by the universal approximation theorem~\cite{Cybenko1989}, we employ a standard fully connected multilayer perceptron (MLP) as a baseline model for comparison in mesh-to-mesh data transfer. 
	A general $L$-layer feedforward neural network is defined as a parameterized mapping $u_{\theta}: \mathbb{R}^{d_0} \to \mathbb{R}^{d_L}$. 
	For each hidden layer $k=1,\dots,L-1$, the transformation is given by 
	\[ 
	\mathbf{x}^{(k)} = \bm{\Phi}^{(k)}(\mathbf{x}^{(k-1)}) =  \sigma\!\left(\bm{W}^{(k)}\mathbf{x}^{(k-1)} + \bm{b}^{(k)}\right), 
	\] 
	where $\sigma$ is an activation function, $\bm{W}^{(k)} \in \mathbb{R}^{d_k \times d_{k-1}}$ and $\bm{b}^{(k)} \in \mathbb{R}^{d_k}$ are trainable weight matrices and bias vectors. 
    The nonlinear activation function $\sigma$ is essential for expanding the class of functions that the network can represent beyond affine models. 
    In the final layer, activation is omitted so as to produce unrestricted real-valued outputs suitable for regression. 
    Let $\mathbf{x}^{(0)} := \mathbf{x}$ denote the input. 
	The network can be expressed as
	\[ 
	\tilde{u}_{\theta}(\mathbf{x}) = \bm{W}^{(L)} \left(\bm{\Phi}^{(L-1)} \circ \cdots \circ \bm{\Phi}^{(2)} \circ \bm{\Phi}^{(1)}(\mathbf{x}) \right) + \bm{b}^{(L)}.
	\]
	Given training samples $\{(\mathbf{x}_i, u(\mathbf{x}_i))\}_{i=1}^m$ drawn from the source mesh 
	(e.g., nodal or quadrature points), the parameters are determined by minimizing the empirical mean squared error
	\[
	\theta^* = \arg\min_{\theta} \frac{1}{m}\sum_{i=1}^m \left(u_{\theta}(\mathbf{x}_i) - u(\mathbf{x}_i)\right)^2 .
	\]
	The structure illustrated in \Cref{fig:NN} reflects the composition described above, comprising $L-1$ hidden layers followed by a linear output layer. 
	For the case of a two-dimensional input $(x,y)$, the first-layer weight matrix and bias vector take the form $\bm{W}^{(1)} = [w_{ij}] \in \mathbb{R}^{n \times 2}$ and $\bm{b}^{(1)} = [b_1, \dots, b_n]^T$, producing the hidden-layer output $\bm{\Phi}^{(1)} = [\phi^1_1(x,y), \dots, \phi^1_n(x,y)]^T$.
	The final output layer yields a weighted combination of the hidden-layer activations, 
	\[
	\tilde{u}_{\theta}(x,y) = \bm{W}^{(L)} \bm{\Phi} + b^{(L)}, 
	\quad \bm{W}^{(L)} = [u_1, \dots, u_n].
	\]
	When the output bias is omitted, so that $b^{(L)} = 0$, this representation is algebraically identical to the coefficient-basis form of a finite element solution. 
	
	Throughout this section and the subsequent ones, interpolation accuracy is assessed on the target mesh $\mathcal{T}_B$ by two discrete error metrics, namely the mean absolute error (MAE) and the relative $L_2$ error (RL2). 
    They are defined by
    \begin{equation*}
    	\text{MAE} = \frac{1}{M} \sum_{i=1}^{M} \left| \tilde{u}_{\theta}(\mathbf{x}_i) - u(\mathbf{x}_i) \right|, \qquad 
    	\text{RL2} = \frac{\sqrt{\sum_{i=1}^M \bigl(\tilde{u}_{\theta}(\mathbf{x}_i) - u(\mathbf{x}_i)\bigr)^2}}{\sqrt{\sum_{i=1}^M \bigl(u(\mathbf{x}_i)\bigr)^2}},
    \end{equation*}
    where $\{\mathbf{x}_i\}_{i=1}^M$ denotes the evaluation points on the target mesh. 
    In all experiments below, the evaluation points are taken to be the nodes of the target mesh, and the reference values are obtained by evaluating the exact solution there. 
	We begin the numerical study with tests on a quasi-uniform mesh and then turn to adaptive mesh settings. 
	
\begin{figure}[htbp]
    \centering
        \begin{subfigure}[b]{0.25\textwidth}
            \includegraphics[width=\textwidth]{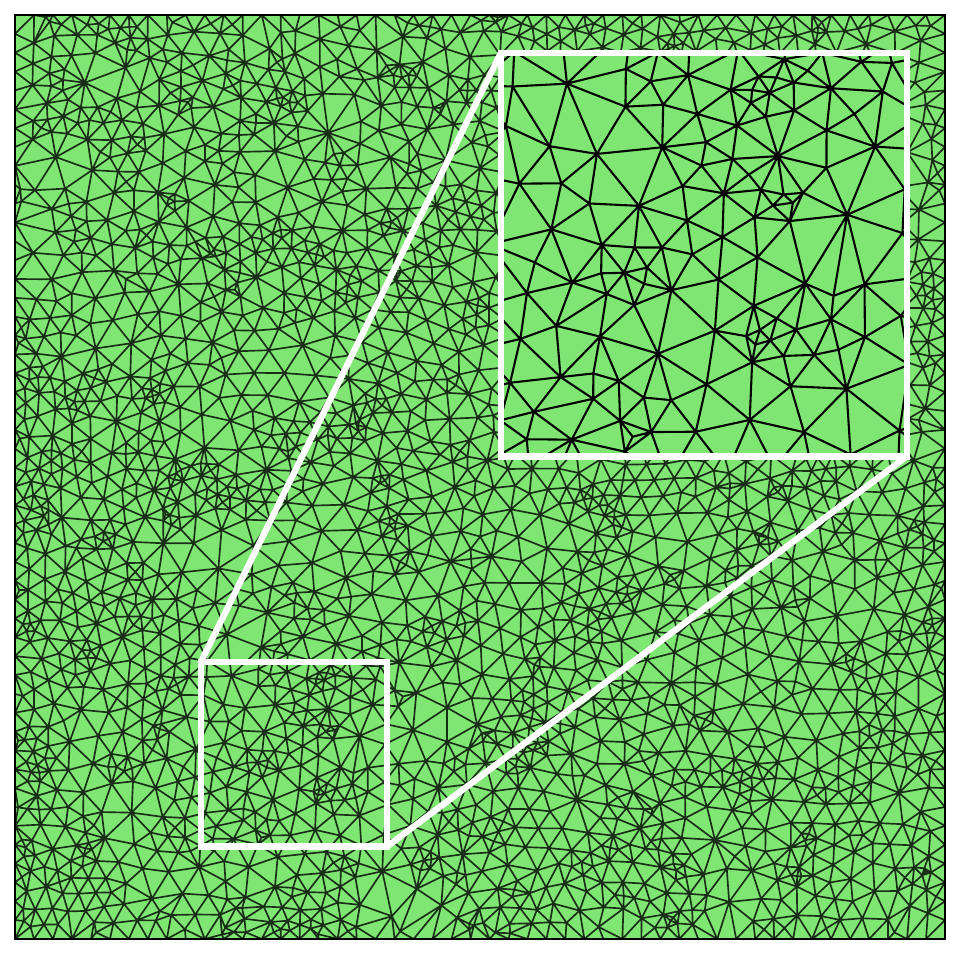}
            \caption{}
            \label{fig:NN-train1-A}
        \end{subfigure}
        \hspace{0.03\textwidth}
        \begin{subfigure}[b]{0.35\textwidth}
            \includegraphics[width=0.85\textwidth]{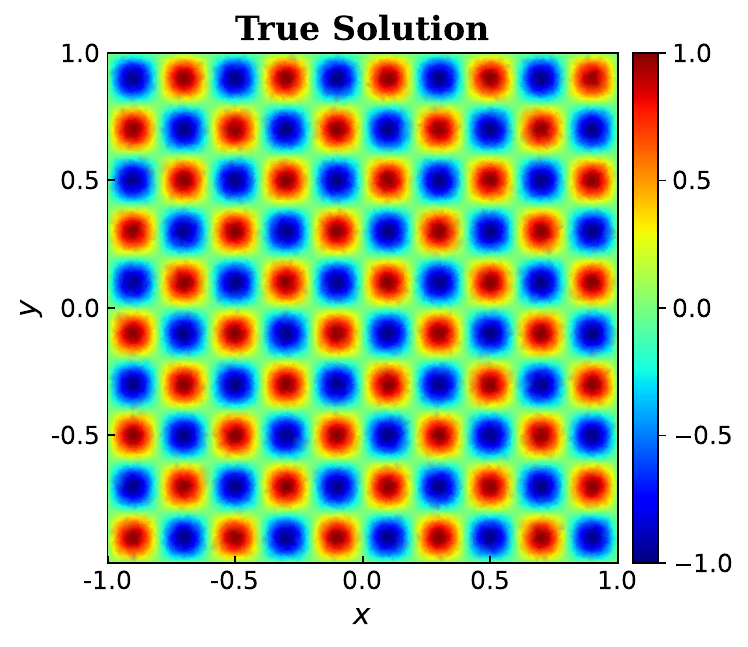}
            \caption{}
            \label{fig:NN-train1-B}
        \end{subfigure}	
        \begin{subfigure}[b]{0.32\textwidth}
            \includegraphics[width=0.9\textwidth]{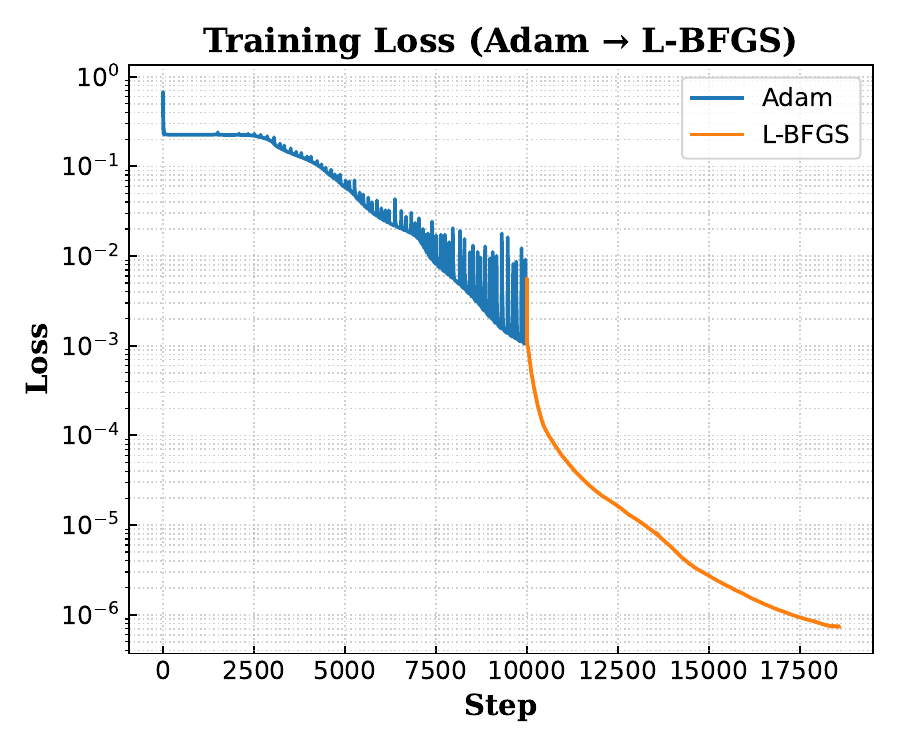}
            \caption{}
            \label{fig:NN-train1-C}
        \end{subfigure}
        \begin{subfigure}[b]{0.25\textwidth}
            \includegraphics[width=\textwidth]{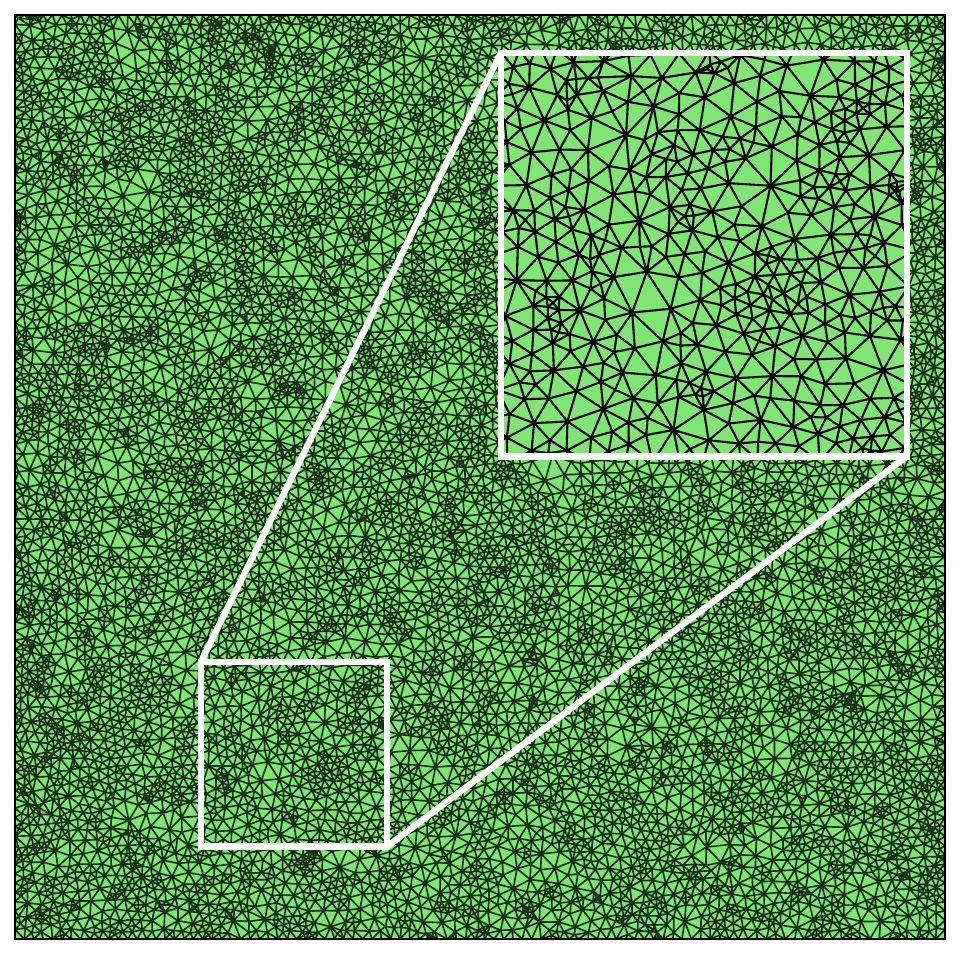}
            \caption{}
            \label{fig:NN-test1-D}
        \end{subfigure}
        \hspace{0.03\textwidth}
        \begin{subfigure}[b]{0.68\textwidth}
            \includegraphics[width=0.92\textwidth]{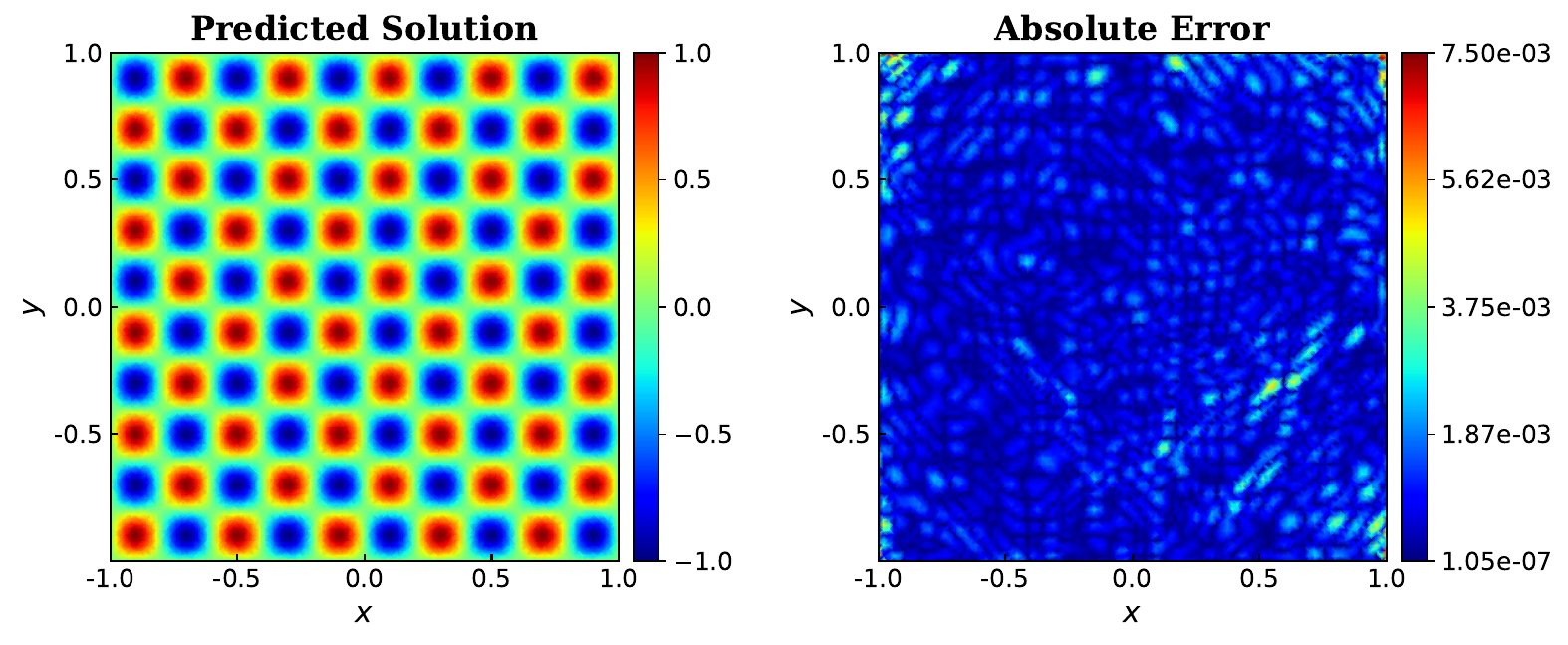}
            \caption{}
            \label{fig:NN-test1-E}
        \end{subfigure}
        \caption{\Cref{ex1}, NN interpolation on a quasi-uniform mesh. 
            (A) training mesh containing 2,500 nodes; 
            (B) exact solution sampled on the training mesh; 
            (C) history of the training loss; 
            (D) test mesh with 10,000 evaluation points; 
            (E) interpolated solution together with the corresponding error distribution.}
        \label{fig:NN-test1}
    \end{figure}

	\begin{example}\label{ex1}
        Consider the high-frequency oscillatory function
		\[
		u(x,y) = \sin(5 \pi x)\sin(5 \pi y), \quad (x,y)\in[-1,1]^2.
		\]
		We use a multilayer perceptron with three hidden layers of width 128, corresponding to the architecture $[2,128,128,128,1]$. 
		In this example, the network uses the $\tanh$ activation function and all trainable parameters are initialized with the Xavier uniform scheme. 
        A hybrid optimization strategy is used for training, with Adam providing rapid initial descent and L-BFGS supplying subsequent refinement. 
        Training is terminated once the loss stagnates or a prescribed iteration budget is reached. 
        This combination is found empirically to give more reliable convergence than either method alone. 
        \Cref{fig:NN-test1} collects the overall training and evaluation results for this example. 
	\end{example}
    The training dataset consists of the nodal points of the quasi-uniform mesh shown in \Cref{fig:NN-train1-A}, which contains approximately $2,500$ nodes. 
    \Cref{fig:NN-train1-B} shows the exact solution sampled on this mesh. 
    Since the mesh resolution is still limited, the sampled field captures the main oscillatory pattern but does not yet fully reveal its finer structure. 
    In \Cref{fig:NN-train1-C}, the blue segment of the loss curve corresponds to the Adam stage, while the orange segment corresponds to the subsequent L-BFGS refinement. 
    It is evident that the oscillations become increasingly pronounced in the later stage of Adam, whereas after switching to L-BFGS the loss resumes a stable decay. 
    The training process reaches a nearly stationary state after about 18,000 iterations. 
    Beyond this point, the loss remains at the level of $10^{-6}$ and further training incurs substantially higher cost while yielding only marginal improvement. 
    Interpolation performance is further assessed on the finer target meshes shown in \Cref{fig:NN-test1-D}. 
    The resulting predictions and error distributions are presented in \Cref{fig:NN-test1-E}. 
    These results show that, even though the training data are available only at comparatively sparse sample points, the learned network yields smooth predictions on finer meshes that closely reproduce the overall profile of the exact solution. 
    The reconstruction is therefore not a simple linear interpolation between neighboring samples, but a smooth approximation of the underlying field. 
    In our tests, the MAE remains at the level of $10^{-3}$ across the hyperparameter settings considered. 
    Further enlarging the training set mainly improves visual smoothness, with little further reduction in the error. 
    
    In contrast to the quasi-uniform mesh considered above, many practical simulations rely on adaptive mesh refinement to balance accuracy and efficiency. 
    Such refinement increases resolution in regions of rapid variation or large estimated error while keeping a coarser discretization elsewhere, thereby producing highly nonuniform node distributions. 
    This setting is particularly relevant for target functions with localized gradients or singular behavior. 
    The next example turns to this case and considers a target function given by a finite element solution on an adaptively refined mesh. 
	
\begin{figure}[htbp]
    \centering
        \begin{subfigure}[b]{0.25\textwidth}            
            \includegraphics[width=\textwidth]{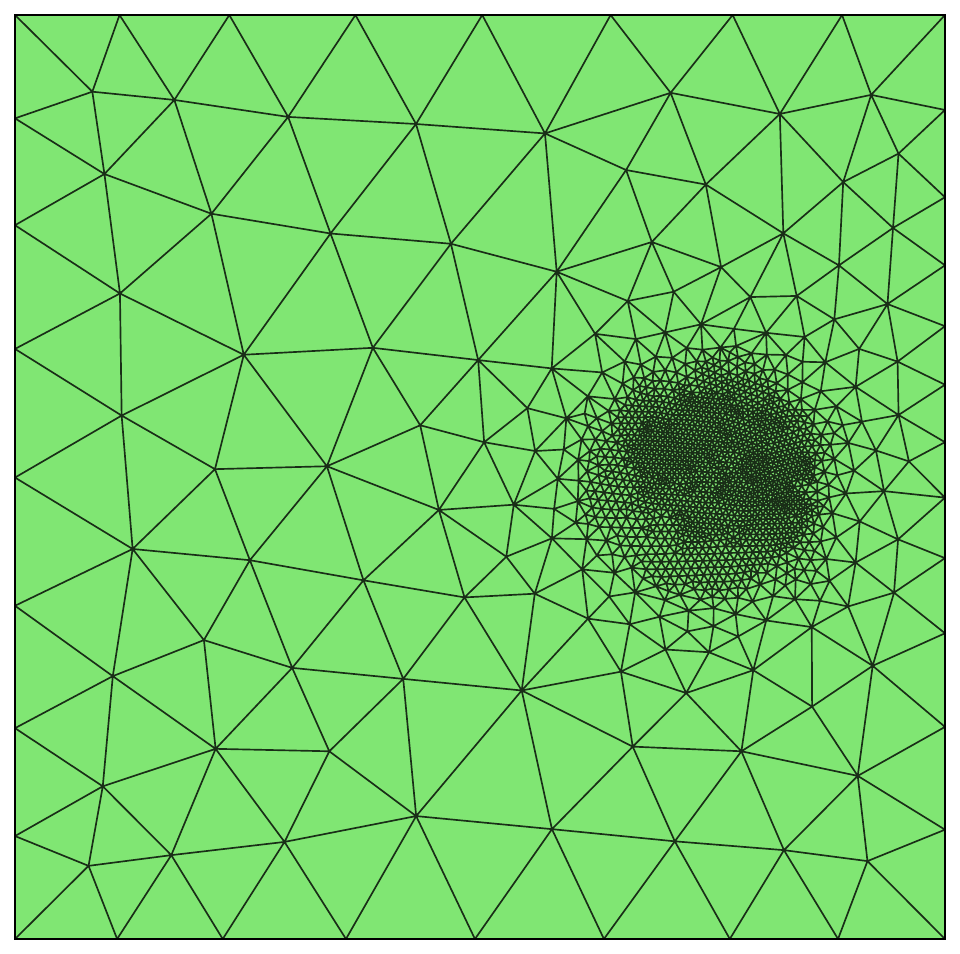}
            \caption{}
            \label{fig:NN-train2-A}
        \end{subfigure}
        \hspace{0.03\textwidth}
        \begin{subfigure}[b]{0.35\textwidth}
            \includegraphics[width=0.85\textwidth]{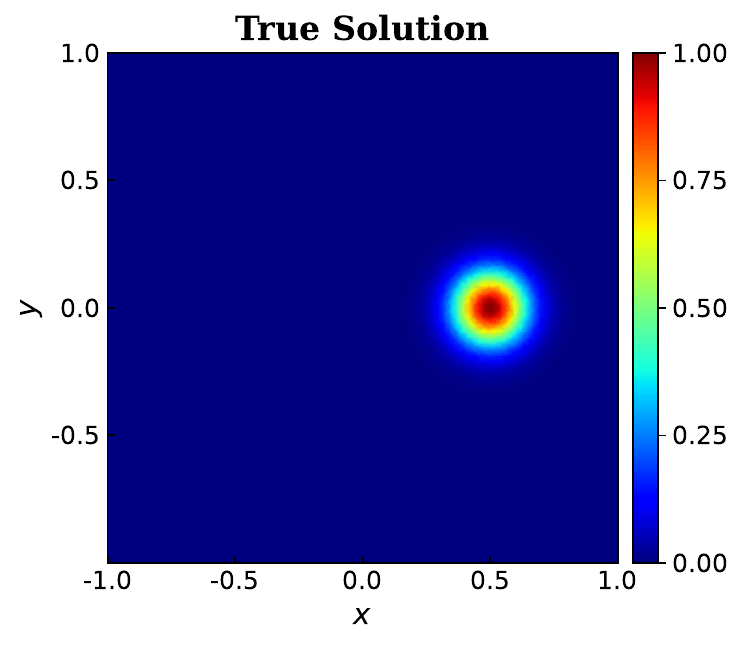}
            \caption{}
            \label{fig:NN-train2-B}
        \end{subfigure}		
        \begin{subfigure}[b]{0.32\textwidth}
            \includegraphics[width=0.9\textwidth]{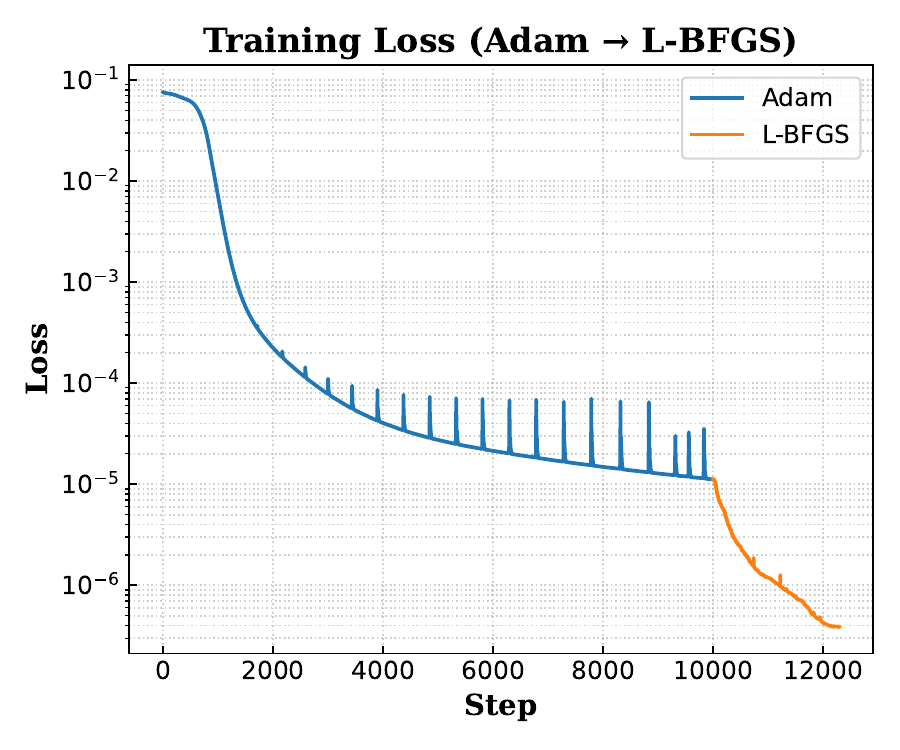}
            \caption{}
            \label{fig:NN-train2-C}
        \end{subfigure}
        \begin{subfigure}[b]{0.25\textwidth}
            \includegraphics[width=\textwidth]{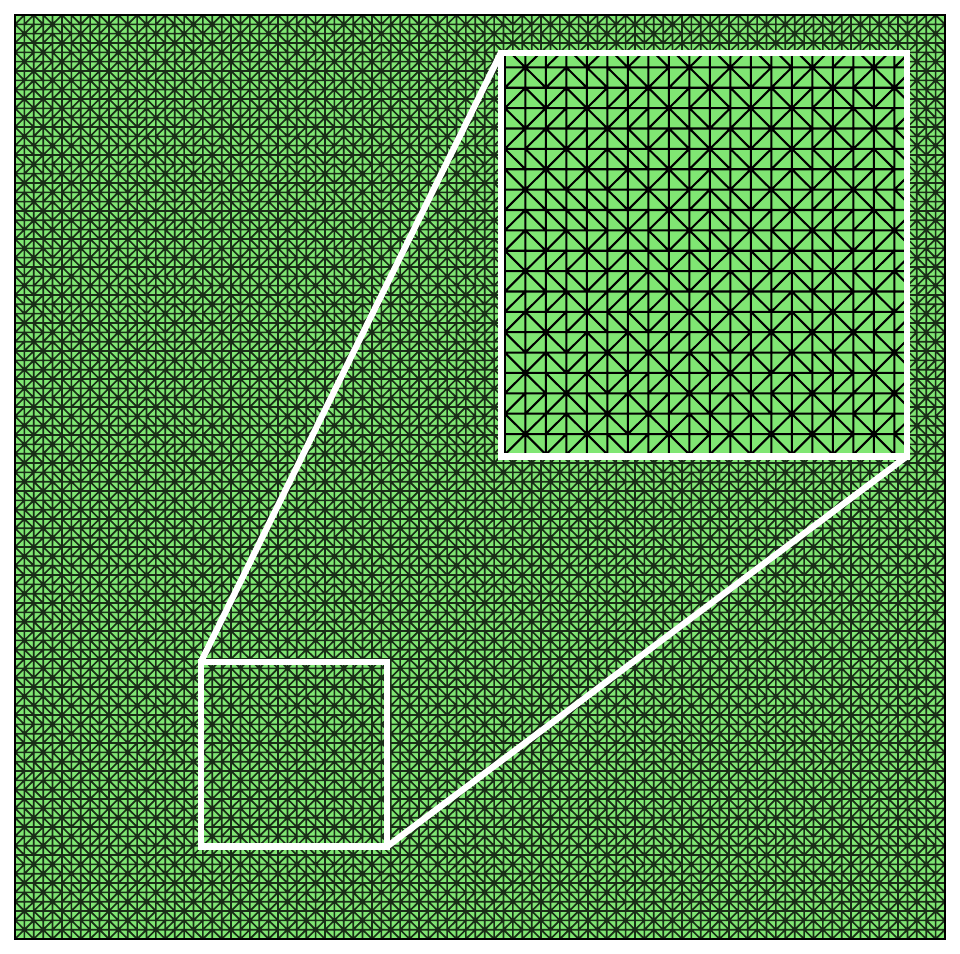}
            \caption{}
            \label{fig:NN-test2-D}
        \end{subfigure}
        \hspace{0.03\textwidth}
        \begin{subfigure}[b]{0.68\textwidth}
            \includegraphics[width=0.92\textwidth]{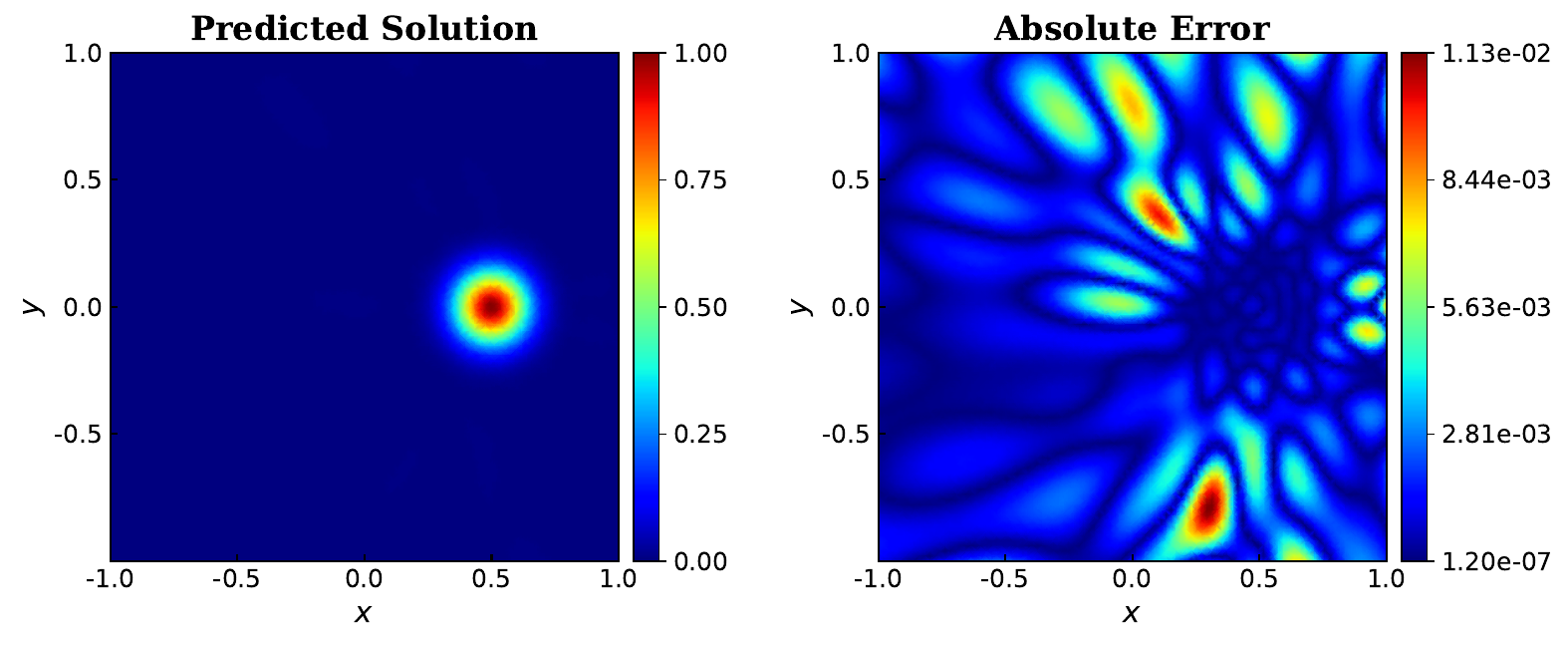}
            \caption{}
            \label{fig:NN-test2-E}
        \end{subfigure}
        \caption{\Cref{ex2}, NN interpolation on an adaptive mesh. 
            (A) training mesh; 
            (B) reference solution; 
            (C) history of the training loss; 
            (D) test mesh with 10,000 evaluation points; 
            (E) interpolated solution together with the corresponding error distribution.}
        \label{fig:NN-test2}
    \end{figure}

	\begin{example}\label{ex2}
		Let $\Omega = [-1,1]^2$ and $I = (0,1]$. 
        Consider the heat equation with Dirichlet boundary conditions
		\[
		\begin{cases}
			\partial_t u - \Delta u = f, & \text{in }\Omega \times I,\\
			u(\mathbf{x},t) = g(\mathbf{x},t), & \text{on }\partial\Omega\times I, \\
			u(\mathbf{x},0) = u_0(\mathbf{x}), & \text{in } \Omega,
		\end{cases}
		\]
		where the forcing term $f$ is chosen so that the exact solution is given by
		\[
		u(x, y, t) = \exp\left(-50\left((x - 0.5 \cos(2\pi t))^2 + (y - 0.5 \sin(2\pi t))^2\right)\right).
		\]
        This exact solution is used as the reference in the subsequent error evaluation. 
        The boundary and initial conditions are taken consistently from the exact solution, namely 
        \[g(x,y,t) = u(x,y,t)\big|_{(x,y)\in\partial\Omega}, \quad  u_0(x,y) = u(x,y,0). \]
        Adaptive finite element discretization of this problem (see, e.g., \cite{hao2025nn}) yields a sequence of numerical solutions together with the corresponding adaptively refined meshes. 
        One representative snapshot, along with its associated mesh $\mathcal{T}_A$, is then extracted to define the training dataset $\mathcal{D}_A$. 
		All hyperparameters are kept the same as in \Cref{ex1}, except that the network architecture is reduced to $[2,64,64,64,1]$. 
        A detailed illustration is shown in \Cref{fig:NN-test2}. 
	\end{example}
    \Cref{fig:NN-train2-A,fig:NN-train2-B} illustrate a representative adaptive mesh and the corresponding reference solution. 
    The nodal values of this discrete solution are used as the training data for the neural-network surrogate, whose loss history is shown in \Cref{fig:NN-train2-C}. 
    After 10,000 Adam iterations, the loss decreases to the order of $10^{-5}$ with noticeable oscillations and subsequent L-BFGS refinement reduces the training error further. 
    The trained model is evaluated on the uniform mesh shown in \Cref{fig:NN-test2-D}, and the corresponding prediction and error distribution are displayed in \Cref{fig:NN-test2-E}. 
    Although the predicted solution is visually close to the reference, the error distribution indicates that the achieved accuracy remains at the level of $10^{-3}$. 
        
    Although MLPs possess universal approximation capability, the experiments above show that their practical use for mesh-based interpolation remains challenging. 
    While the networks produce smooth reconstructions, the achieved accuracy typically stagnates near $10^{-3}$, leaving limited room for further improvement. 
    More generally, for multilayer fully trainable deep networks, the nonconvex nature of the underlying optimization problem makes gradient-based training computationally demanding, and increased depth offers no guaranteed performance gain. 
    In addition, practical deployment is complicated by the need to tune multiple hyperparameters, including network width, depth, initialization method, optimization algorithm and learning rate. 
    For interpolation tasks, however, accuracy, efficiency and robustness are all essential. 
    These considerations motivate the alternative strategy introduced in the next section, where we turn to extreme learning machines for faster training and more stable performance. 
	
\subsection{Data Transfer via Extreme Learning Machine}\label{sec:ELM}
    As a computationally efficient alternative to fully trainable deep neural networks, the extreme learning machine (ELM) is adopted to construct the surrogate function $\tilde{u}_{\theta}(\mathbf{x})$. 
    Its defining feature is that the hidden-layer parameters are randomly initialized and then fixed, so that only the output weights need to be determined. 
    As a result, training is reduced from a nonconvex optimization problem over the full parameter space to a convex least-squares problem for the output layer. 
    This avoids iterative optimization over the full parameter space and makes ELM an attractive choice for mesh-to-mesh interpolation tasks where efficiency and simplicity are important. 
	
\begin{figure}[htbp]
		\includegraphics[width=0.7\textwidth]{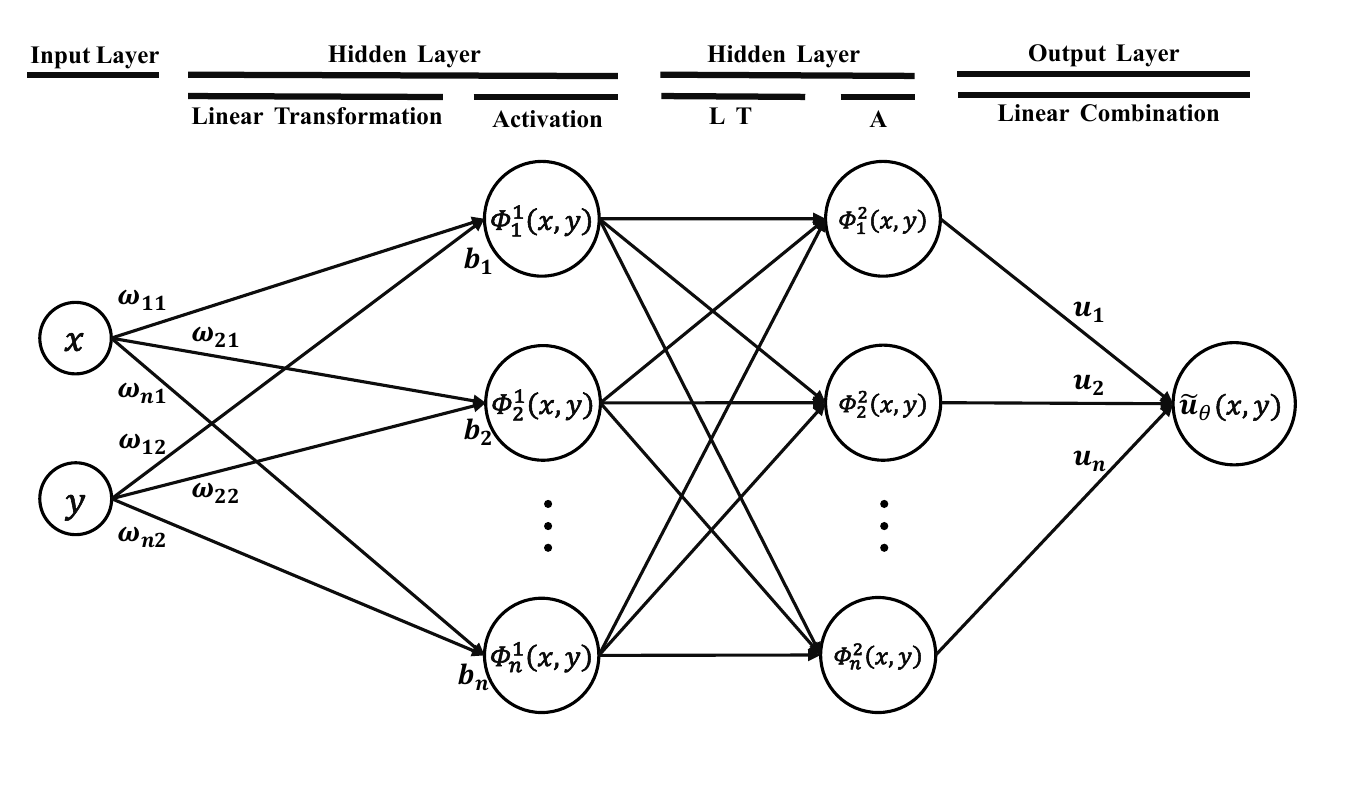}
		\caption{Schematic architecture of an extreme learning machine network.}
		\label{fig:ELM}
	\end{figure}

    For a given input $\mathbf{x} \in \mathbb{R}^{d_0}$, the model performs a sequence of nonlinear transformations followed by a linear output mapping:
	\begin{align*}
		\bm{\Phi}^{(1)}(\mathbf{x}) &= \sigma\left(\bm{W}^{(1)} \mathbf{x} + \bm{b}^{(1)}\right), \\
		\bm{\Phi}^{(2)}(\mathbf{x}) &= \sigma\left(\bm{W}^{(2)} \bm{\Phi}^{(1)}(\mathbf{x}) + \bm{b}^{(2)}\right), \\
		\tilde{u}_{\theta}(\mathbf{x}) &= \bm{\beta}^\top \bm{\Phi}^{(2)}(\mathbf{x}),
	\end{align*}
	where $\bm{W}^{(1)} \in \mathbb{R}^{d_1 \times d_0}$ and $\bm{W}^{(2)} \in \mathbb{R}^{d_2 \times d_1}$ denote the weight matrices, and $\bm{b}^{(1)} \in \mathbb{R}^{d_1}$ and $\bm{b}^{(2)} \in \mathbb{R}^{d_2}$ denote the bias vectors. 
	Nonlinearity in the hidden layers is introduced through the element-wise activation function $\sigma(\cdot)$, which acts on each affine transformation to generate nonlinear hidden features. 
    At the output layer, the coefficient vector $\bm{\beta}\in\mathbb{R}^{d_2}$ forms the final linear combination of the hidden features $\bm{\Phi}^{(2)}(\mathbf{x})$.
    Given training samples $\{(\mathbf{x}_i, u(\mathbf{x}_i))\}_{i=1}^m$ drawn from the source mesh (e.g., nodal or quadrature points), the parameters $\{\bm{W}^{(1)}, \bm{W}^{(2)}, \bm{b}^{(1)}, \bm{b}^{(2)}\}$ are randomly initialized once and then kept fixed throughout training.
    Consequently, training no longer involves nonconvex optimization over all network parameters, but reduces to a convex optimization problem for $\bm{\beta}$.     
    As shown in \Cref{fig:ELM}, the architecture retains the neural-network layout introduced previously, while its depth is now taken to be two hidden layers. 
    \begin{remark}
		While traditional ELMs typically employ a single hidden layer, the network in this work adopts two hidden layers in order to enhance approximation accuracy while maintaining training stability. 
		This design is motivated by empirical observations rather than theoretical guarantees. 
		In practice, single-layer networks were insufficiently expressive for mesh-based interpolation, whereas additional depth beyond two layers provided no consistent benefit and in some cases led to unstable behavior in training. 
	\end{remark}

    Let $\bm{\Phi}(\mathbf{x}) \in \mathbb{R}^{d_2}$ denote the composite nonlinear mapping of the input $\mathbf{x}$ through the hidden layers. 
    The ELM representation then takes the simplified form
	\begin{equation*}
		\tilde{u}_{\theta}(\mathbf{x}) = \bm{\beta}^\top \bm{\Phi}(\mathbf{x}).
	\end{equation*}  
    This representation is again in the coefficient-basis form of a finite element solution, as noted earlier for the MLP case in \Cref{sec:MLP}. 
    Under this representation, evaluating the feature map $\bm{\Phi}(\mathbf{x}_i)$ at the training points yields the hidden-layer output matrix $H \in \mathbb{R}^{N \times d_2}$, while the corresponding sampled values form the target vector $U \in \mathbb{R}^N$, namely
	\begin{equation*}
		H = \begin{bmatrix}
			\bm{\Phi}(\mathbf{x}_1)^\top \\
			\bm{\Phi}(\mathbf{x}_2)^\top \\
			\vdots \\
			\bm{\Phi}(\mathbf{x}_N)^\top
		\end{bmatrix}, \quad
		U = \begin{bmatrix}
			u_h(\mathbf{x}_1) \\
			u_h(\mathbf{x}_2) \\
			\vdots \\
			u_h(\mathbf{x}_N)
		\end{bmatrix}.
	\end{equation*}
	Accordingly, the training of ELM reduces to finding the output weights $\bm{\beta}$ that minimize the discrepancy between $H\bm{\beta}$ and $U$. 
    This is formulated as the standard least-squares problem:
	\begin{equation*}
		\bm{\beta} = \arg\min_{\bm{\beta}} \|H \bm{\beta} - U\|_2^2.
	\end{equation*}
	If $H^\top H$ is nonsingular, this problem admits the explicit solution
    \begin{equation*}
    	\bm{\beta} = (H^\top H)^{-1} H^\top U.
    \end{equation*}
	In practice, however, the matrix $H^\top H$ may become ill-conditioned or even singular when the number of hidden nodes $d_2$ is large or when the columns of $H$ are strongly correlated. 
	Standard remedies include the Moore-Penrose pseudoinverse \cite{Huang2006}, regularization techniques such as ridge regression \cite{Hoerl2000, MartnezMartnez2011} and robust variants of ELM \cite{Horata2013}. 
    Among these approaches, the pseudoinverse is often the simplest and most commonly used in practice, leading to the representation
    \begin{equation} \label{eq:beta}
    	\bm{\beta} = H^\dagger U,
    \end{equation}
    where $H^\dagger$ denotes the Moore-Penrose pseudoinverse of $H$. 
    This direct matrix-based computation avoids iterative optimization and keeps both training time and evaluation cost low, making ELM a particularly attractive surrogate model when efficiency is critical.
	In the examples below, the ELM models are trained using this formulation. 
    
\begin{figure}[htbp]
    \centering
        \begin{subfigure}[b]{0.25\textwidth}         
            \includegraphics[width=\textwidth]{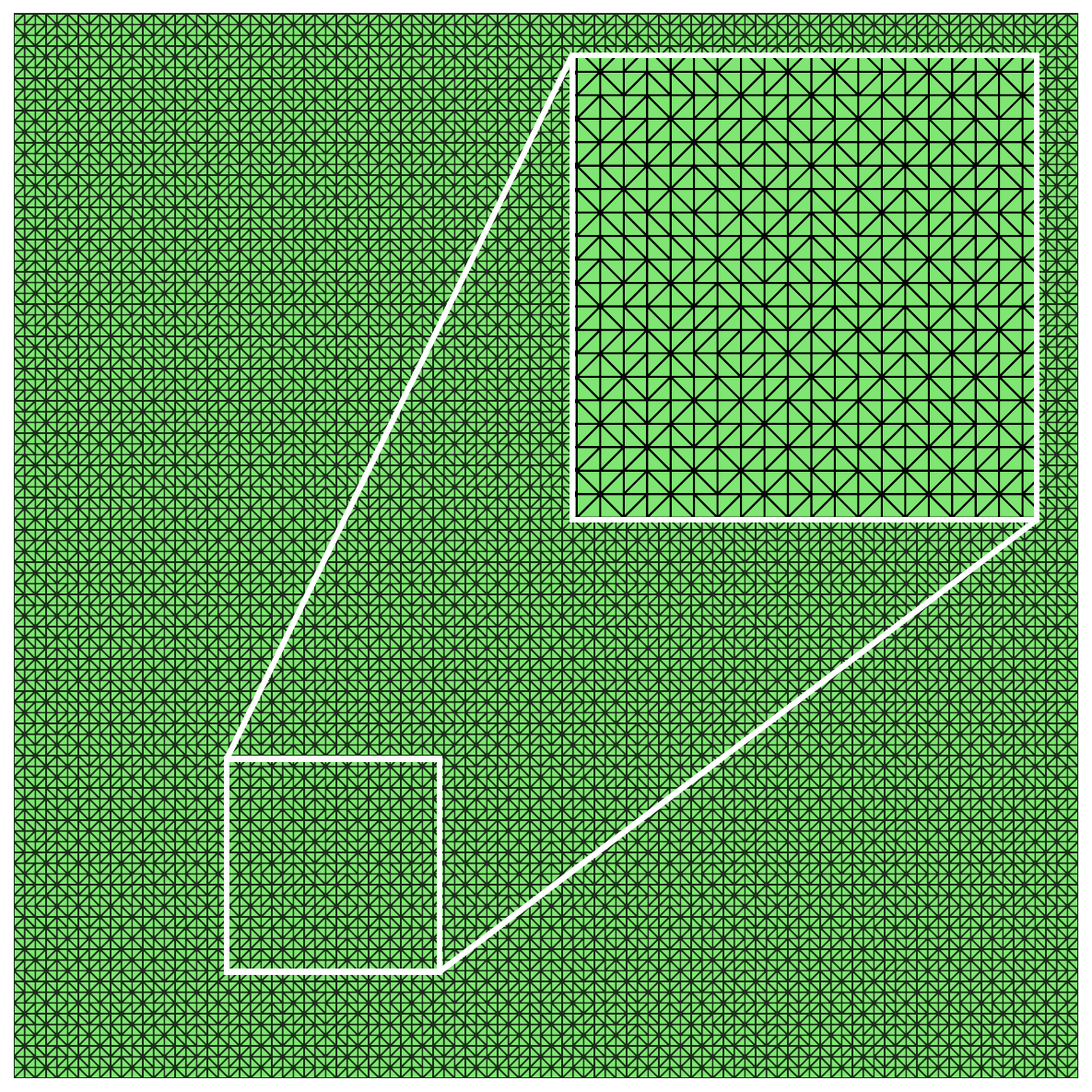}
            \caption{}
            \label{fig:ELM-test-A1}
        \end{subfigure}
        \hspace{0.03\textwidth}
        \begin{subfigure}[b]{0.68\textwidth}
            \includegraphics[width=0.92\textwidth]{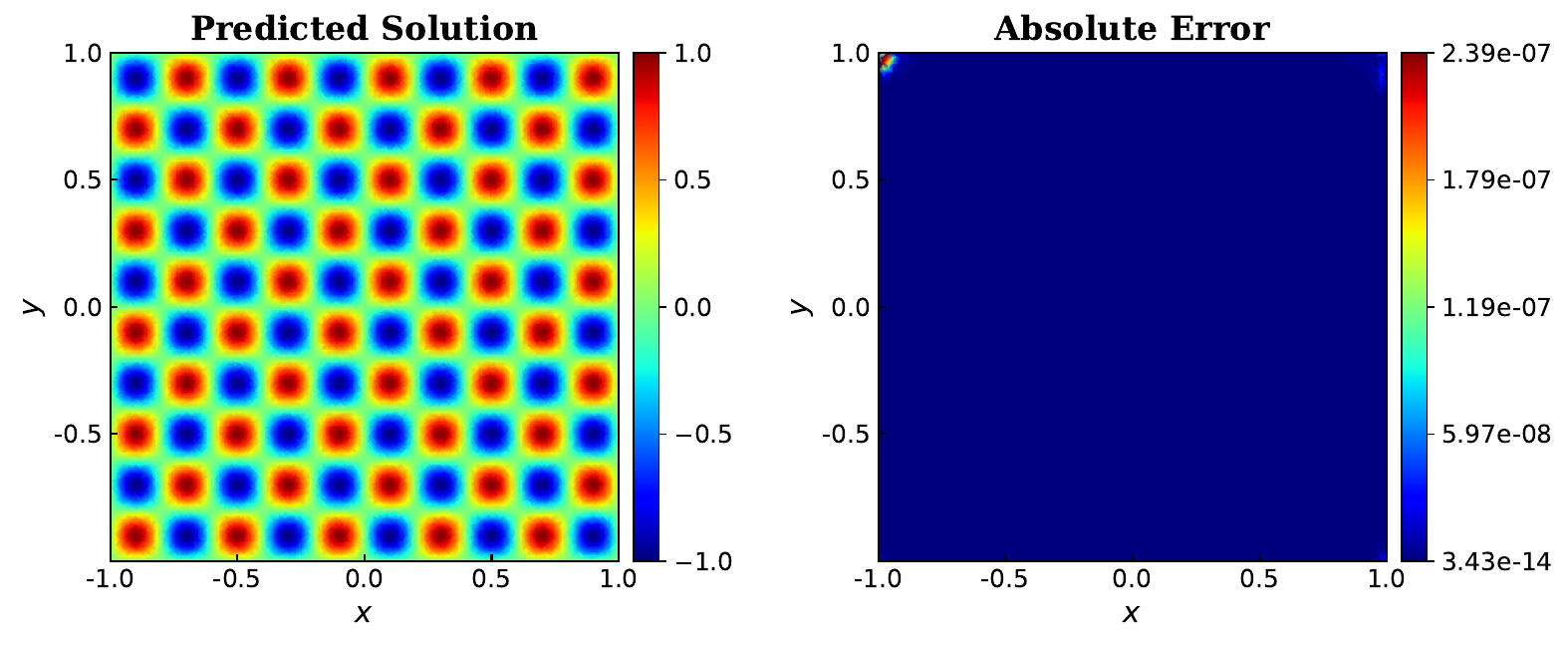}
            \caption{}
            \label{fig:ELM-test-B1}
        \end{subfigure}	
        \caption{\Cref{ex3}, ELM interpolation on a quasi-uniform mesh. 
            (A) Test mesh; (B) interpolated solution and error.}
        \label{fig:ELM-test1}
    \end{figure}

	\begin{example} \label{ex3}
		We again consider the benchmark function $u(\mathbf{x})$ introduced in \Cref{ex1}. 
		The ELM configuration is $[2,256,2500,1]$ with activation function $\sigma(x)=\sin x$. 
        The weight matrices $\{\bm{W}^{(1)}, \bm{W}^{(2)}\}$ are sampled uniformly from $[-0.4,0.4]$, while the bias vectors $\{\bm{b}^{(1)}, \bm{b}^{(2)}\}$ are initialized to zero. 
        Training is performed on all source-mesh nodes $\{\mathbf{x}_i\}_{i=1}^N$ with the corresponding function values $\{u(\mathbf{x}_i)\}_{i=1}^N$. 
	\end{example}
    \Cref{fig:ELM-test1} presents the interpolation results on a uniform triangular mesh with 10,000 nodes. 
    \Cref{fig:ELM-test-A1} shows the test mesh, while \Cref{fig:ELM-test-B1} shows the interpolated solution together with the absolute error field. 
    The relative $L^2$ error is $6.11 \times 10^{-10}$ and the mean absolute error is $2.14 \times 10^{-11}$. 
    Both error metrics are several orders of magnitude smaller than those reported for the MLP baseline in \Cref{ex1}. 
    
	It should be emphasized that the reported errors should not be interpreted as optimal, since the ELM performance still depends on several hyperparameters including the numbers of hidden units in the two hidden layers and the distribution used to initialize the random weights. 
    In particular, the initialization distribution can substantially affect the final approximation quality \cite{Wang2017}. 
    Although the reported errors should not be regarded as optimal, the present example still demonstrates that ELM is capable of delivering extremely accurate and rapid transfer on the quasi-uniform benchmark. 
	As in the previous section, the next example moves from the quasi-uniform setting to an adaptively refined source mesh. 
    If high accuracy can also be maintained in that setting, it would provide especially strong evidence for the practical value of ELM-based transfer. 

\begin{figure}[htbp]
    \centering
        \begin{subfigure}[b]{0.25\textwidth}            	
            \includegraphics[width=\textwidth]{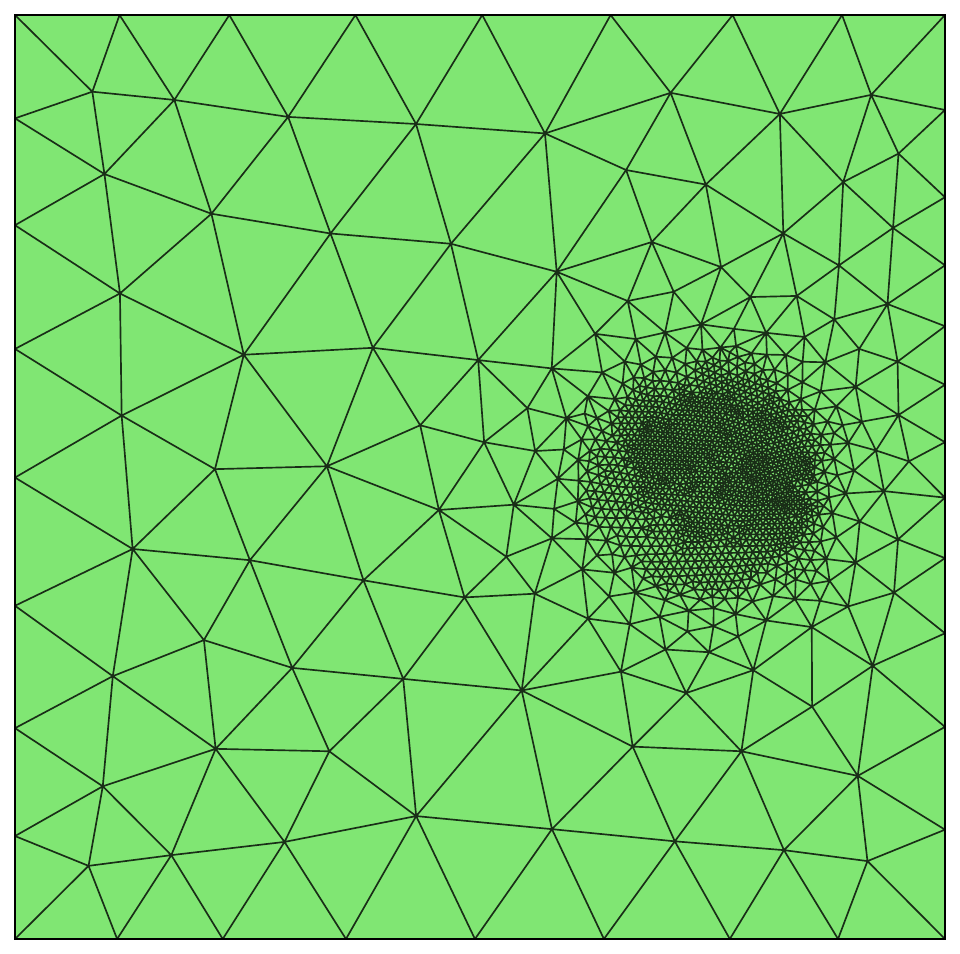}
            \caption{}
            \label{fig:ELM-test-A2}
        \end{subfigure}
        \hspace{0.03\textwidth}
        \begin{subfigure}[b]{0.68\textwidth}
            \includegraphics[width=0.92\textwidth]{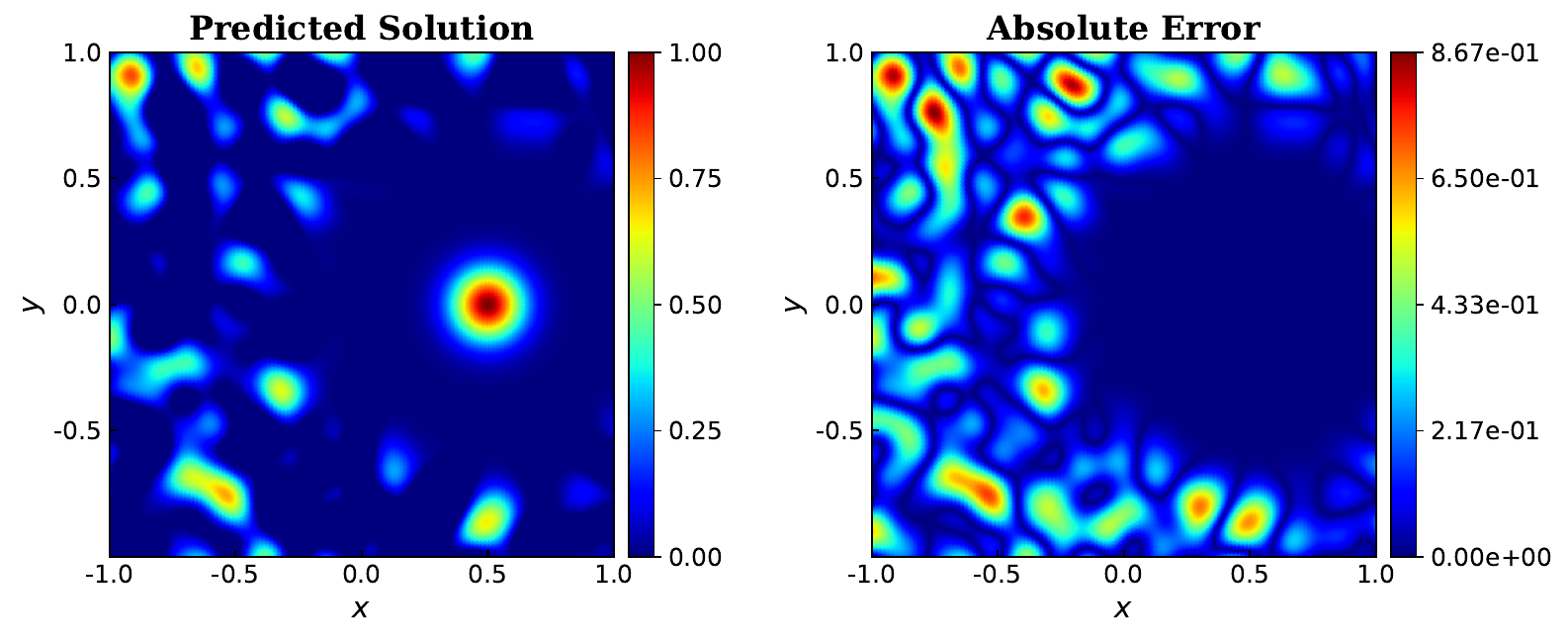}
            \caption{}
            \label{fig:ELM-test-B2}
        \end{subfigure}

        \begin{subfigure}[b]{0.25\textwidth}
            \includegraphics[width=\textwidth]{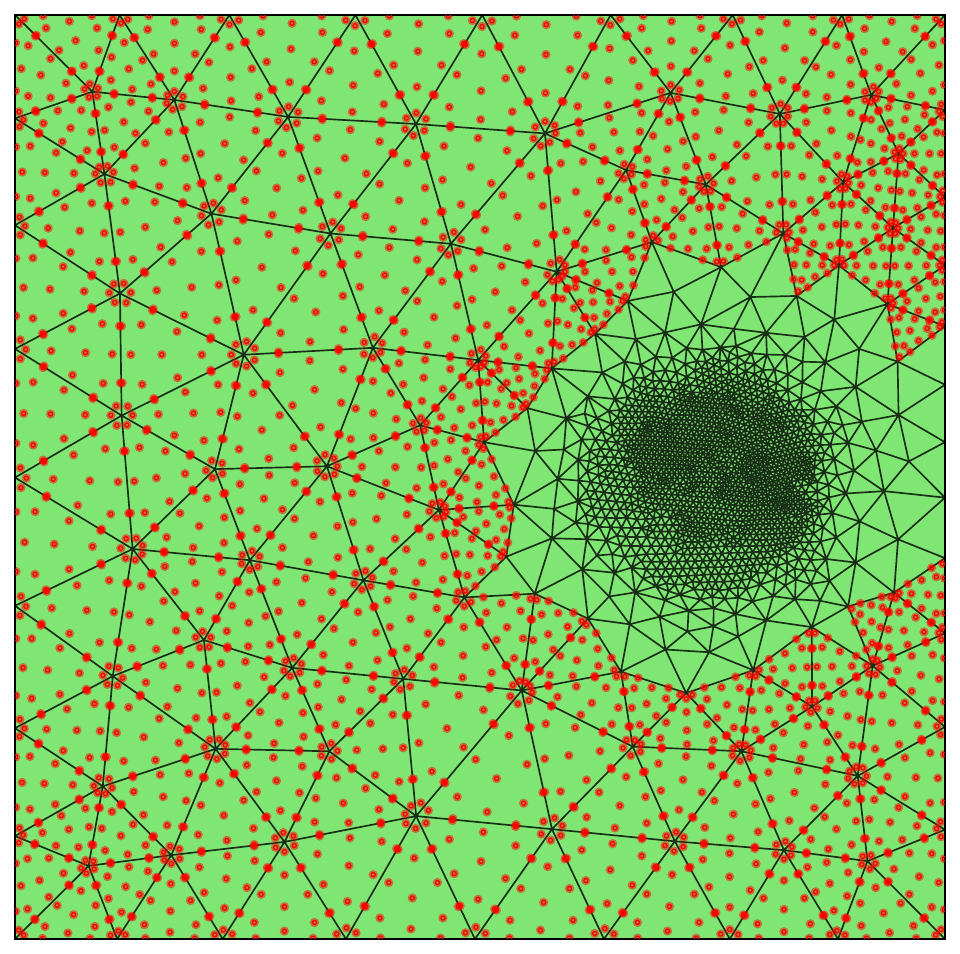}
            \caption{}
            \label{fig:ELM-test-C2}
        \end{subfigure}	
        \hspace{0.03\textwidth}
        \begin{subfigure}[b]{0.68\textwidth}
            \includegraphics[width=0.92\textwidth]{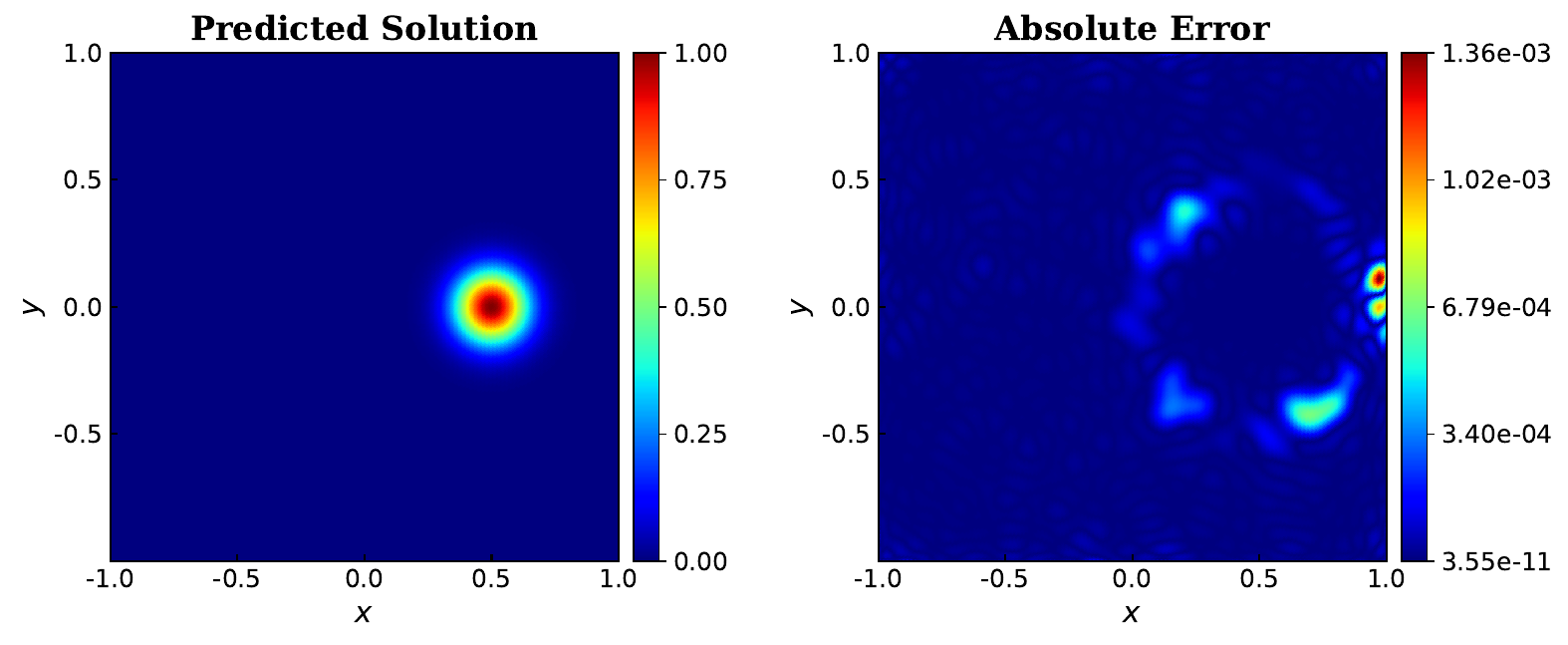}
            \caption{}
            \label{fig:ELM-test-D2}
        \end{subfigure}	
        \caption{\Cref{ex4}, ELM interpolation on an adaptive mesh. 
            (A) training mesh; 
            (B) interpolated solution and absolute error on a uniform $100 \times 100$ test grid; 
            (C) training mesh with an additional Gaussian point (red); 
            (D) interpolated solution and absolute error.}
        \label{fig:ELM-test2}
    \end{figure}

	\begin{example} \label{ex4}
        Repeating the heat-equation example introduced in \Cref{ex2}, we apply ELM to the same adaptive-mesh transfer problem. 
        The training samples remain nodal values extracted from the source mesh and are taken from a relatively sparse adaptive mesh characterized by highly nonuniform spatial resolution across the domain (\Cref{fig:ELM-test-A2}). 
        Here the ELM adopts a two-hidden-layer architecture $[2,128,1024,1]$, with activation function $\sigma(x)=\sin x$, weight matrices $\{\bm{W}^{(1)},\bm{W}^{(2)}\}$ sampled independently from the uniform distribution on $[-0.4,0.4]$, and bias vectors initialized to zero. 
        Interpolation results obtained with this setting are shown in \Cref{fig:ELM-test2}. 
	\end{example}
	
    After an extensive grid search over the main hyperparameters, the resulting interpolation still fails to achieve satisfactory accuracy, with a mean absolute error of approximately $1.19\times10^{-1}$. 
    As shown in \Cref{fig:ELM-test-B2}, the largest errors are concentrated mainly in regions where the mesh nodes are sparse. 
    When the training nodes are highly concentrated in a limited part of the domain, the model tends to fit those densely sampled regions much more effectively, while leaving sparsely sampled regions underrepresented. 
    Even though the target function remains smooth in the coarse-mesh regions, such imbalance still disrupts the recovery of the global structure and gives rise to substantial local errors. 
    These results indicate that the approximation quality of ELM is strongly influenced by the sampling distribution of the source mesh, and in this adaptive-mesh setting its accuracy is markedly worse than that of the MLP in \Cref{ex2}. 
    
    A natural remedy is therefore to enrich the training set by introducing additional Gaussian quadrature points in the finite element mesh. 
    This choice is convenient in the present setting, since such points can be generated elementwise and their values can be obtained directly by finite element interpolation. 
    \Cref{fig:ELM-test-C2} shows the enriched training mesh, where the red markers denote the added Gaussian points. 
    These points are inserted into elements whose areas exceed a prescribed threshold. 
    The corresponding interpolation result is displayed in \Cref{fig:ELM-test-D2}. 
    The region of largest error shifts away from the originally under-sampled area, and the overall accuracy improves by approximately four orders of magnitude, with a mean absolute error of $2.52\times10^{-5}$. 
    This confirms that the poor performance in the original setting is closely related to the severe nonuniformity of the sampling distribution. 
    Nevertheless, the improved accuracy is achieved only under carefully chosen conditions. 
    The improved result is obtained only after substantial manual adjustment, including the use of the network architecture $[2,128,1024,1]$, weight matrices initialized uniformly in $[-0.7,0.7]$, and the insertion of sixteen Gaussian points into each element whose area exceeds $0.007$. 
    Turning this procedure into a genuinely adaptive enrichment strategy is difficult in practice, since there is no reliable criterion for selecting elements, controlling the growth of training points, or deciding when the enrichment process should terminate. 
    The deeper issue is that such sample enrichment does not increase the expressive capacity of the network itself, so genuine adaptivity may need to be built into the model rather than introduced through external point enrichment alone. 
    Along this line, the next section moves from external sample enrichment to model-level adaptivity.

\subsection{Data Transfer via Radial Basis Function-ELM}\label{sec:RBF}
    In addition to globally supported activations, the ELM framework can be equipped with radial basis functions (RBFs) to obtain a localized hidden-layer representation.
    Owing to this locality, each input predominantly activates nearby basis units, which is often advantageous for capturing local structure in scattered data.
    In 1986, Micchelli established a foundational theoretical framework for classical RBF interpolation \cite{Micchelli1986}.
    Classical RBF interpolation places the centers at the data sites, and determines the coefficients from the interpolation conditions.
    More generally, regularized fitting variants replace the interpolation equations with a regularized fitting system.
    Accordingly, classical RBF methods are primarily concerned with multivariate interpolation and approximation rather than with a neural-network architecture.
    Broomhead and Lowe \cite{Broomhead1988} first reinterpreted the classical RBF construction as a three-layer feedforward network and showed that, once the centers are fixed, the output weights can be determined by a linear least-squares procedure.
    Moody and Darken \cite{Moody1989} subsequently developed this viewpoint into a practical single-hidden-layer RBF network with locally responsive hidden units.
    In this setting, only a finite number of radial basis units is used, rather than assigning one center to each data site; the centers and widths are treated as design or training parameters, while the output layer remains linear.
    Poggio and Girosi \cite{Girosi1990} later provided further theoretical support for the approximation capability of RBF networks for continuous functions.
    In this section, RBF-ELM refers specifically to the neural-network form of RBF whose training follows the ELM framework, as opposed to classical RBF. 
	The RBF-ELM feature map is defined as follows. 

    Given an input point $\mathbf{x}\in\mathbb{R}^{d_0}$, the RBF-ELM constructs its hidden features through two hidden stages followed by a linear output mapping:
    \begin{align*}
        \bm{\Phi}^{(1)}(\mathbf{x})
        &= \sigma_1\!\left(\bm{W}^{(1)}\mathbf{x}+\bm{b}^{(1)}\right), \\
        \bm{\Phi}^{(2)}(\mathbf{x})
        &= \sigma_2\!\left(\bm{W}^{(2)}\bm{\Phi}^{(1)}(\mathbf{x})\right), \\
        \tilde{u}_{\theta}(\mathbf{x})
        &= \bm{\beta}^{\top}\bm{\Phi}^{(2)}(\mathbf{x}),
    \end{align*}
    where $\bm{W}^{(1)} \in \mathbb{R}^{d_0 d_1 \times d_0}$ and $\bm{W}^{(2)} \in \mathbb{R}^{d_1 \times d_0 d_1}$ denote the fixed weight matrices of the first and second hidden stages, and $\bm{b}^{(1)} \in \mathbb{R}^{d_0 d_1}$ denotes the fixed bias vector in the first hidden stage. 
    Nonlinearity in the hidden stages is introduced through the element-wise activation functions $\sigma_1(\cdot)$ and $\sigma_2(\cdot)$, defined by
    \[
    \sigma_1(t)=t^2, \qquad \sigma_2(t)=e^{-t}.
    \]
    With appropriate choices of the fixed parameters $\bm{W}^{(1)}$, $\bm{W}^{(2)}$, and $\bm{b}^{(1)}$, the $j$th component of $\bm{\Phi}^{(2)}(\mathbf{x})$ can be written in the Gaussian form
    \begin{equation}\label{eq:rbf_gaussian_feature}
        \phi^{(2)}_j(\mathbf{x})
        =
        \exp\!\left(-\frac{\|\mathbf{x}-\mathbf{c}_j\|^2}{\varepsilon_j^2}\right),
        \qquad j=1,\dots,d_1. 
    \end{equation}
    Here, $\mathbf{c}_j\in\mathbb{R}^{d_0}$ is the center of the $j$th Gaussian basis function, and $\varepsilon_j>0$ is the width parameter controlling the spatial spread of its response. 
    In practice, the centers may be chosen by uniform placement, random sampling, or clustering-based initialization (e.g., $k$-means), whereas the width parameters are commonly fixed either as a uniform constant or according to the distances between neighboring centers. 
    Each Gaussian kernel yields a localized response around its center, analogous in this respect to a finite element basis function, although Gaussian RBFs are not compactly supported in general. 
    At the output layer, the coefficient vector $\bm{\beta}\in\mathbb{R}^{d_1}$ forms the final linear combination of the hidden features $\bm{\Phi}^{(2)}(\mathbf{x})$. 
    \begin{figure}[htbp]
        \centering
        \includegraphics[width=0.7\textwidth]{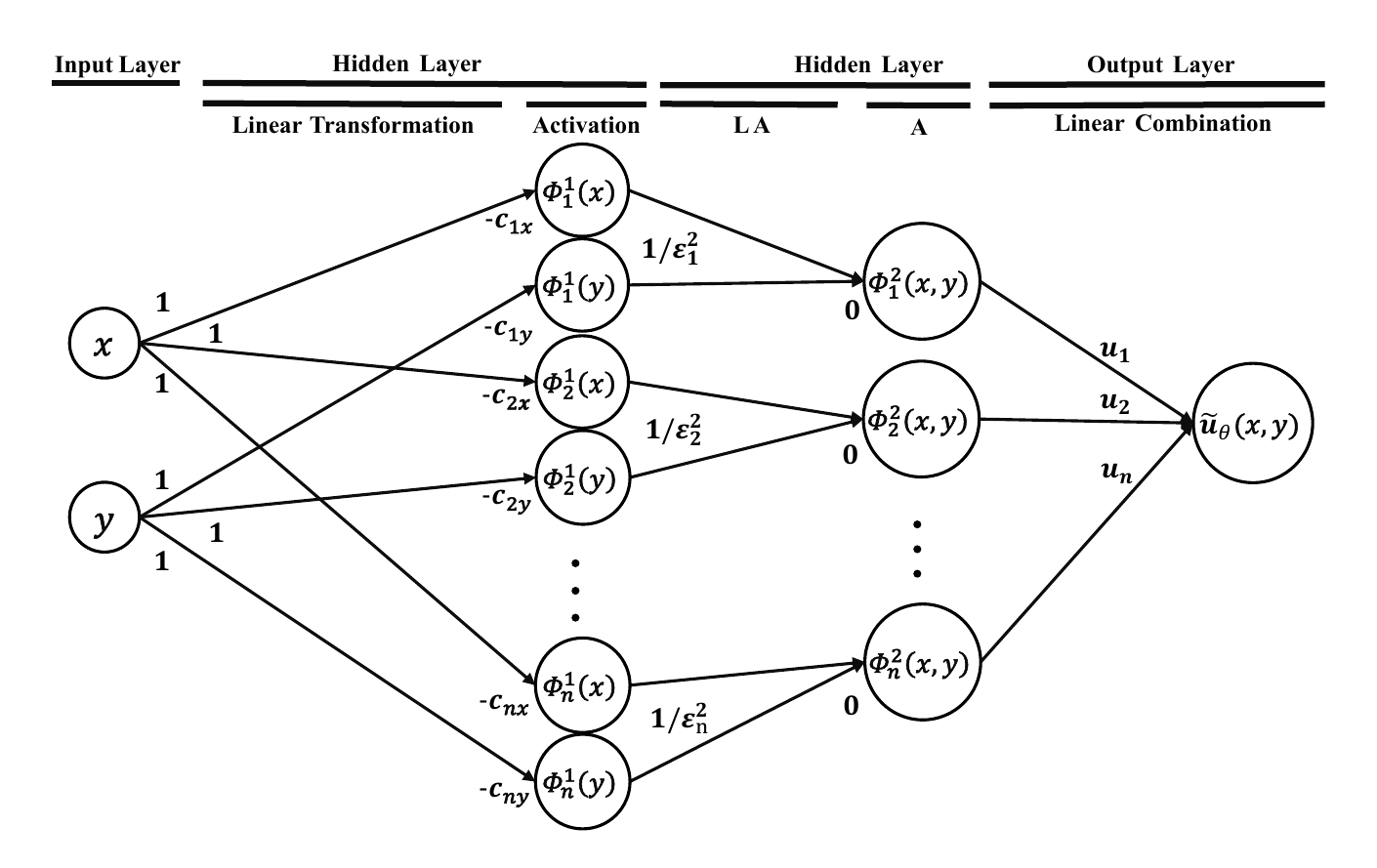}
        \caption{Schematic architecture of a RBF-ELM network.}
        \label{fig:RBF}
    \end{figure}
    \Cref{fig:RBF} shows a two-dimensional realization with input $\mathbf{x}=(x,y)^\top$. 
    Accordingly, the first hidden stage has width $2n$, whereas the second hidden stage has width $n$. 
    The resulting architecture is non-fully-connected, since the $j$th unit in the second hidden stage is connected only to the two first-stage units associated with the $j$th basis function. 
    The matrix $\bm{W}^{(1)}\in\mathbb{R}^{2n\times 2}$ is formed by stacking $n$ copies of the identity matrix $\bm{I}_2$,
    \[
    \bm{W}^{(1)}=\mathbf{1}_{n}\otimes \bm{I}_2,
    \]
    and the corresponding bias vector is
    \[
    \bm{b}^{(1)}
    =
    (-c_{1x},-c_{1y},-c_{2x},-c_{2y},\dots,-c_{nx},-c_{ny})^\top
    \in\mathbb{R}^{2n},
    \]
    where $\mathbf{1}_{n}\in\mathbb{R}^{n}$ denotes the all-ones vector, $\bm{I}_2\in\mathbb{R}^{2\times 2}$ is the $2\times 2$ identity matrix, and $\otimes$ denotes the Kronecker product. 
    At the second hidden stage, $\bm{W}^{(2)}\in\mathbb{R}^{n\times 2n}$ takes the form
    \[
    \bm{W}^{(2)}
    =
    \operatorname{diag}\!\bigl(\varepsilon_1^{-2},\dots,\varepsilon_{n}^{-2}\bigr)\otimes
    \begin{bmatrix}
    1 & 1
    \end{bmatrix},
    \]
    where $\operatorname{diag}(\cdot)$ denotes a diagonal matrix. 
    These choices yield the coordinate representation of \eqref{eq:rbf_gaussian_feature}:
    \[
    \phi_j^{(2)}(\mathbf{x})
    =
    \exp\!\left(
    -\frac{(x-c_{jx})^2+(y-c_{jy})^2}{\varepsilon_j^2}
    \right),
    \qquad j=1,\dots,n.
    \]
    Let $\bm{\Phi}(\mathbf{x})$ denote the hidden feature vector produced by the two hidden stages. 
    The RBF-ELM approximation can then be written in the same linear-output form as ELM,
    \[
    \tilde{u}_{\theta}(\mathbf{x})=\bm{\beta}^{\top}\bm{\Phi}(\mathbf{x}).
    \]
    Its main structural difference from ELM lies in the construction of the hidden basis. 
    This distinction motivates a return to the adaptive-mesh data-transfer example, where standard ELM failed to achieve satisfactory accuracy. 
    
\begin{figure}[ht!]
    \centering
        \includegraphics[width=0.98\textwidth]{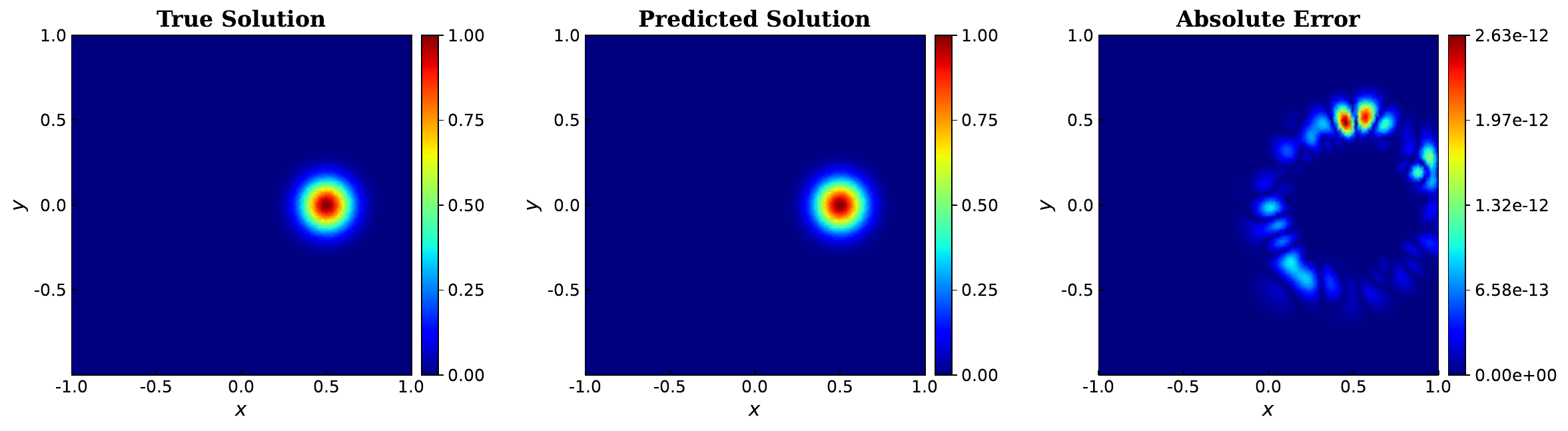}
        \caption{\Cref{ex5}, RBF-ELM Interpolation on adaptive mesh and corresponding absolute error distribution.}
        \label{fig:RBF-test1}
    \end{figure}

	\begin{example} \label{ex5}
		Repeating the heat-equation example introduced in \Cref{ex2}, we apply RBF-ELM to the same adaptive-mesh transfer problem.
        The training samples consist of all nodal coordinates from the adaptive mesh, with the corresponding finite element solution values taken as labels.
        Here the RBF-ELM adopts the architecture $[2,300,1]$, where the $300$ centers are randomly selected from the training nodes and the width parameters satisfy $\varepsilon_1^2=\cdots=\varepsilon_{300}^2=1/60$. 
	\end{example}
    \Cref{fig:RBF-test1} shows the interpolation results obtained with this setting. 
    Compared with the standard ELM, the RBF-ELM reconstruction is markedly more accurate and reproduces the solution profile much more faithfully over the entire domain. 
    The resulting errors are negligible for the present computation, with a mean absolute error of $4.09\times 10^{-14}$ and a relative $\ell^2$ error of $1.97\times 10^{-12}$.
    This improvement comes from the fact that the two models organize their basis interactions in very different ways. 
    In the ELM, the hidden features are globally entangled through the network structure, are strongly coupled across the whole domain, and the resulting coefficient matrix is globally dense. 
    By contrast, the hidden layer of RBF-ELM consists of Gaussian basis functions centered at selected nodes, each of which responds primarily to its local neighborhood. 
    This is similar to the local-support mechanism of finite element basis functions: a basis function interacts mainly with its neighboring ones, rather than with the entire approximation space. 
    As a result, the corresponding coefficient matrix in RBF-ELM exhibits weaker long-range correlation (i.e., entries near the diagonal are significant, while those far from the diagonal are nearly zero). 
    This does not, however, imply a noticeable improvement in conditioning, since the condition number of the ELM feature matrix is \(2.86\times 10^{16}\), while that of the RBF-ELM feature matrix is \(5.63\times 10^{16}\). 
    The gain in accuracy should therefore be understood primarily as a consequence of the localized representation itself. 
    The approximation at any point is dominated by nearby basis functions, while contributions from distant centers are negligible. 
    Such a mechanism preserves local shape information more faithfully and suppresses the spurious nonlocal effects induced by globally coupled hidden representations. 
    
    \begin{remark}\label{rbfAnalysis}
        In the reported RBF-ELM experiments, the centers are randomly selected from the training nodes and all basis functions share a common width parameter. 
        As noted earlier, in a fully trainable neural-network formulation of RBF, these hidden-layer parameters may also be optimized. 
        The decision not to train these parameters in this work is based on a practical observation. 
        We conducted extensive experiments over a broad range of initial parameter choices, allowing the width parameters to vary across centers during training. 
        \Cref{fig:error-analysis-rbf} illustrates the test results on the target mesh \(\mathcal{T}_B\), with the three panels from left to right showing the error before training, the error after training, and their ratio (on a logarithmic scale). 
        The horizontal axis represents the number of centers \(N_c\), and the vertical axis characterizes the initial width parameter in the form of \(\gamma = 1/\varepsilon^2\). 
        It can be observed that the practical effect of training is unpredictable: it may help, have no effect, or even degrade performance. 
        Even when improvement occurs, the gain in accuracy is limited. 
        It can also be observed that a large number of initial parameter combinations already yield results approaching an accuracy of \(10^{-12}\), and further optimization often leads to negative effects rather than improvement. 
        For this reason, fixing the parameters without further training is a viable strategy. 
        This paper reports results obtained with empirically chosen parameters, which are likely improvable with more systematic tuning. 
        \begin{figure}[ht!]
        \centering
            \includegraphics[width=\textwidth]{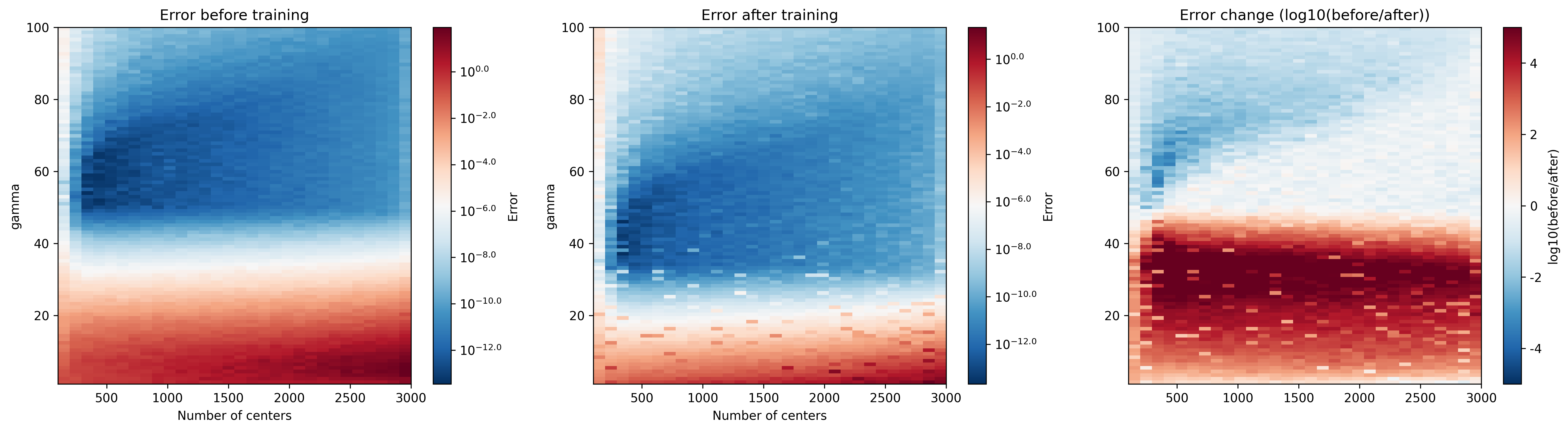}
            \caption{RBF-ELM generalization error: initial error (left), post-training error (middle), and their log ratio (right)}
            \label{fig:error-analysis-rbf}
        \end{figure}
    \end{remark}
    
    Having examined the three interpolation models individually, we now compare their respective strengths and limitations. 
    Table~\ref{tab:comparison} summarizes the main characteristics of MLP, ELM and RBF-ELM for mesh-based data transfer. 
    MLP is applicable to arbitrary meshes and is relatively insensitive to sampling irregularity, since its parameters are optimized through iterative training over the whole dataset. 
    This global optimization makes the model robust, but it also leads to a relatively high training cost and, in the present experiments, the attainable accuracy typically saturates around a relative $L_2$ error of $10^{-3}$. 
    ELM, by contrast, is extremely efficient because the output weights are obtained from a single linear solve once the hidden-layer parameters are fixed. 
    However, this efficiency is accompanied by a strong dependence on the hidden representation: under uniform sampling, ELM can achieve near machine precision, whereas under irregular sampling its globally coupled hidden features no longer provide a reliable transfer mechanism and the accuracy deteriorates rapidly. 
    RBF-ELM inherits the low training cost of ELM while replacing the globally mixed hidden representation by a localized radial-basis layer. 
    As a result, the model remains efficient, yet is much less sensitive to mesh irregularity because nearby basis functions mainly determine nearby values.  
    Although its performance still depends on the choice of the number of centers and the width parameter, this dependence is considerably milder than the sensitivity of ELM to initialization and hidden-feature construction. 
    Taken together, these observations explain why RBF-ELM offers the most balanced compromise among applicability, efficiency, accuracy and stability for the transfer problems considered here. 

    \begin{table}[!ht]
    \centering
    \caption{Comparison of interpolation models for mesh-based data transfer.}
    \label{tab:comparison}
    \renewcommand{\arraystretch}{1.35}
    \begin{NiceTabular}{
        >{\raggedright\arraybackslash}m{3.7cm}
        >{\centering\arraybackslash}m{2.8cm}
        >{\centering\arraybackslash}m{2.8cm}
        >{\centering\arraybackslash}m{3.0cm}
        }
        \toprule
         & \textbf{MLP} & \textbf{ELM} & \textbf{RBF-ELM} \\
        \midrule
        \textbf{Applicability}
        & \texttt{+++} & \texttt{+} & \texttt{+++} \\
        \midrule
        \textbf{Training efficiency}
        & \texttt{+} & \texttt{+++} & \texttt{+++} \\
        \midrule
        \textbf{Accuracy}
        & \texttt{+} & \texttt{++} & \texttt{+++} \\
        \midrule
        \textbf{Stability}
        & \texttt{+++} & \texttt{+} & \texttt{++} \\
        \bottomrule
        \Block[l]{1-4}{\footnotesize More ``\texttt{+}'' indicates better performance.}
    \end{NiceTabular}
    \end{table}

\section{Application}\label{sec:APPLIC}
		
    The favorable performance of RBF-ELM in the preceding tests motivates its examination in a more representative application setting. 
    Repeated data transfer between meshes of different resolutions or topologies is common in adaptive and multiphysics simulations. 
    In this setting, controlling error accumulation become particularly important. 
    Without loss of generality, we consider interpolation between two non-nested meshes, where the field data are alternately transferred back and forth. 

    Let \(u\) denote the reference solution. 
    The initial nodal datasets on the two meshes are defined as
	\[
	\mathcal{D}_A^{(0)}=\{(\mathbf{x},u(\mathbf{x})):\mathbf{x}\in\mathcal{N}_A\},\qquad
	\mathcal{D}_B^{(0)}=\{(\mathbf{x},u(\mathbf{x})):\mathbf{x}\in\mathcal{N}_B\}.
	\]
    Starting from \(\mathcal{D}_A^{(0)}\), one transfer cycle for \(k=0,1,2,\ldots\) consists of the following two steps: 
	\begin{enumerate}
		\item[(i)] \textbf{\(\mathcal{T}_A \to \mathcal{T}_B\)}. Train \(\tilde u_{\theta_A^{(k)}}\) on \(\mathcal{D}_A^{(k)}\), then evaluate it at the node set \(\mathcal{N}_B\) to obtain
		\[
		\mathcal{D}_B^{(k+1)}=\{(\mathbf{x},\,\tilde u_{\theta_A^{(k)}}(\mathbf{x})):\mathbf{x}\in\mathcal{N}_B\},
		\]
		\item[(ii)] \textbf{\(\mathcal{T}_B \to \mathcal{T}_A\)}. Train \(\tilde u_{\theta_B^{(k+1)}}\) on \(\mathcal{D}_B^{(k+1)}\), then evaluate it at \(\mathcal{N}_A\) to obtain
		\[
		\mathcal{D}_A^{(k+1)}=\{(\mathbf{x},\,\tilde u_{\theta_B^{(k+1)}}(\mathbf{x})):\mathbf{x}\in\mathcal{N}_A\}.
		\]
	\end{enumerate}
	Repeating steps (i) and (ii) produces the alternating sequence
	\[
	\mathcal{D}_A^{(0)} \longrightarrow \mathcal{D}_B^{(1)} \longrightarrow \mathcal{D}_A^{(1)} \longrightarrow \cdots \longrightarrow \mathcal{D}_B^{(k)} \longrightarrow \mathcal{D}_A^{(k)}.
	\]
    This iterative transfer process acts as a stress test, amplifying the intrinsic behavior of data transfer methods under repeated application. 
    Taking standard finite element interpolation as the baseline, a series of numerical tests is carried out below to examine the robustness and error propagation behavior of the RBF-ELM method. 
    The test cases are arranged in increasing order of complexity, ranging from one-dimensional smooth fields to two-dimensional oscillatory and singular fields. 
    Based on the discussion in \Cref{rbfAnalysis}, a single width parameter \(\varepsilon\) is used for all centers, with its value chosen separately for different examples. 
    The mesh configuration and the corresponding RBF-ELM parameters are detailed in what follows.
    We begin with a simple one-dimensional test case.
    
\begin{figure}[htbp]
    \centering
		\begin{subfigure}[b]{0.245\textwidth}            
			\includegraphics[width=\textwidth]{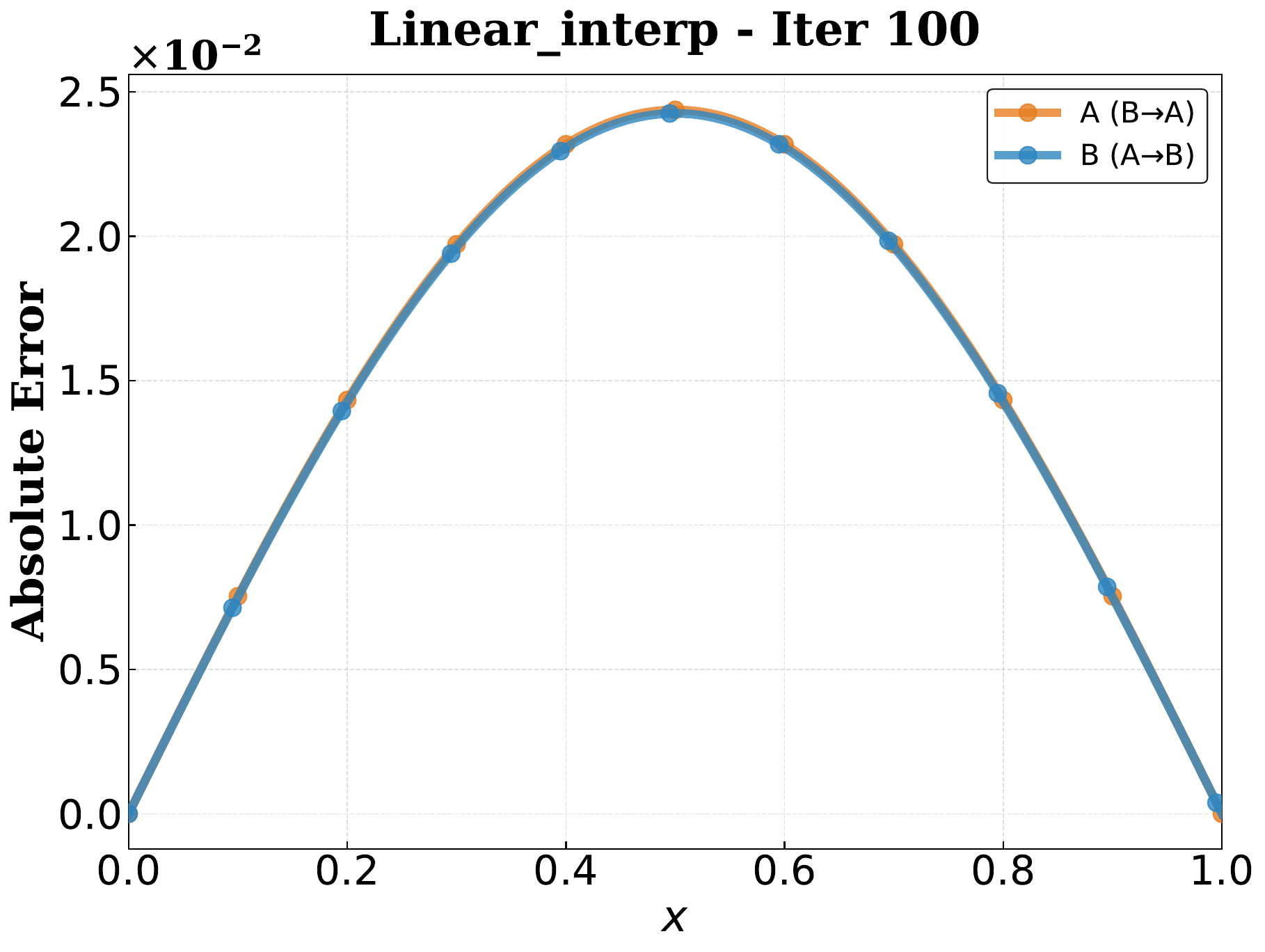}
			\caption{}
			\label{fig:APP1Dsmooth-A}
		\end{subfigure}
		\begin{subfigure}[b]{0.245\textwidth}
			\includegraphics[width=\textwidth]{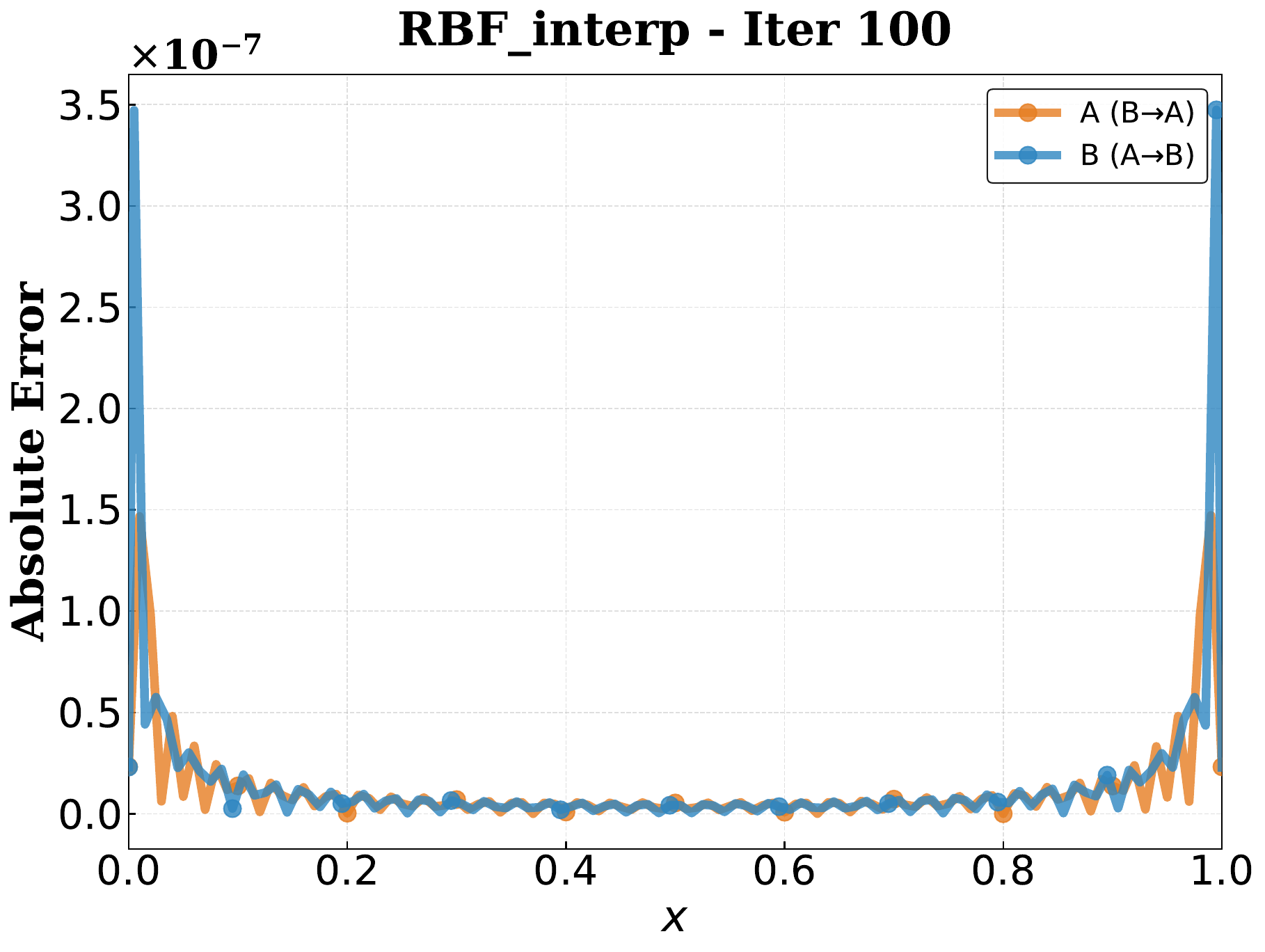}
			\caption{}
			\label{fig:APP1Dsmooth-B}
		\end{subfigure}
        \begin{subfigure}[b]{0.49\textwidth}
			\includegraphics[width=0.98\textwidth]{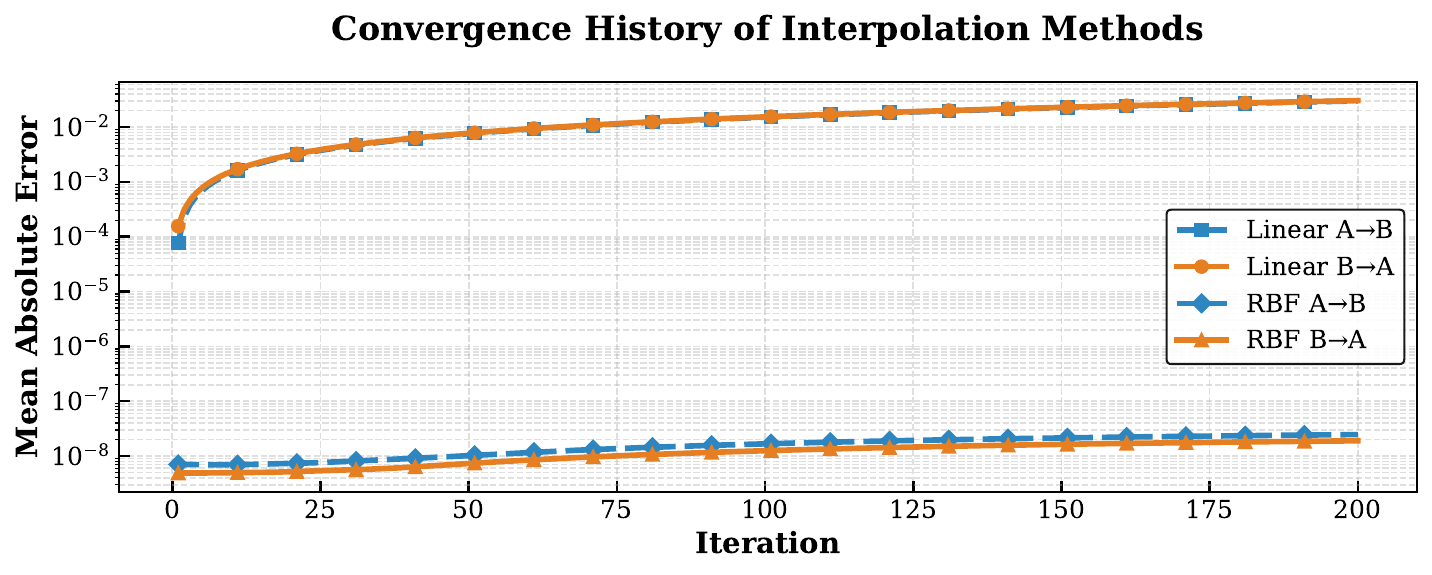}
			\caption{}
			\label{fig:APP1Dsmooth-C}
		\end{subfigure}
		\caption{\Cref{exm1d1}
			Comparison of piecewise linear interpolation and the RBF-ELM interpolation model for a smooth 1D function.
            (A) Error of piecewise linear interpolation after 100 iterations; 
            (B) Error of RBF-ELM interpolation after 100 iterations;
            (C) Error evolution history over transfer iterations.}
        \label{fig:APP1Dsmooth}
	\end{figure}

    \begin{example}\label{exm1d1}
        Consider the smooth one-dimensional field
        \[
        u(x)=\sin(\pi x), \qquad x\in[0,1].
        \]
        The source mesh \(\mathcal{T}_A\) is taken to be a uniform partition of \([0,1]\) with \(100\) nodes. 
        The target mesh \(\mathcal{T}_B\) consists of the midpoints of all elements of \(\mathcal{T}_A\) together with the boundary nodes \(x=0\) and \(x=1\). 
        In the RBF-ELM approximation, \(N_c=80\) centers are randomly selected from the training nodes, and the common width parameter is chosen as \(\varepsilon^2=1/100\). 
    \end{example}
    \Cref{fig:APP1Dsmooth} demonstrates the pointwise absolute errors for bidirectional transfer, where the blue curve shows the error from \( \mathcal{T}_A \) to \( \mathcal{T}_B \) and the orange curve corresponds to the reverse transfer. 
    A well-known deficiency of standard finite element interpolation in repeated-transfer settings is the gradual loss of solution variation over successive transfers. 
    The impact of this progressive smoothing on solution accuracy is further illustrated in \Cref{fig:APP1Dsmooth-A}. 
    After 100 alternating transfers, the maximum error of the standard finite element interpolation is on the order of \(10^{-2}\), while that of the RBF-ELM method is on the order of \(10^{-7}\) (\Cref{fig:APP1Dsmooth-B}). 
    The long-term effect of repeated transfer is shown more clearly in \Cref{fig:APP1Dsmooth-C}, where the error is plotted against the number of transfer cycles. 
    Although both methods exhibit error growth, the increase is much slower for RBF-ELM, indicating that it is considerably more resistant to error accumulation under repeated interpolation. 

\begin{figure}[htbp]
    \centering
		\begin{subfigure}[b]{0.49\textwidth}           
			\includegraphics[width=0.49\textwidth]{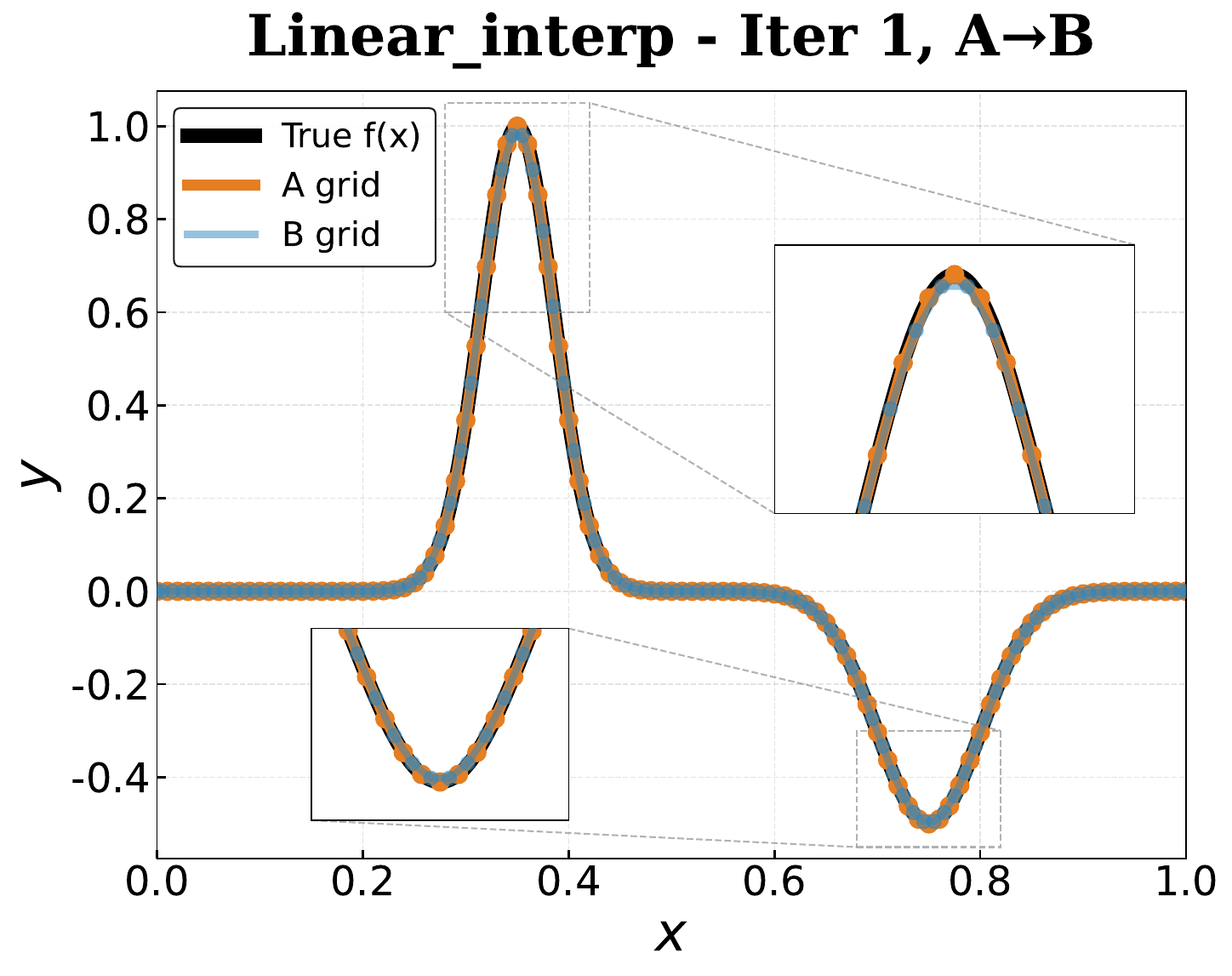}
			\includegraphics[width=0.49\textwidth]{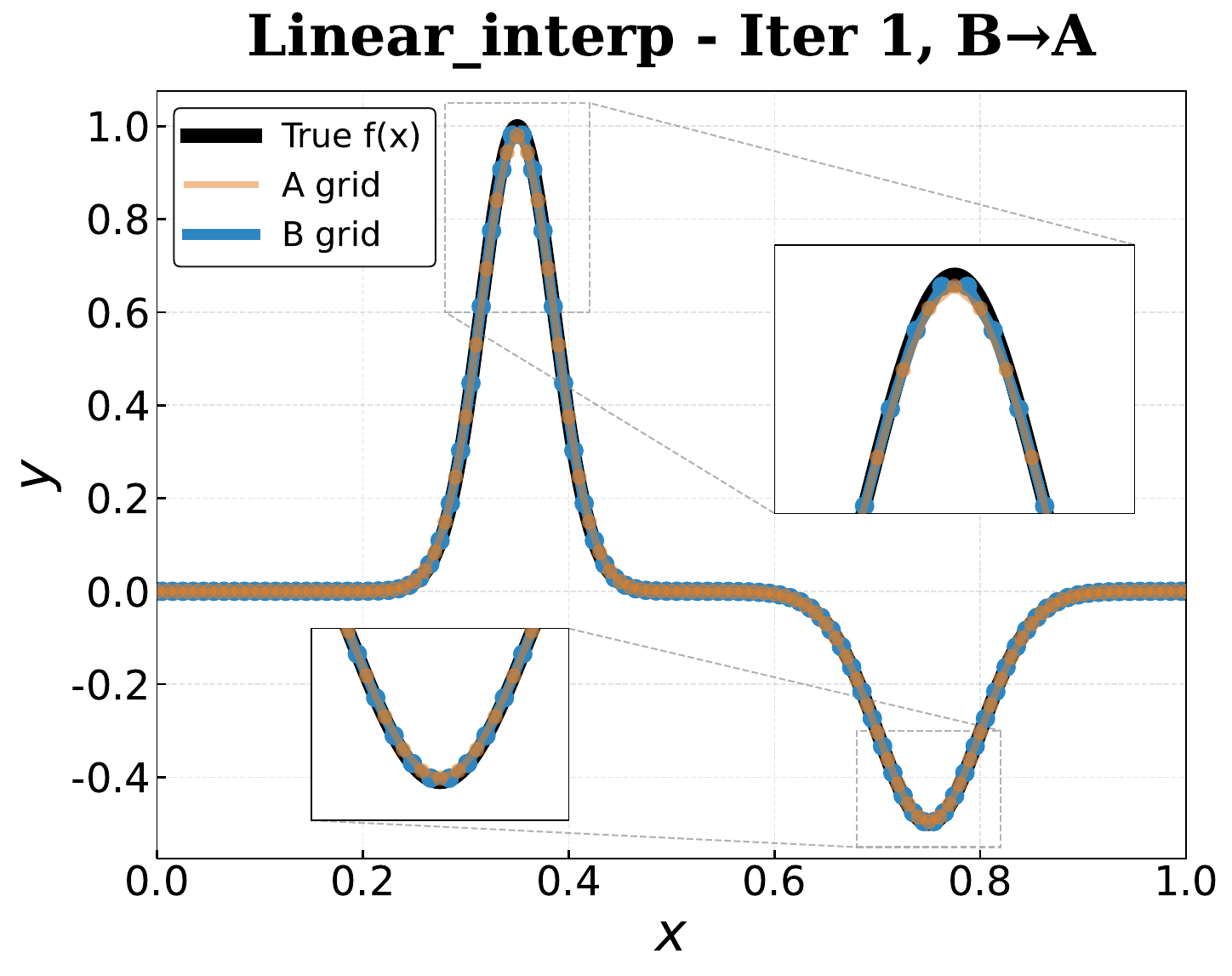}
			\caption{}
			\label{fig:APP1DpeakL-A}
		\end{subfigure}
		\begin{subfigure}[b]{0.49\textwidth}
			\includegraphics[width=0.49\textwidth]{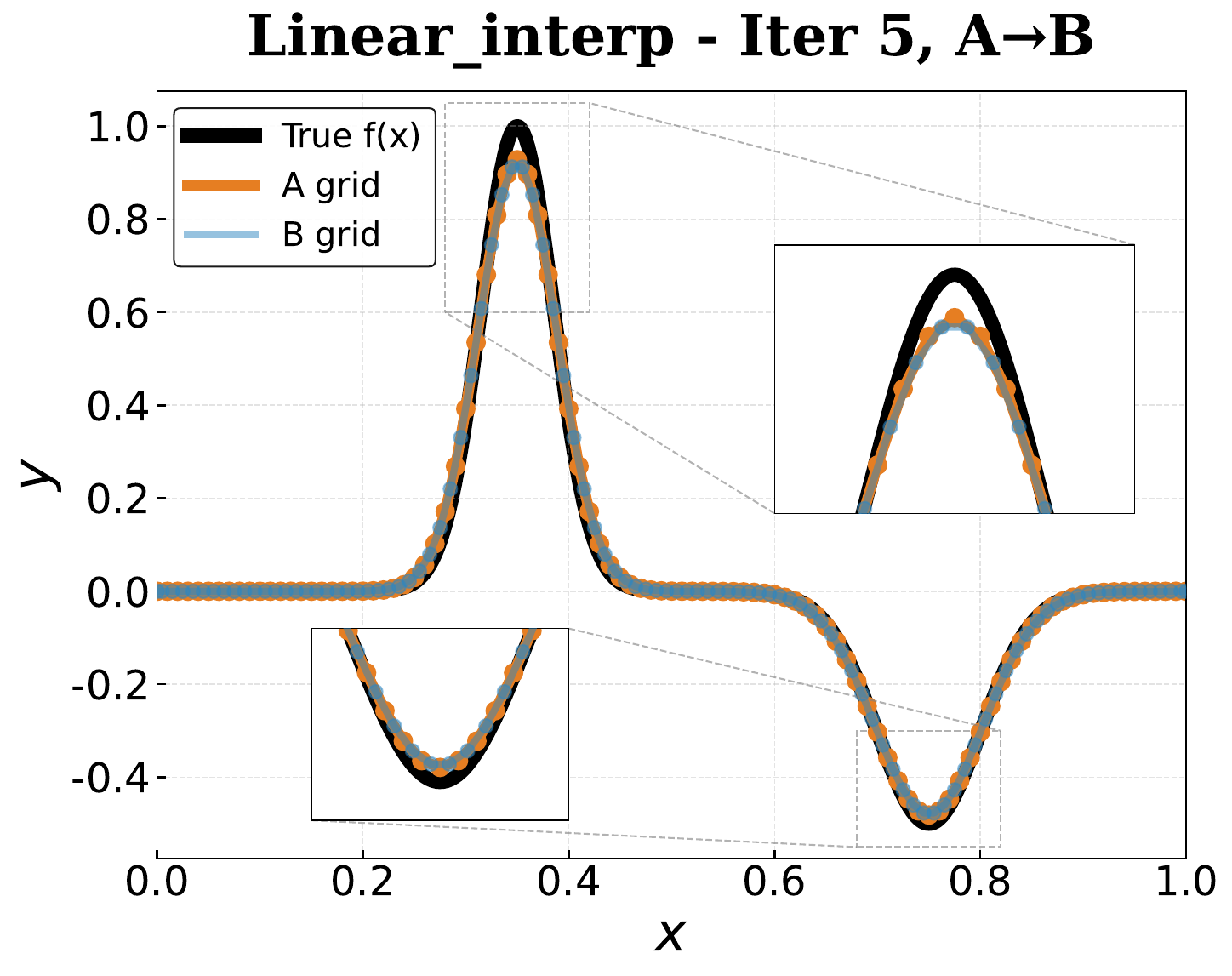}
			\includegraphics[width=0.49\textwidth]{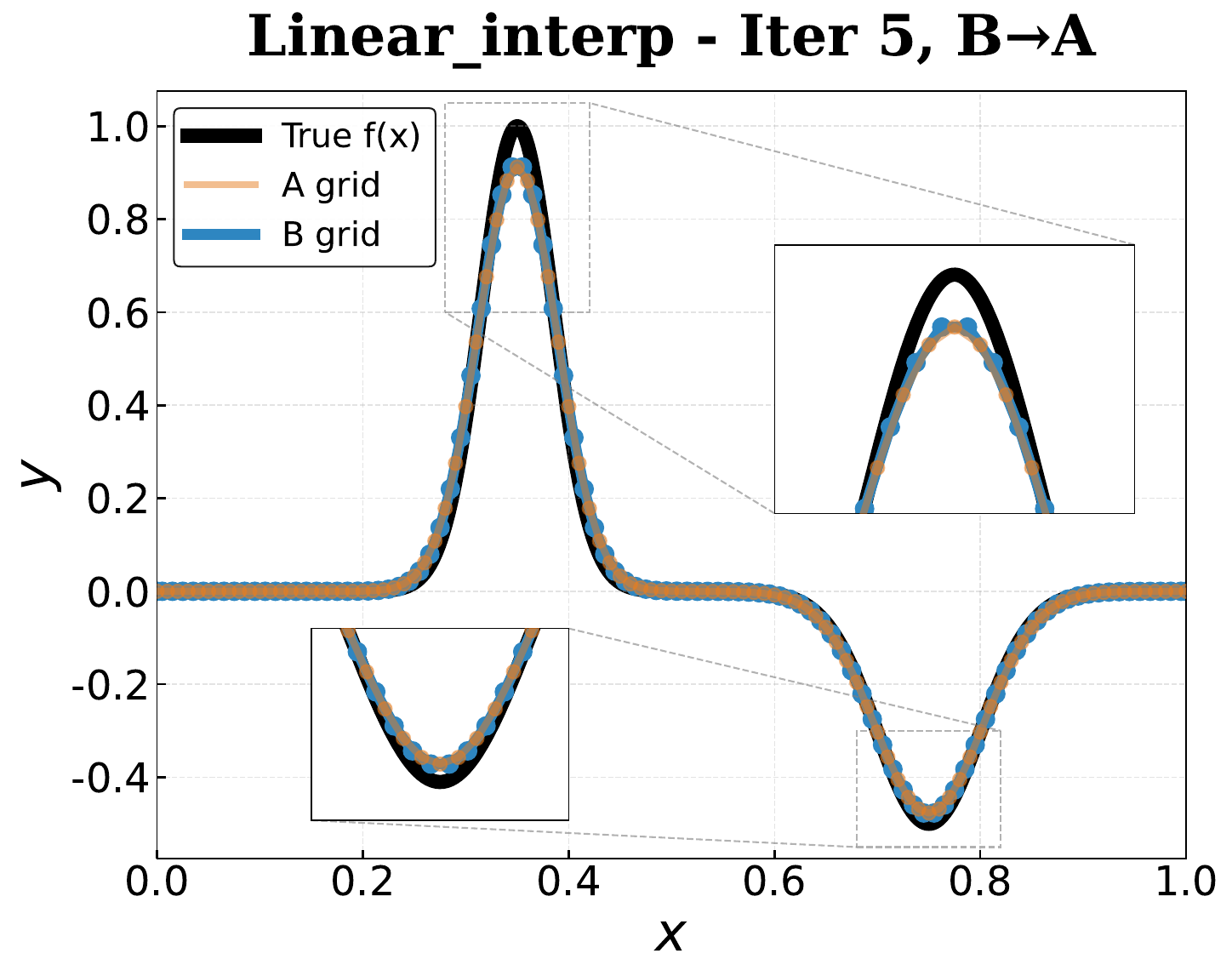}
			\caption{}
			\label{fig:APP1DpeakL-B}
		\end{subfigure}	
		\begin{subfigure}[b]{0.49\textwidth}
			\includegraphics[width=0.49\textwidth]{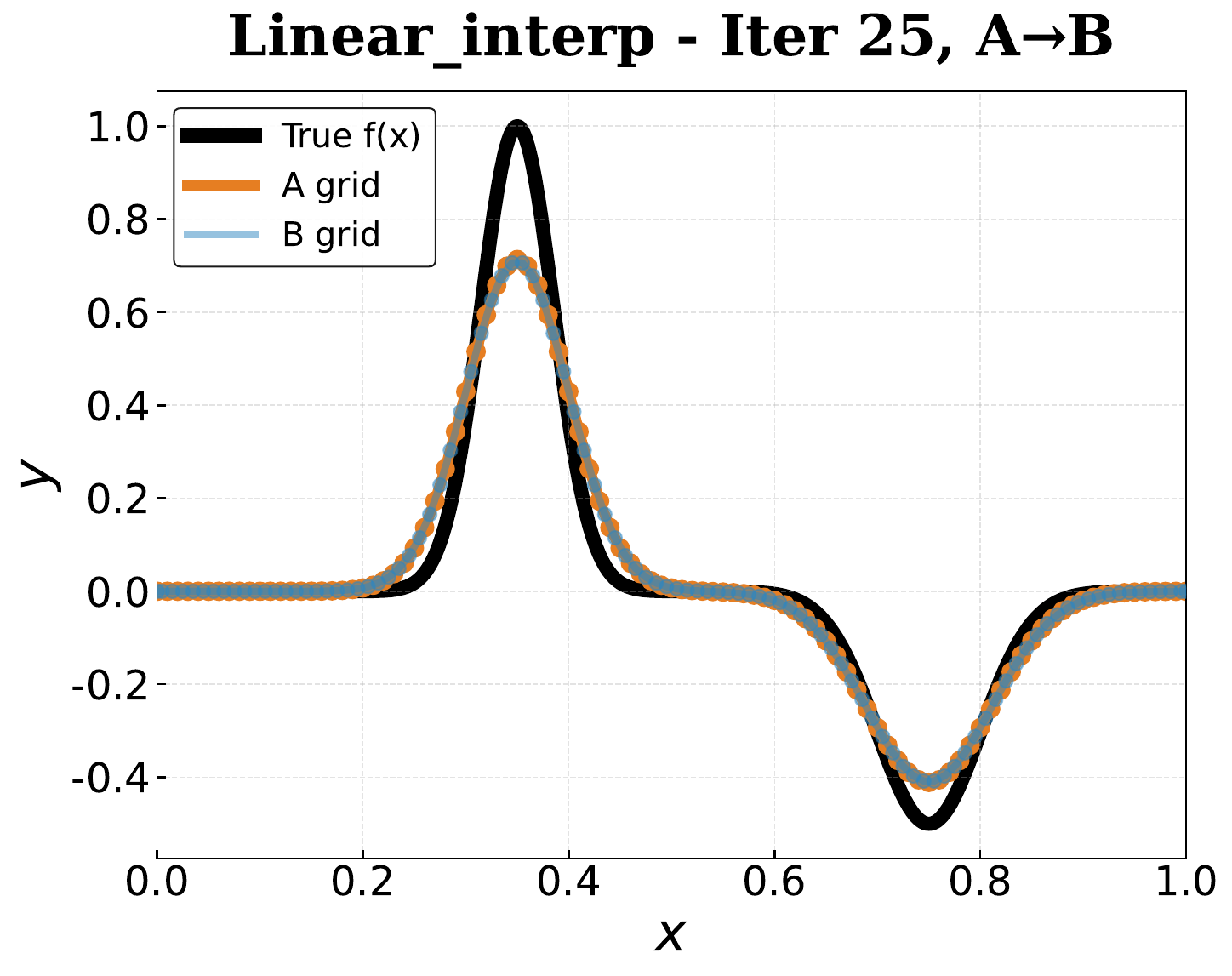}
			\includegraphics[width=0.49\textwidth]{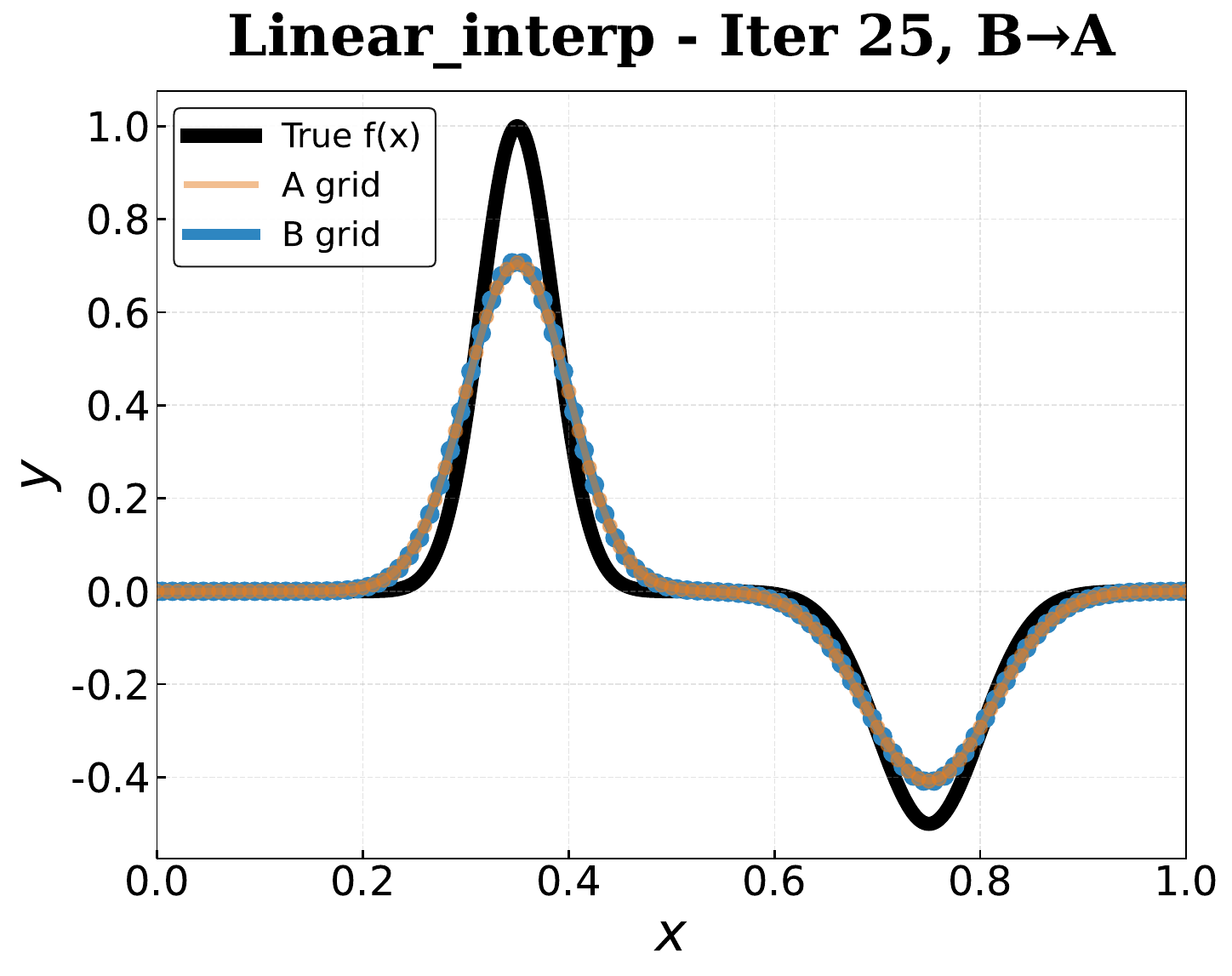}
			\caption{}
			\label{fig:APP1DpeakL-C}
		\end{subfigure}	
		\begin{subfigure}[b]{0.49\textwidth}
			\includegraphics[width=0.49\textwidth]{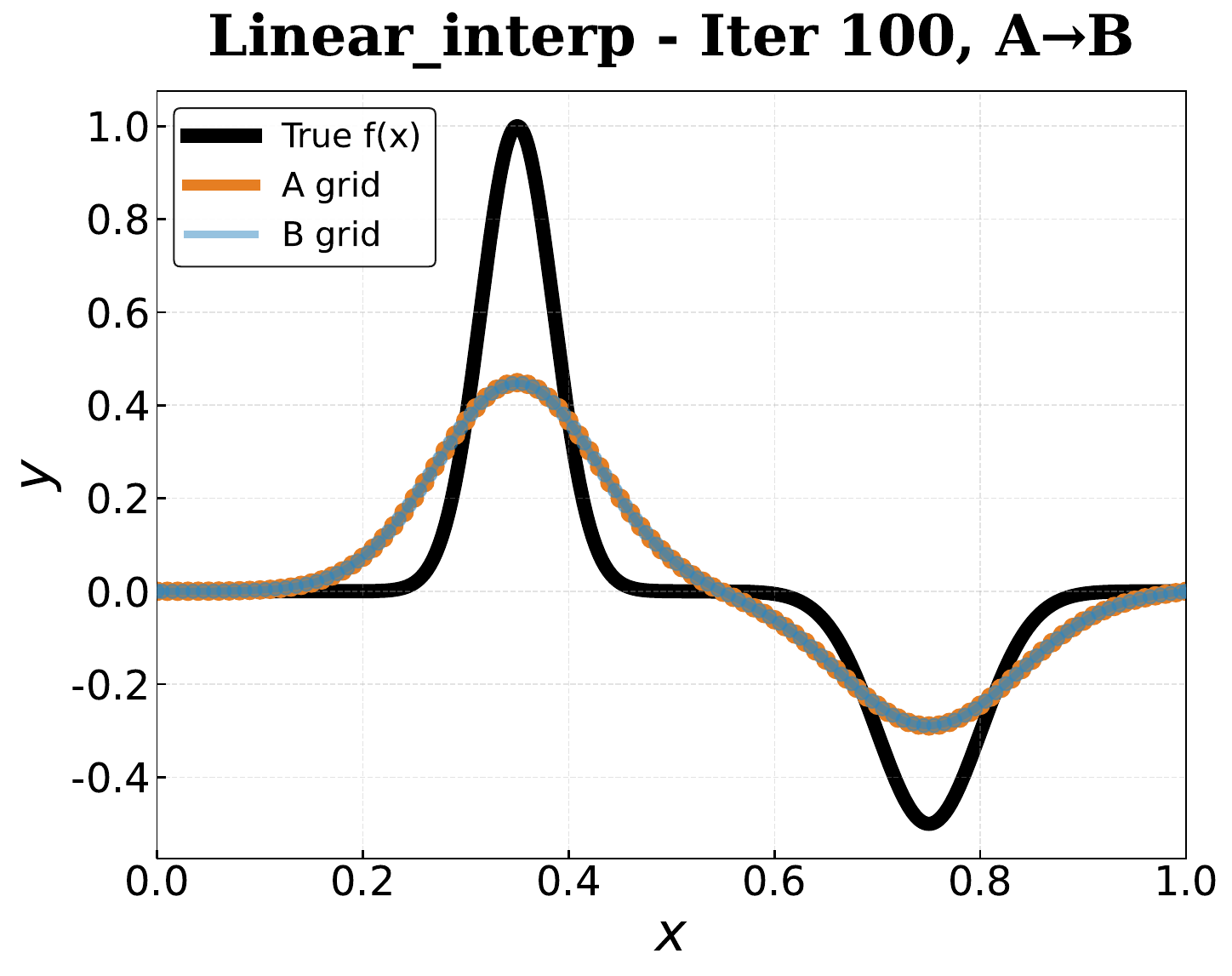}
			\includegraphics[width=0.49\textwidth]{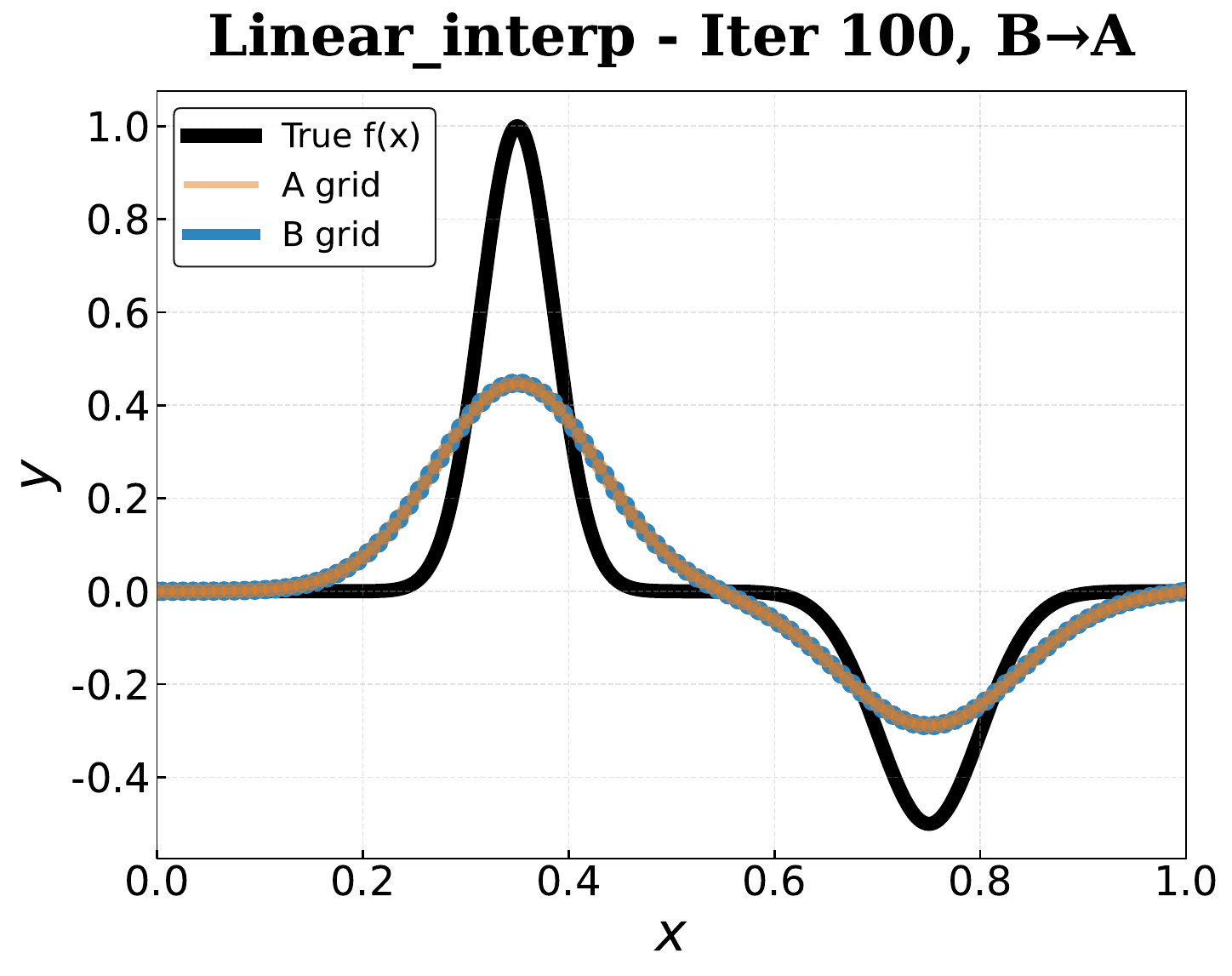}
			\caption{}
			\label{fig:APP1DpeakL-D}
		\end{subfigure}	
		\caption{\Cref{exm1d2}
				Piecewise linear interpolation results for a 1D multi-peak problem.
                (A)–(D) Interpolated solutions after 1, 5, 25 and 100 transfer iterations.}
        \label{fig:APP1DpeakL}
	\end{figure}

\begin{figure}[htbp]
    \centering
		\begin{subfigure}[b]{0.49\textwidth}
			\includegraphics[width=0.49\textwidth]{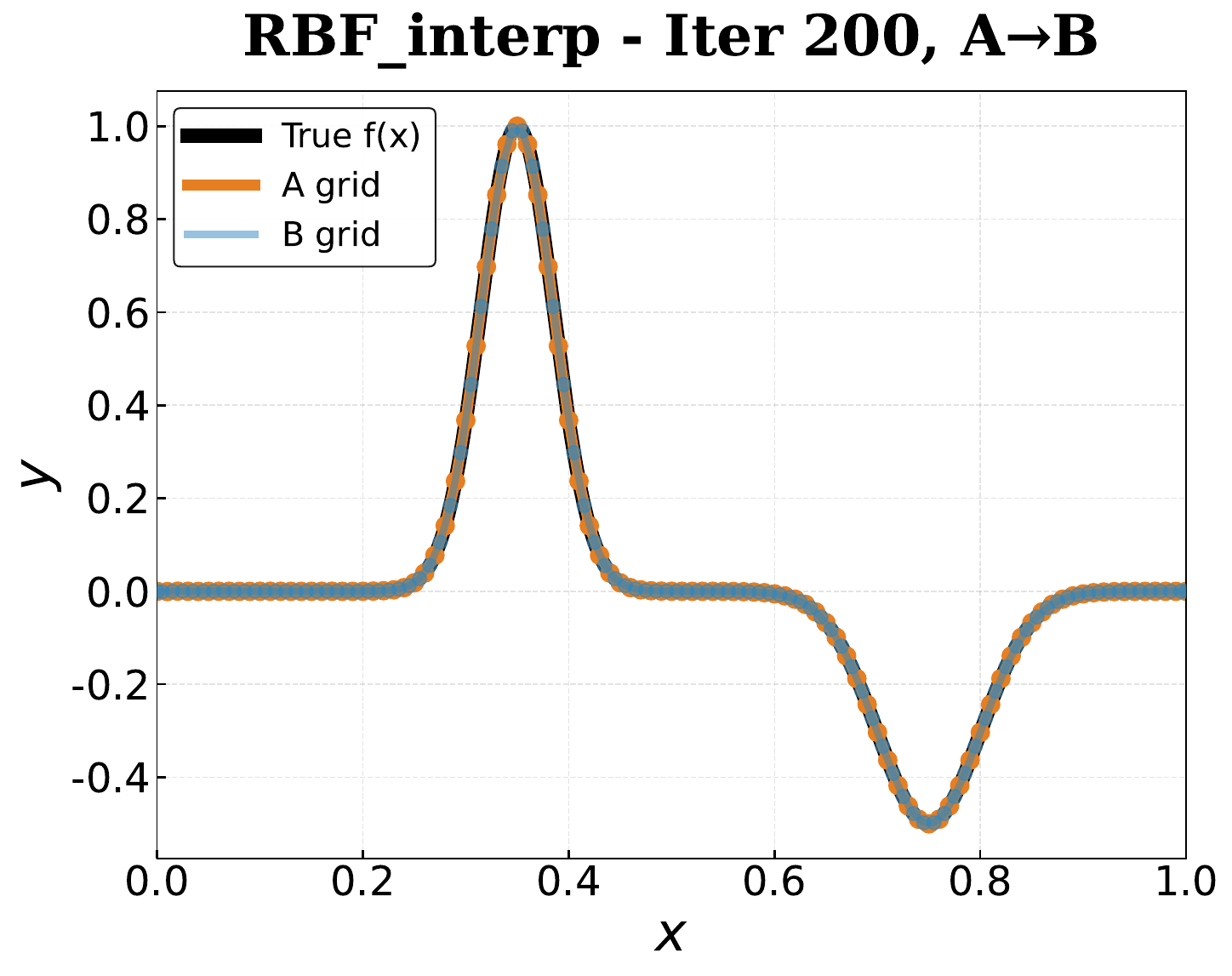}
			\includegraphics[width=0.49\textwidth]{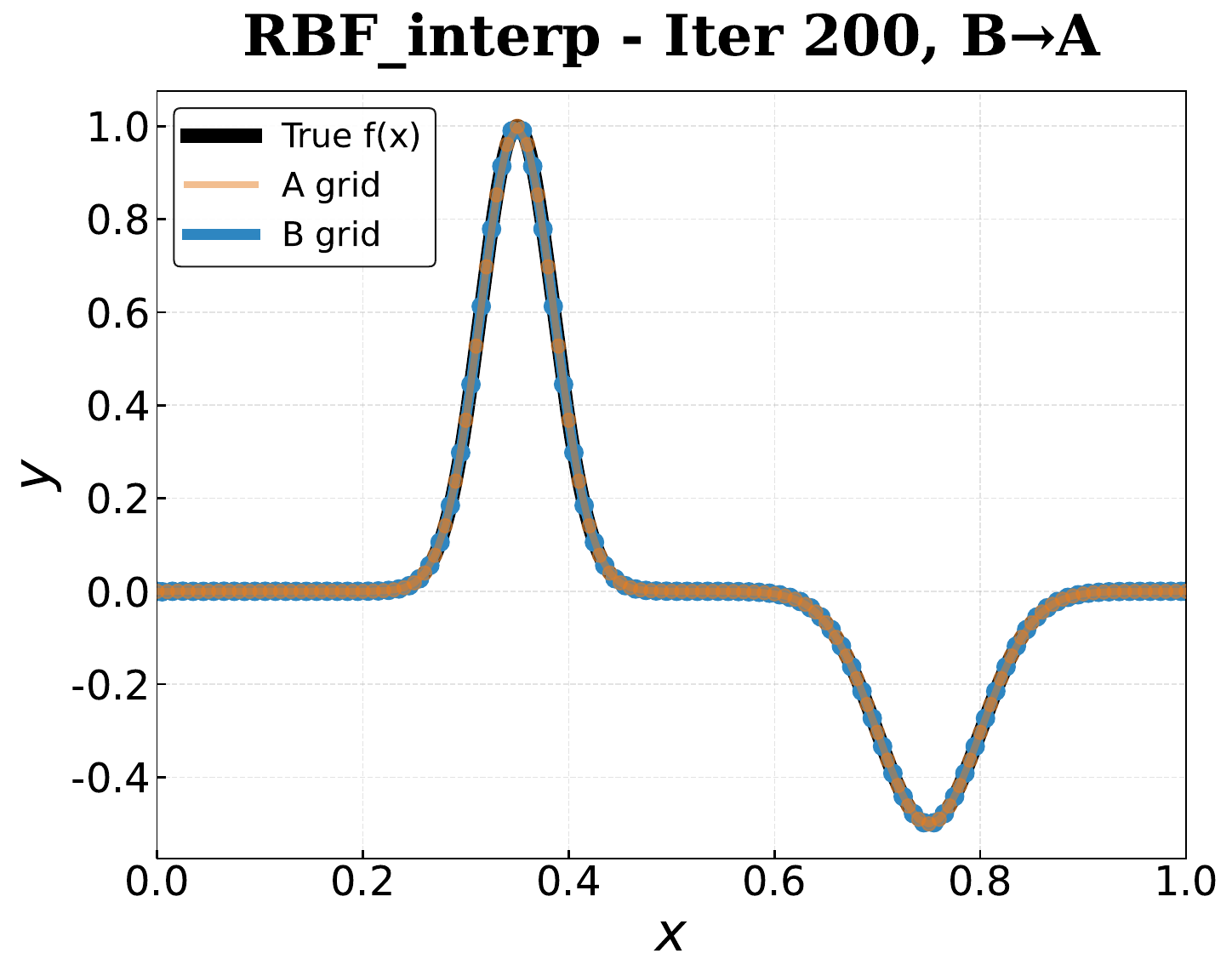}
			\caption{}
			\label{fig:APP1DpeakP-A}
		\end{subfigure}
		\begin{subfigure}[b]{0.49\textwidth}
			\includegraphics[width=\textwidth]{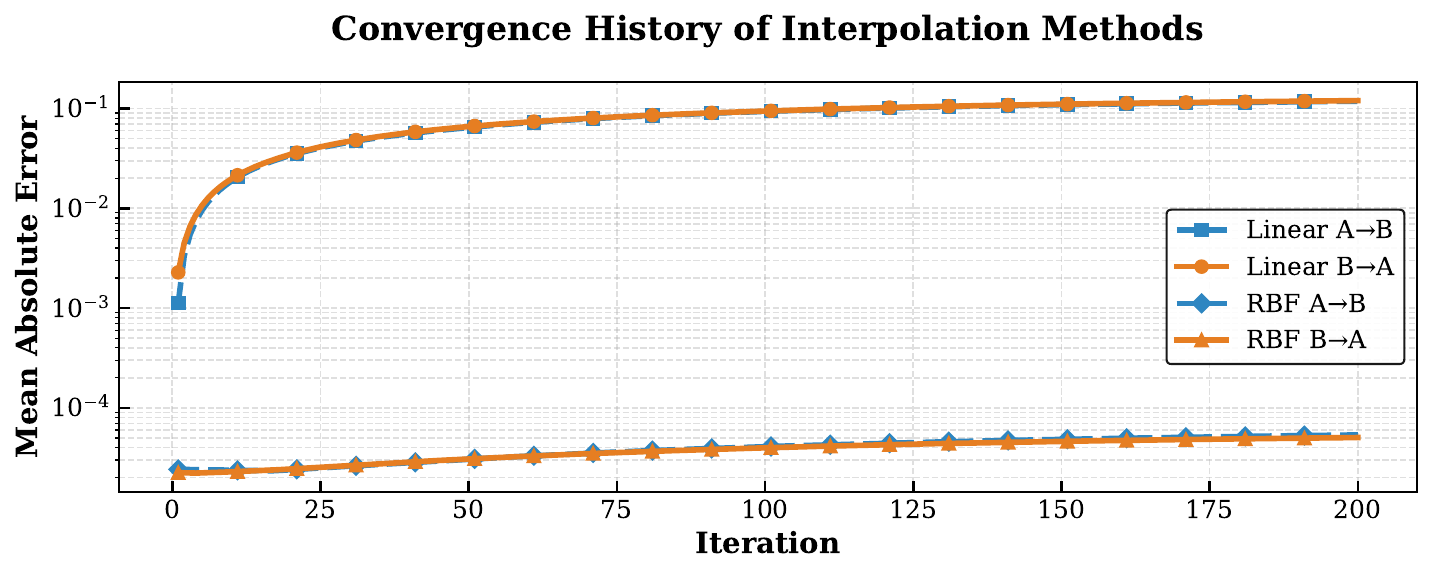}
			\caption{}
			\label{fig:APP1DpeakP-B}
		\end{subfigure}	
		\caption{\Cref{exm1d2}
				RBF-ELM interpolation performance and iterative error on multi-peak function.
                (A) Interpolated solution after 200 iterations;
                (B) Error evolution history over transfer iterations.}
        \label{fig:APP1DpeakP}
	\end{figure}

    The effect of repeated transfer becomes more pronounced when the reference solution contains stronger localized variation. 
    \begin{example}\label{exm1d2}
        Consider the one-dimensional multi-peak field 
        \[
        u(x)=\exp\!\left(-400(x-0.35)^2\right)-\frac{1}{2}\exp\!\left(-200(x-0.75)^2\right).
        \]
        Both the mesh and RBF-ELM parameter settings remain the same as in \Cref{exm1d1}. 
    \end{example}
    The detailed evolution of the interpolation process for this multi-peak function is illustrated in \Cref{fig:APP1DpeakL}.
    In this figure, the black curve represents the reference solution, while the orange and blue markers correspond to the nodal points of the meshes \(\mathcal{T}_A\) and \(\mathcal{T}_B\), respectively. 
    The thick and thin curves show the profile before and after interpolation, respectively. 
    For standard finite element interpolation, repeated transfer progressively distorts the profile and weakens the local curvature near the two peaks. 
    After only five iterations (\Cref{fig:APP1DpeakL-B}), the field is already visibly flatter. 
    As the number of iterations increases to \(25\) and \(100\) (\Cref{fig:APP1DpeakL-C}-\Cref{fig:APP1DpeakL-D}), the peak structure is almost completely lost and the solution tends toward an essentially linear profile. 
    By contrast, the RBF-ELM interpolation preserves the profile much more accurately. 
    Even after \(200\) iterations, the interpolated solution remains visually indistinguishable from the reference field, as shown in \Cref{fig:APP1DpeakP-A}. 
    This qualitative observation is further supported by \Cref{fig:APP1DpeakP-B}, where the error evolution confirms that RBF-ELM accumulates transfer error much more slowly than standard finite element interpolation. 

\begin{figure}[htbp]
    \centering
        \begin{subfigure}[b]{0.3\textwidth}
            \includegraphics[width=\linewidth]{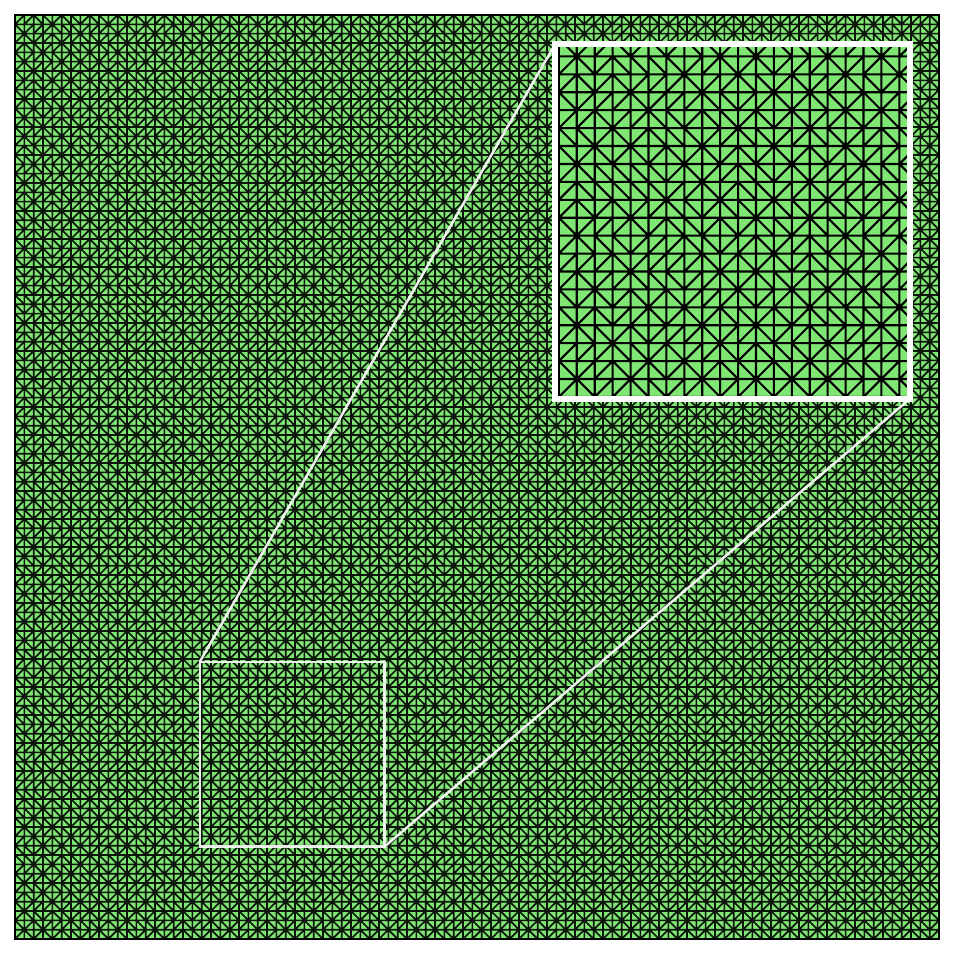}
            \caption{}
            \label{fig:APP2Dmesh-A}
        \end{subfigure}
        \hspace{2.5cm}
        \begin{subfigure}[b]{0.3\textwidth}   
            \includegraphics[width=\linewidth]{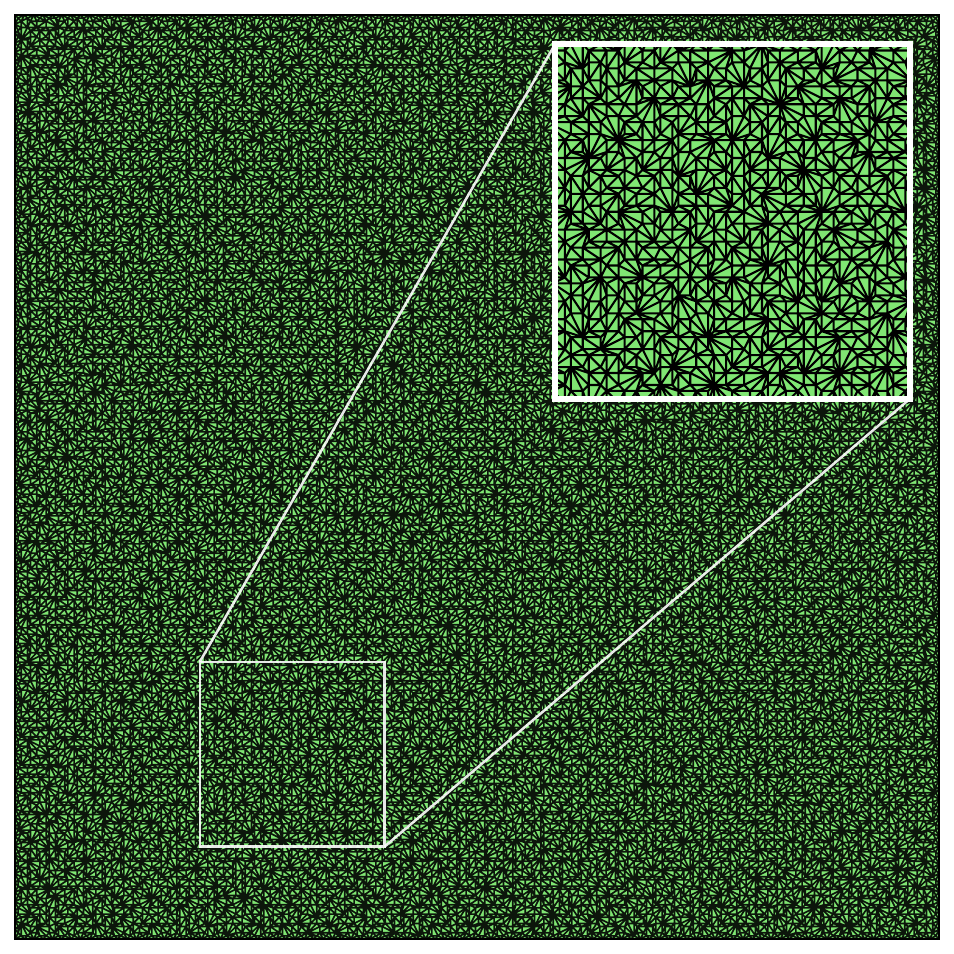}
            \caption{}
            \label{fig:APP2Dmesh-B}
        \end{subfigure}
        \caption{\Cref{exm2d1}
                Mesh configuration for 2D interpolation.
                (A) Uniform triangular mesh \(\mathcal{T}_A\); 
                (B) Delaunay triangular mesh \(\mathcal{T}_B\).}
        \label{fig:APP2Dmesh}
    \end{figure}

\begin{figure}[htbp]
    \centering
        \begin{subfigure}[b]{0.39\textwidth}        
			\includegraphics[width=\textwidth]{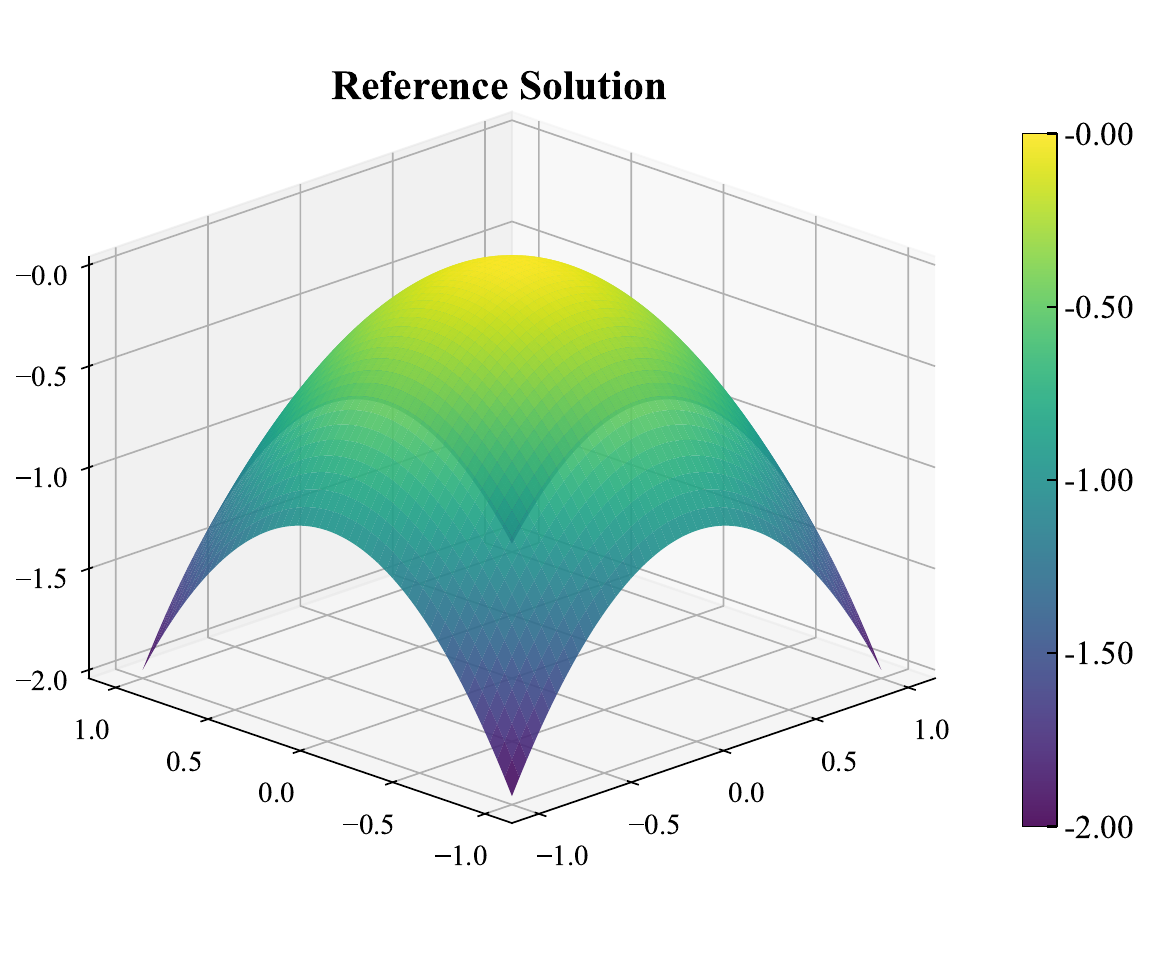}
			\caption{}
			\label{fig:APP2Dsmooth-A}
		\end{subfigure}
		\begin{subfigure}[b]{0.59\textwidth}
			\includegraphics[width=\textwidth]{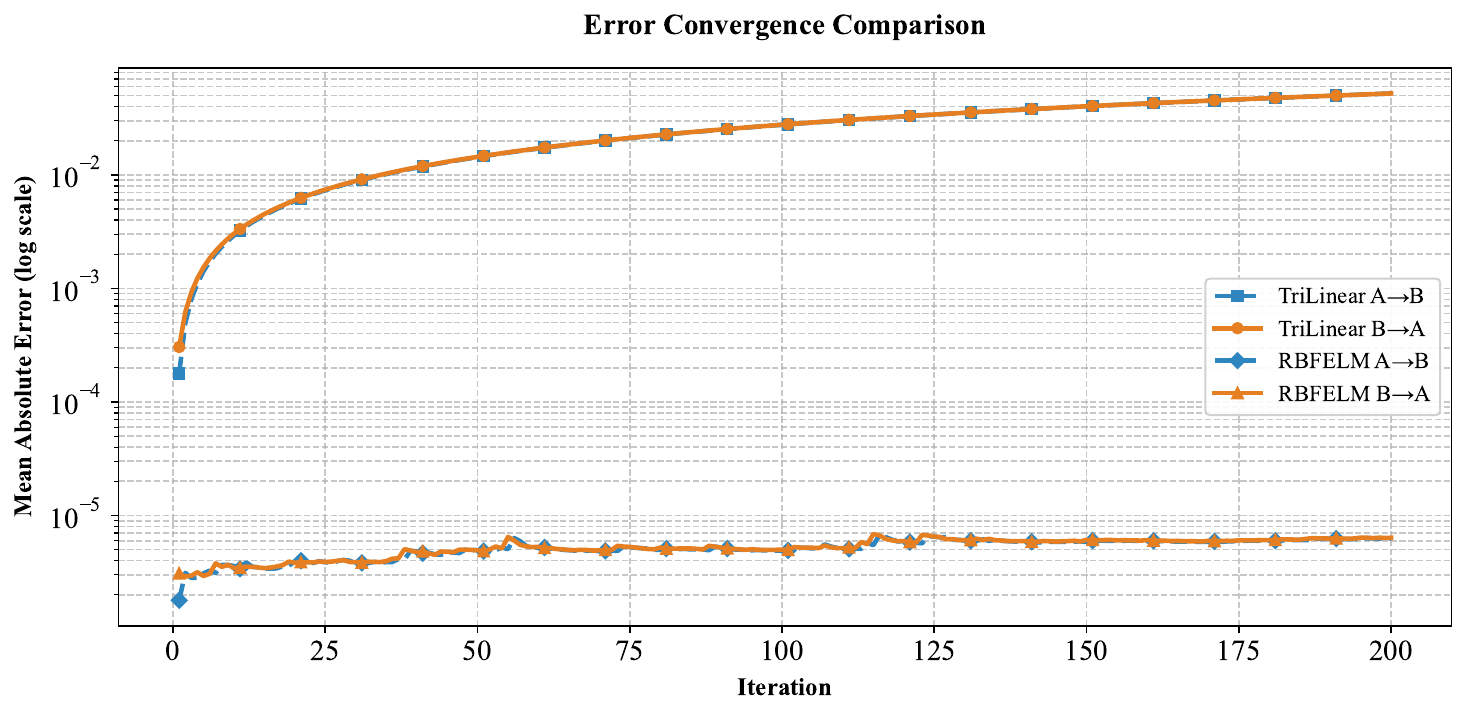}
			\caption{}
			\label{fig:APP2Dsmooth-B}
		\end{subfigure}	
		\caption{\Cref{exm2d1}
				Convergence behavior of interpolation errors for a smooth 2D problem.
				(A) Reference solution; 
				(B) Error evolution history over transfer iterations.}
        \label{fig:APP2Dsmooth1}
	\end{figure}

\begin{figure}[htbp]
    \centering
		\begin{subfigure}[b]{0.49\textwidth}            
			\includegraphics[width=\textwidth]{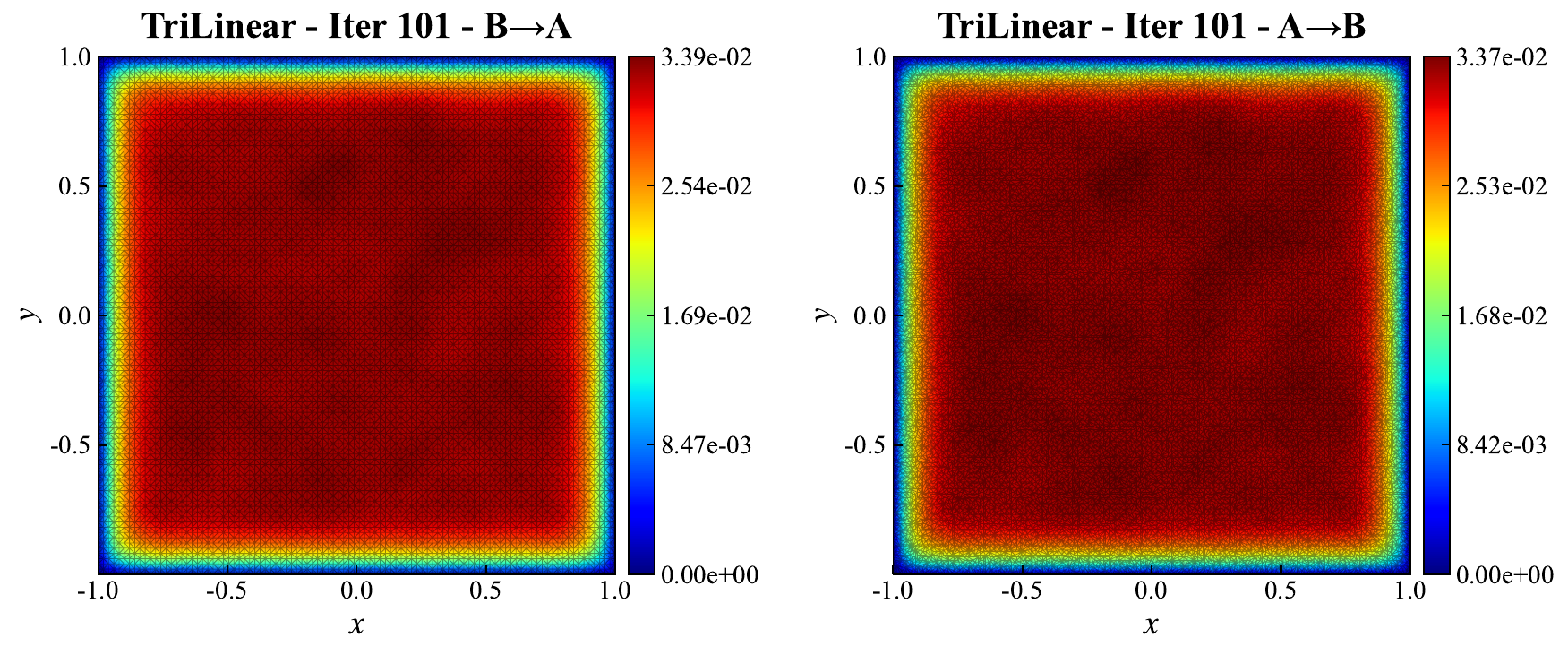}
			\caption{}
			\label{fig:APP2Dsmooth-C}
		\end{subfigure}
		\begin{subfigure}[b]{0.49\textwidth}
			\includegraphics[width=\textwidth]{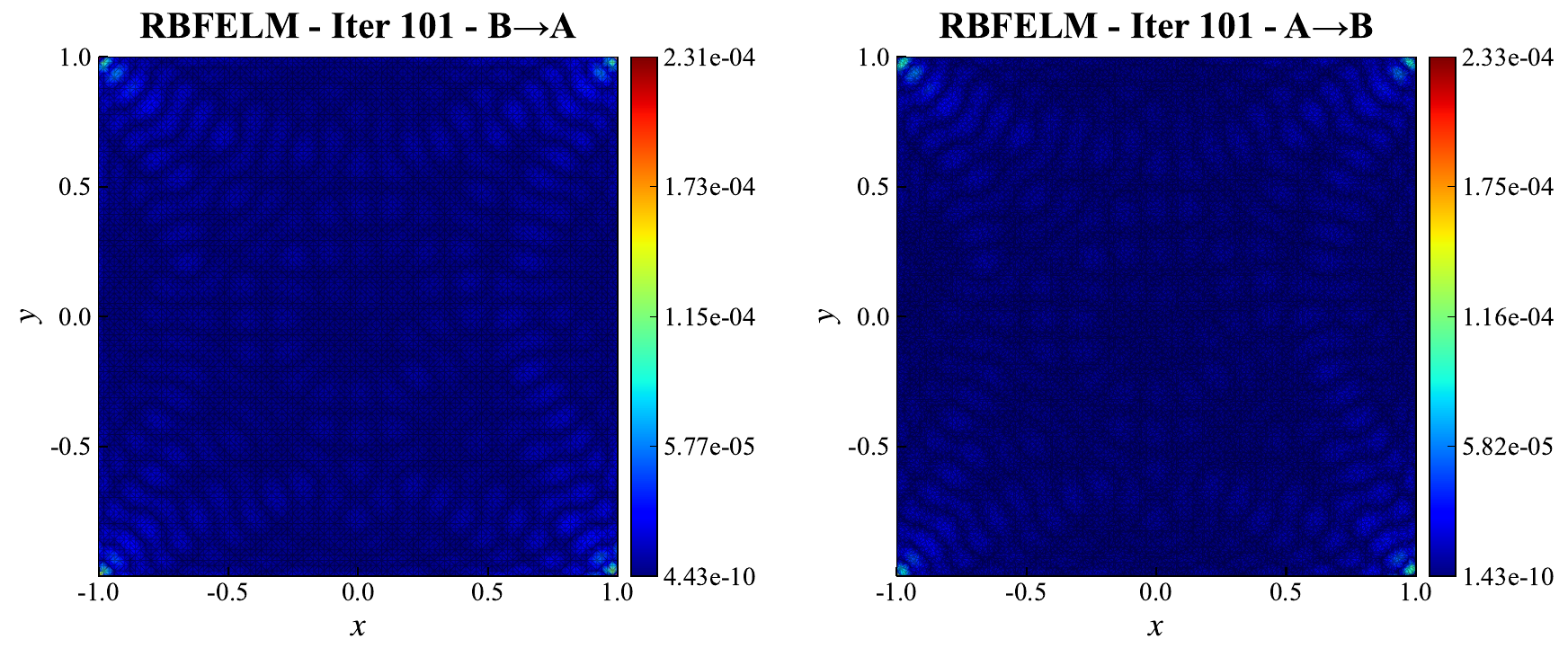}
			\caption{}
			\label{fig:APP2Dsmooth-D}
		\end{subfigure}	
		\caption{\Cref{exm2d1}
				spatial distributions of interpolation errors for a smooth 2D problem.
                (A) Piecewise linear interpolated error after 100 iterations;
                (B) RBF-ELM interpolated error after 100 iterations.}
        \label{fig:APP2Dsmooth2}
	\end{figure}

    The numerical study is next carried out in two dimensions on a pair of triangular meshes. 
    \begin{example}\label{exm2d1}
        Consider the smooth reference solution (\Cref{fig:APP2Dsmooth-A})
        \[
        u(x,y)=-x^2-y^2, \qquad (x,y)\in[-1,1]^2.
        \]
        The mesh \(\mathcal{T}_A\) is a uniform triangulation with \(10{,}000\) nodes (\Cref{fig:APP2Dmesh-A}), and \(\mathcal{T}_B\) is a Delaunay triangulation constructed from the element midpoints of \(\mathcal{T}_A\) together with the boundary nodes (\Cref{fig:APP2Dmesh-B}). 
        By construction, \(\mathcal{T}_B\) typically contains about twice as many nodes as \(\mathcal{T}_A\). 
        Additionally, the RBF-ELM model is configured with \(N_c=1000\) centers randomly selected from the training nodes and a width parameter of \(\varepsilon^2=1/10\). 
    \end{example}
    After 100 transfer iterations, the error map produced by piecewise linear finite element interpolation is shown in \Cref{fig:APP2Dsmooth-C}. 
    The error remains relatively small near the boundary but becomes pronounced in the interior, leading to a global mean absolute error on the order of \(10^{-2}\). 
    By contrast, the error map in \Cref{fig:APP2Dsmooth-D} exhibits a much more localized pattern for RBF-ELM. 
    Most of the error is confined to the domain corners, and the global mean absolute error remains on the order of \(10^{-6}\). 
    This distinction persists throughout the repeated-transfer process and is further confirmed by the error history shown in \Cref{fig:APP2Dsmooth-B}. 

\begin{figure}[htbp]
    \centering
        \begin{subfigure}[b]{0.39\textwidth}       
			\includegraphics[width=\textwidth]{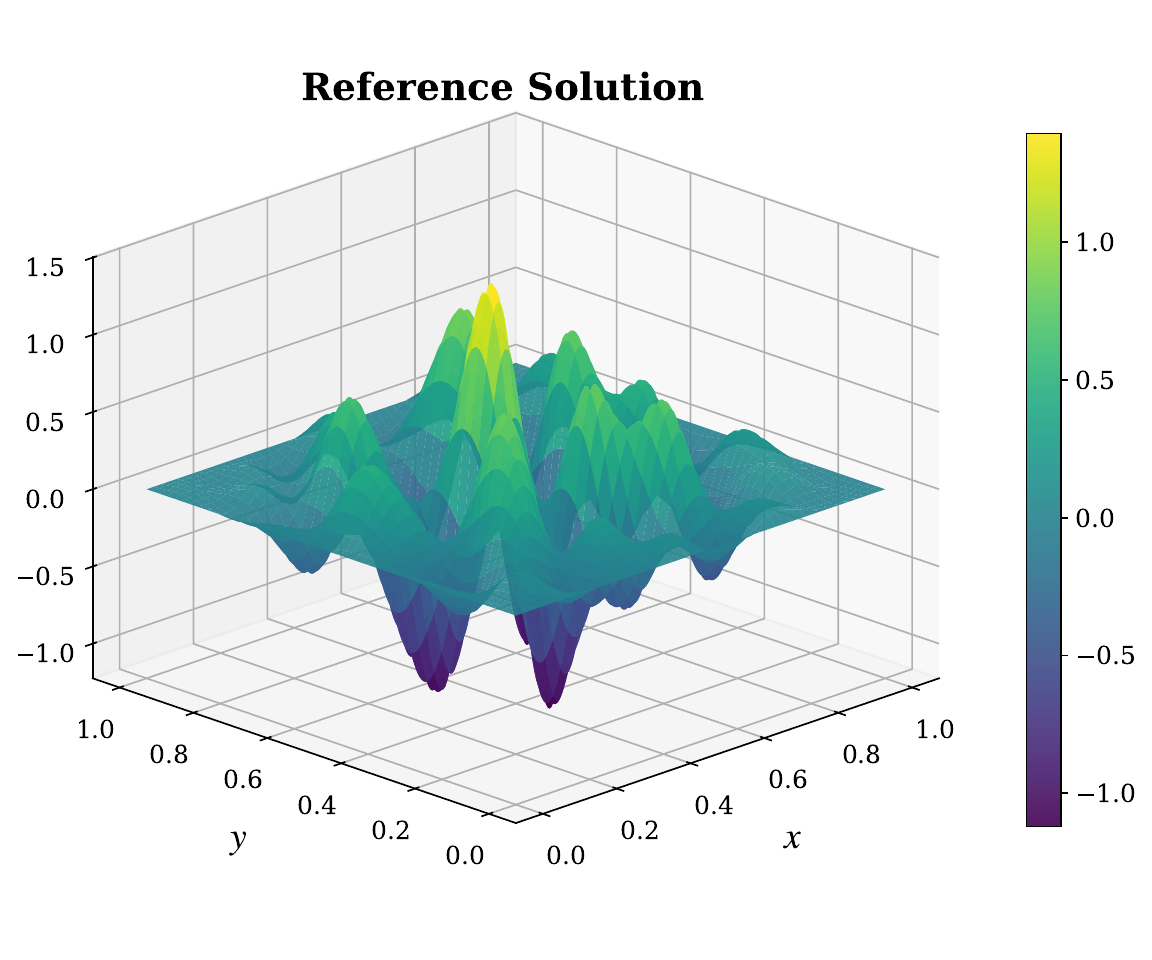}
			\caption{}
			\label{fig:APP2Dpeak-true}
		\end{subfigure}
		\begin{subfigure}[b]{0.59\textwidth}
			\includegraphics[width=\textwidth]{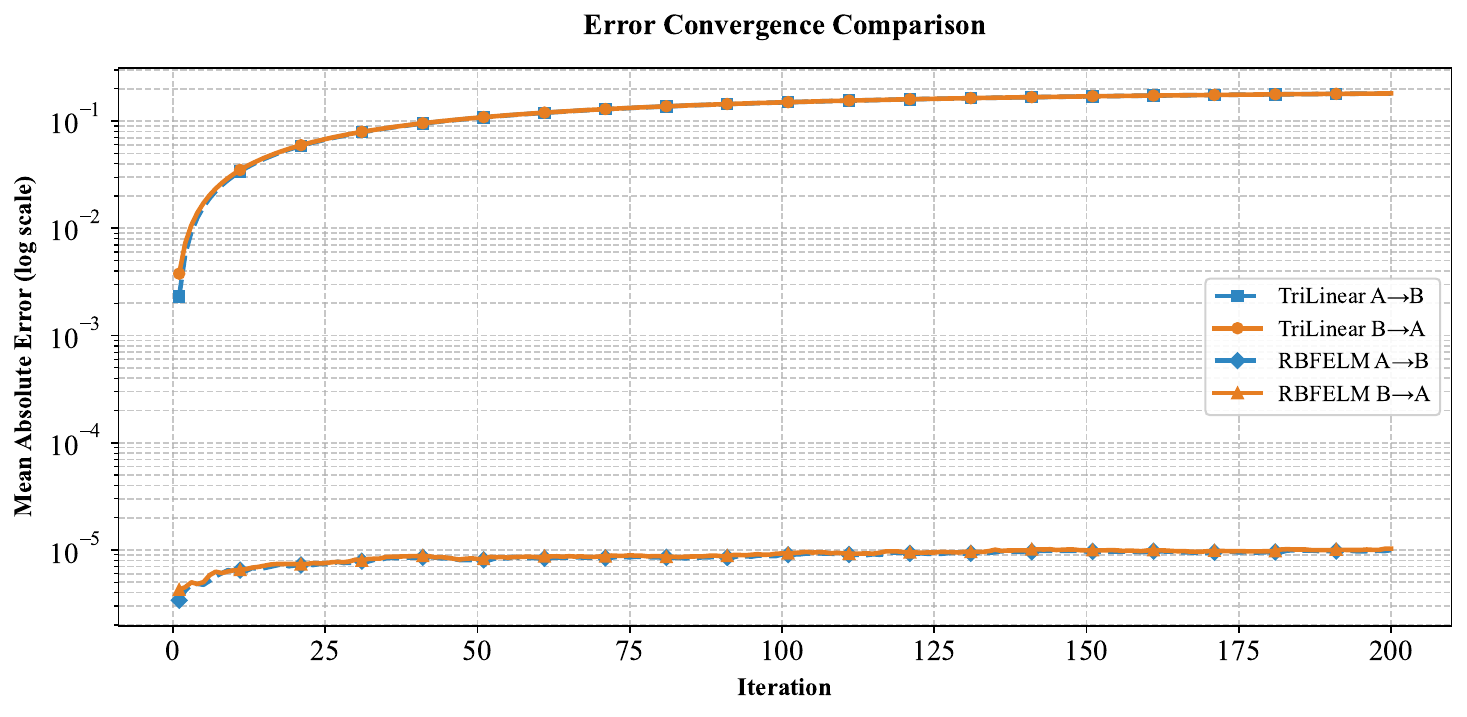}
			\caption{}
			\label{fig:APP2Dpeak-error}
		\end{subfigure}
		\caption{\Cref{exm2d2}
			Convergence behavior of interpolation errors for a 2D multi-peak problem.
			(A) Reference solution; 
			(B) Error evolution history over transfer iterations.}
        \label{fig:APP2Dpeak}
	\end{figure}

\begin{figure}[htbp]
    \centering
		\begin{subfigure}[b]{0.49\textwidth}            
			\includegraphics[width=0.49\textwidth]{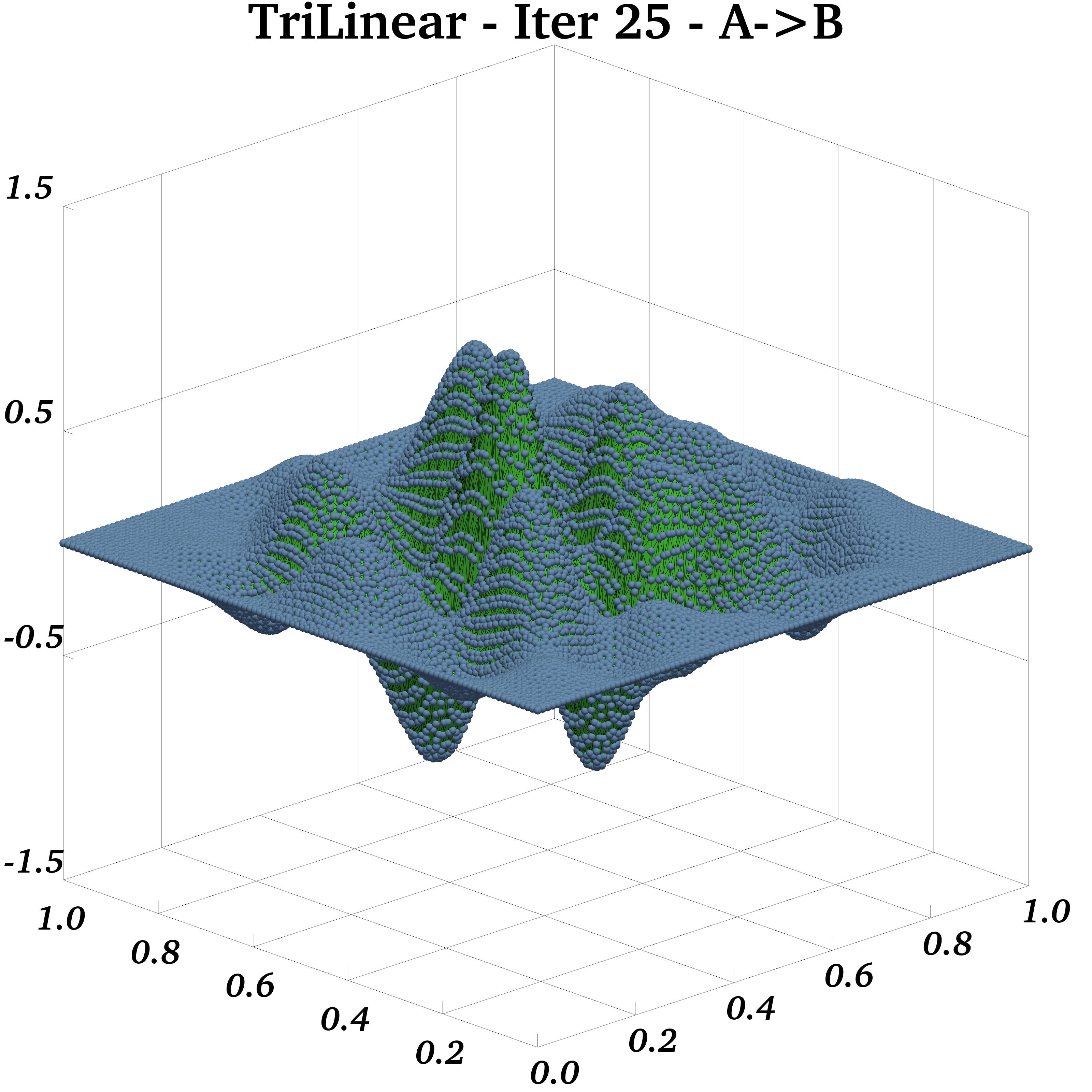}
			\includegraphics[width=0.49\textwidth]{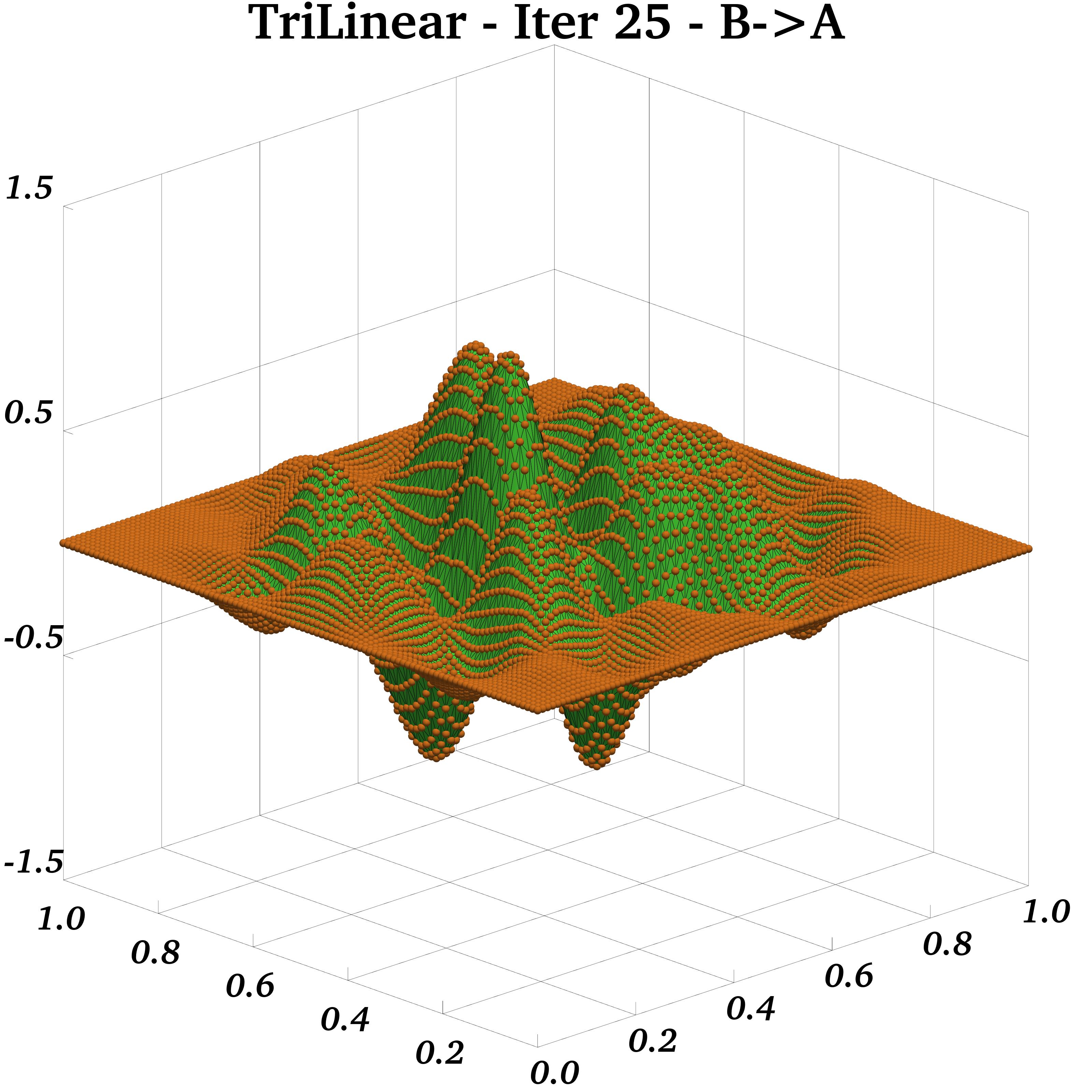}
			\caption{}
			\label{fig:APP2DpeakLR-A}
		\end{subfigure}
		\begin{subfigure}[b]{0.49\textwidth}
			\includegraphics[width=0.49\textwidth]{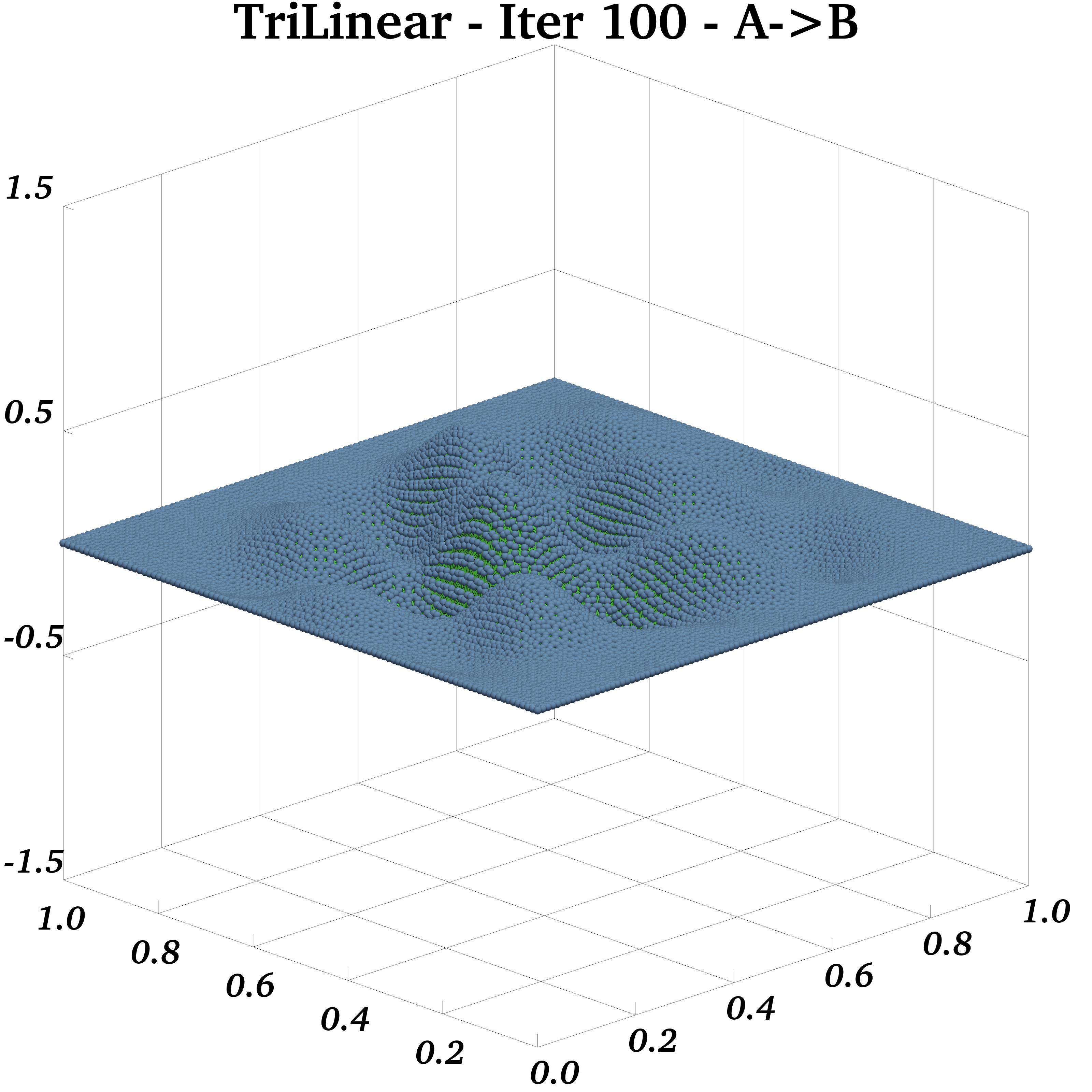}
			\includegraphics[width=0.49\textwidth]{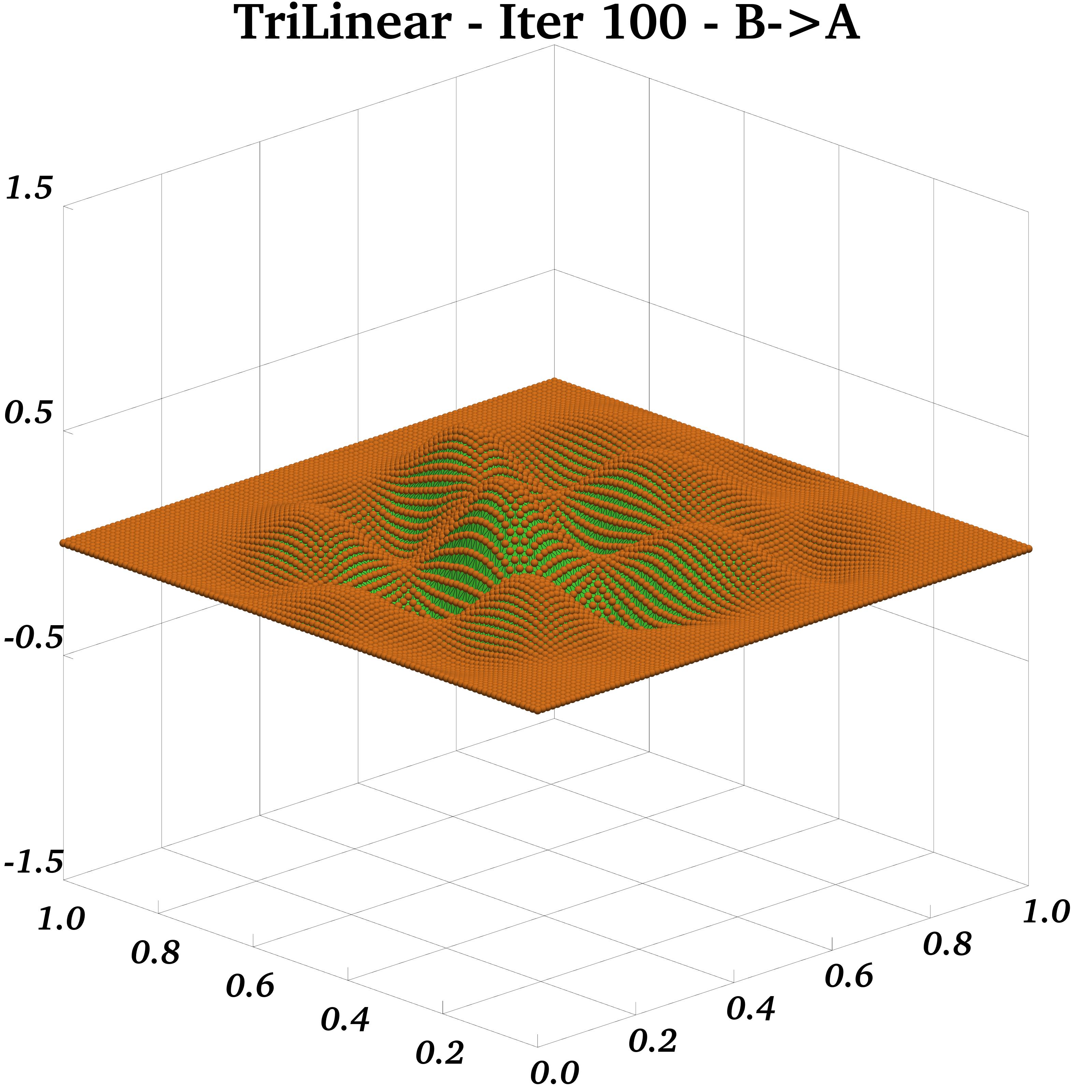}
			\caption{}
			\label{fig:APP2DpeakLR-B}
		\end{subfigure}
		\begin{subfigure}[b]{0.49\textwidth}
			\includegraphics[width=0.49\textwidth]{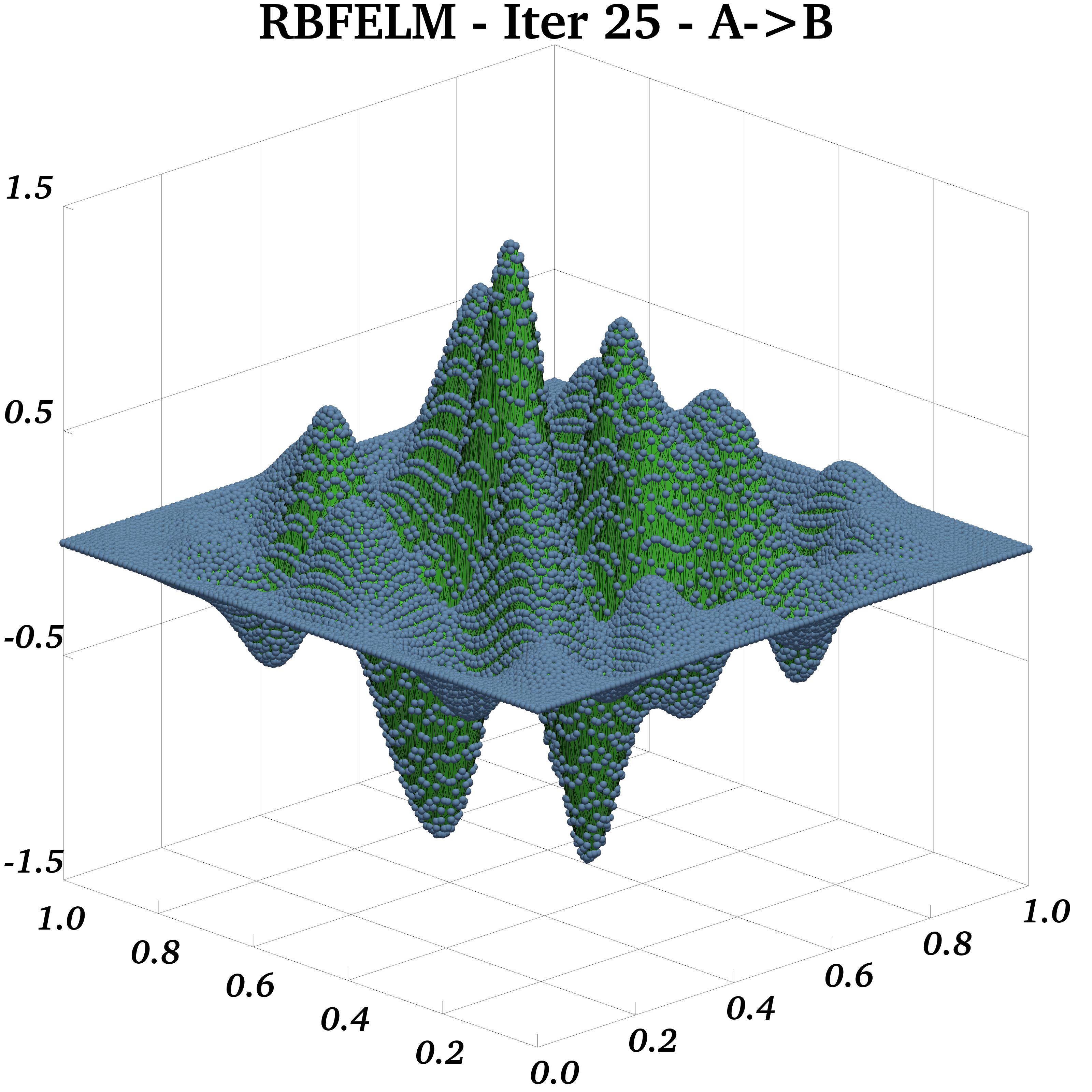}
			\includegraphics[width=0.49\textwidth]{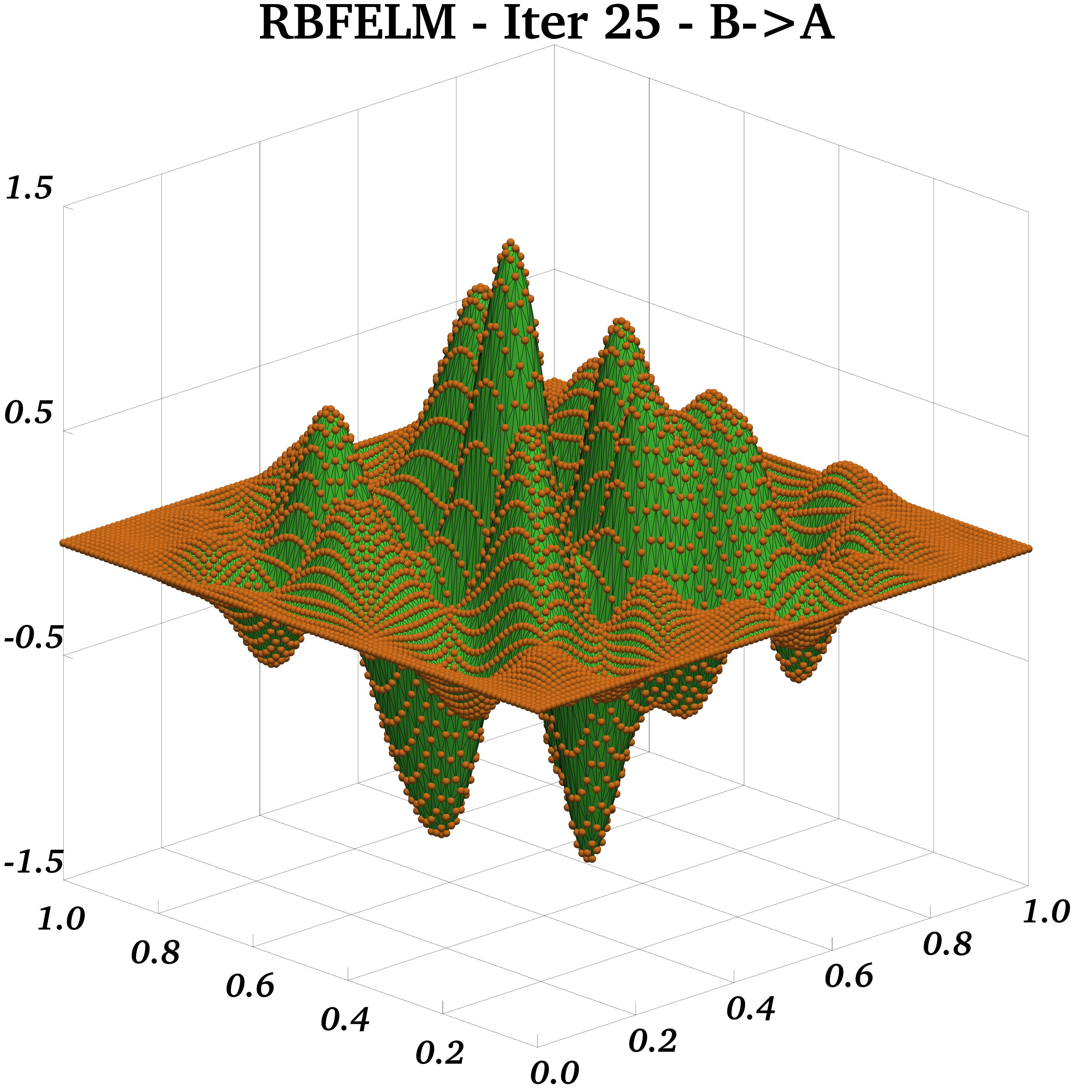}
			\caption{}
			\label{fig:APP2DpeakLR-C}
		\end{subfigure}
		\begin{subfigure}[b]{0.49\textwidth}
			\includegraphics[width=0.49\textwidth]{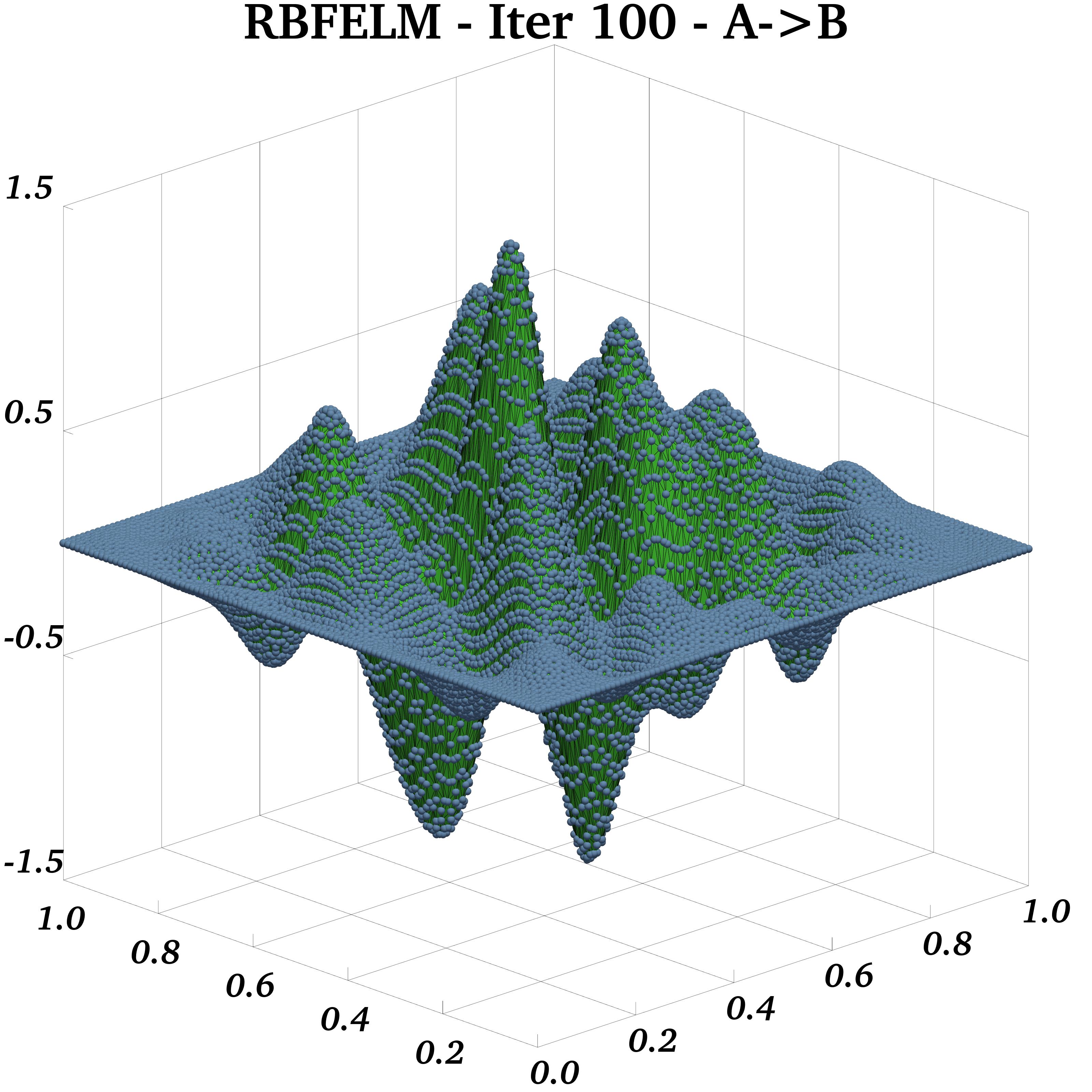}
			\includegraphics[width=0.49\textwidth]{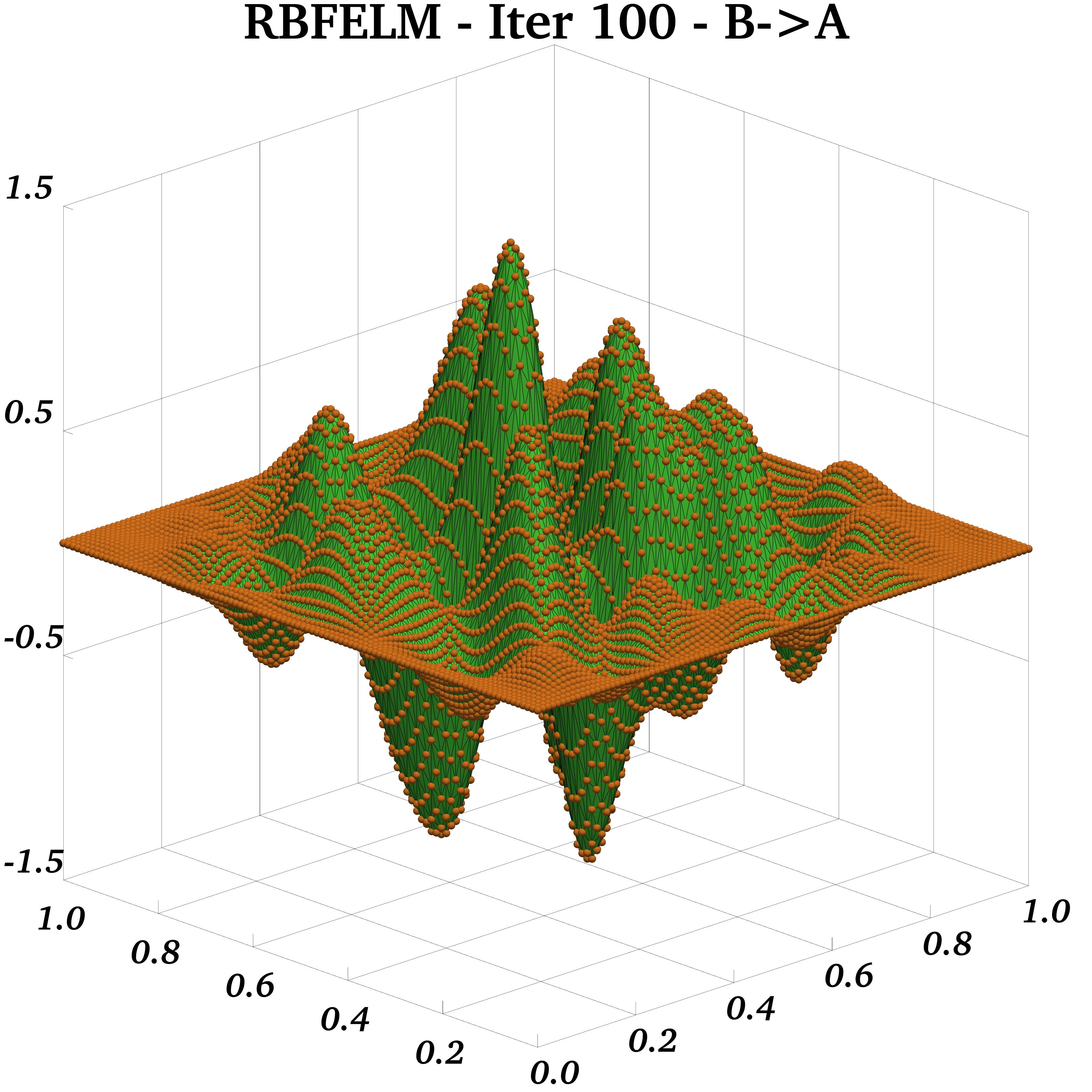}
			\caption{}
			\label{fig:APP2DpeakLR-D}
		\end{subfigure}
		\caption{\Cref{exm2d2}
			Comparison of piecewise linear and RBF-ELM interpolation results. 
			(A)-(B) piecewise linear interpolation results after 25 and 100 iterations; left/right subplots correspond to transfers \(\mathcal{T}_A \to \mathcal{T}_B\) and \(\mathcal{T}_B \to \mathcal{T}_A\), respectively;
			(C)-(D) RBF-ELM results after 25 and 100 iterations.}
		\label{fig:APP2DpeakLR}
	\end{figure}

    \begin{example}\label{exm2d2}
        Consider the two-dimensional oscillatory field defined by the reference function (\Cref{fig:APP2Dpeak-true})
        \[
        u(x, y) = \sin(\pi x)\sin(\pi y)\,\mathbf{a}^\top \mathbf{s}(x, y),
        \]
        where
        \[
        \mathbf{a} =
        \begin{pmatrix}
        1 \\
        0.6 \\
        -0.4
        \end{pmatrix},
        \qquad
        \mathbf{s}(x, y) =
        \begin{pmatrix}
        \sin(6\pi x)\cos(4\pi y) \\
        \sin\big(8\pi(x-0.3)\big)\cos\big(6\pi(y+0.2)\big) \\
        \cos(10\pi x)\sin(8\pi y)
        \end{pmatrix}.
        \]
        The mesh configuration and parameter settings are the same as those in \Cref{exm2d1}, except for the common width parameter, which is set to \(\varepsilon^2=1/50\).
    \end{example}
    The results produced by piecewise linear finite element interpolation after \(25\) and \(100\) transfer iterations are displayed in \Cref{fig:APP2DpeakLR-A,fig:APP2DpeakLR-B}. 
    A noticeable loss of structural detail is already observed after \(25\) iterations, and after \(100\) iterations much of the oscillatory structure has been smeared out. 
    The corresponding RBF-ELM results are shown in \Cref{fig:APP2DpeakLR-C,fig:APP2DpeakLR-D}. 
    Even after \(25\) and \(100\) iterations, the reconstructed field still preserves both the global structure and the fine-scale variations of the reference solution. 
    A quantitative comparison of error growth is given in \Cref{fig:APP2Dpeak-error}. 
    Although the errors of both methods increase under repeated transfer, the growth associated with piecewise linear finite element interpolation is markedly faster than that of RBF-ELM. 

\begin{figure}[htbp]
    \centering
        \begin{subfigure}[b]{0.39\textwidth}            
			\includegraphics[width=\textwidth]{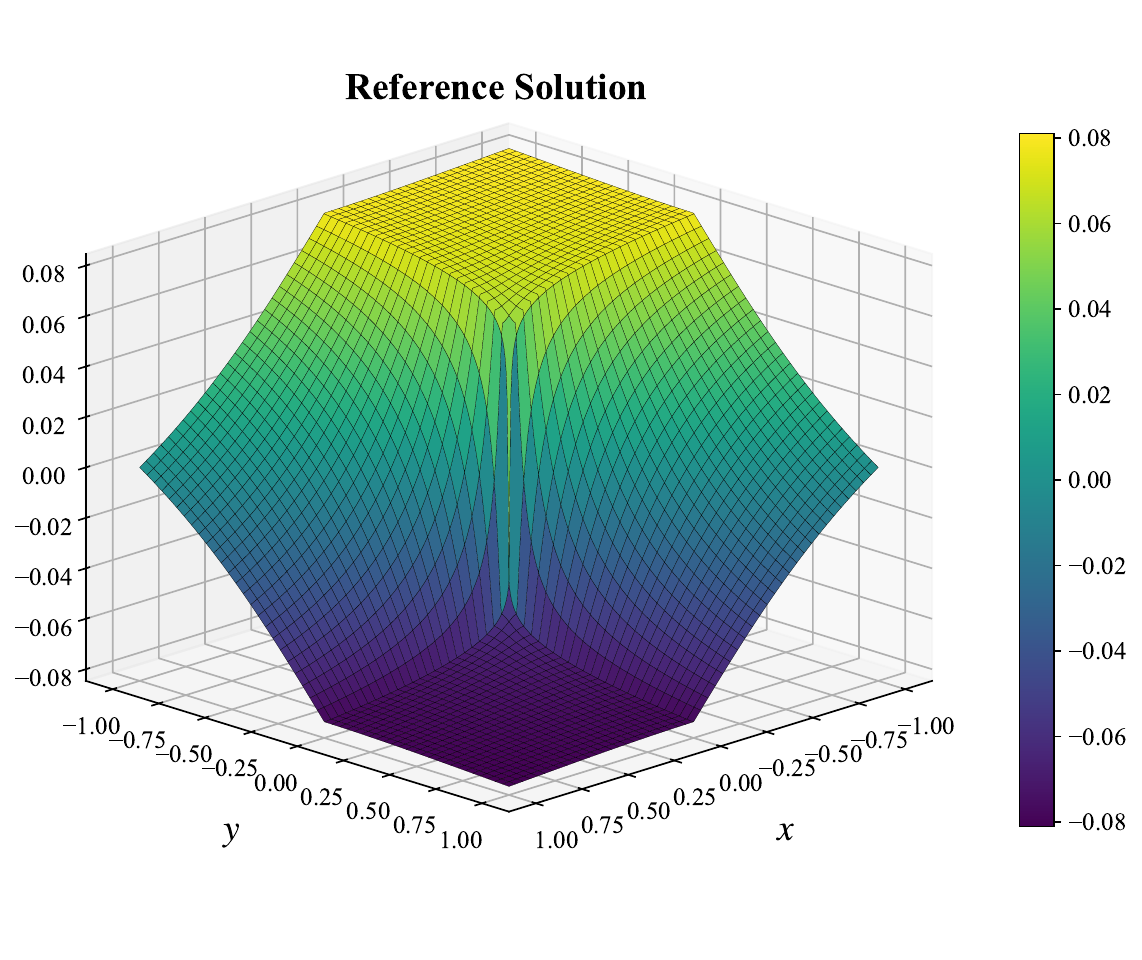}
			\caption{}
			\label{fig:APPbenchmark-A}
		\end{subfigure}
		\begin{subfigure}[b]{0.59\textwidth}
			\includegraphics[width=\textwidth]{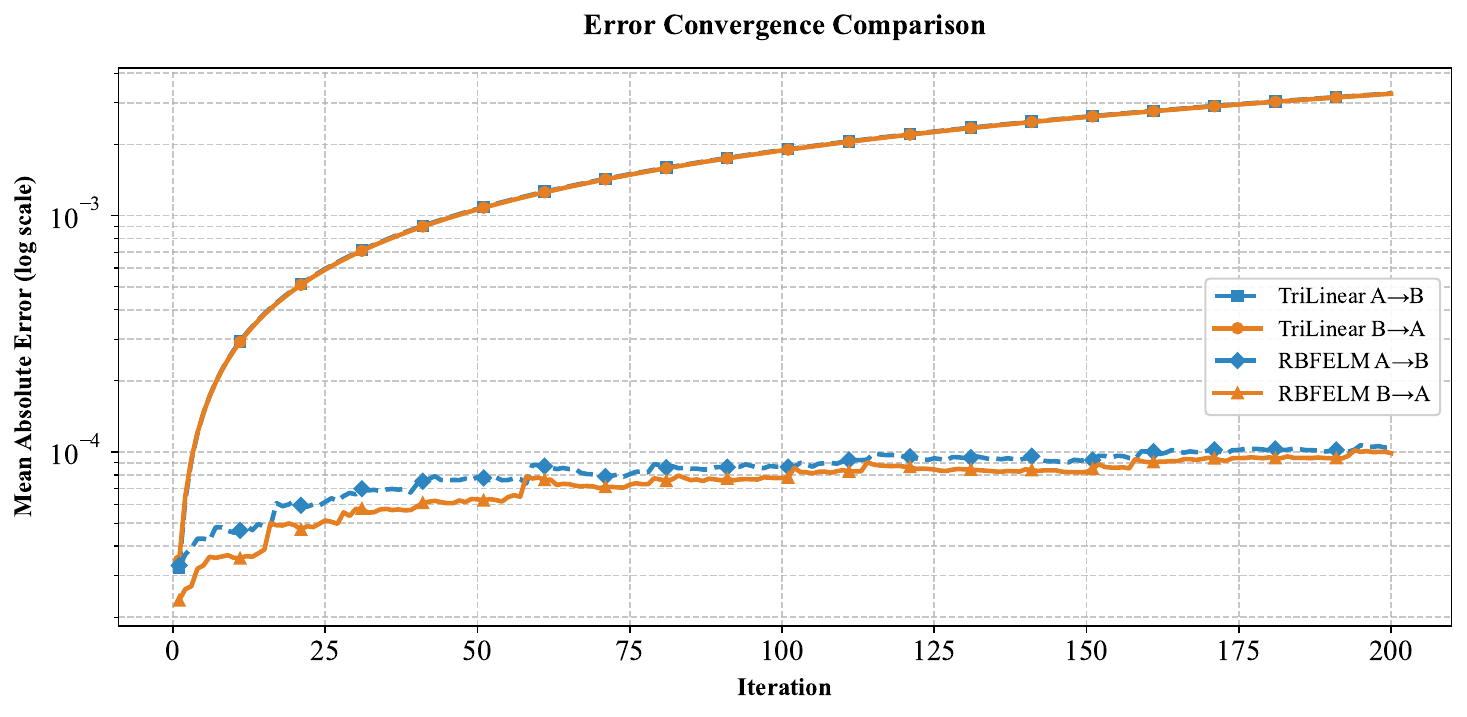}
			\caption{}
			\label{fig:APPbenchmark-B}
		\end{subfigure}	
		\caption{
				Convergence behavior of interpolation errors for a singular 2D problem.
				(A) Reference solution; 
				(B) Error evolution history over transfer iterations.}
        \label{fig:APPbenchmark1}
	\end{figure}

\begin{figure}[htbp]
    \centering
		\begin{subfigure}[b]{0.49\textwidth}            
			\includegraphics[width=\textwidth]{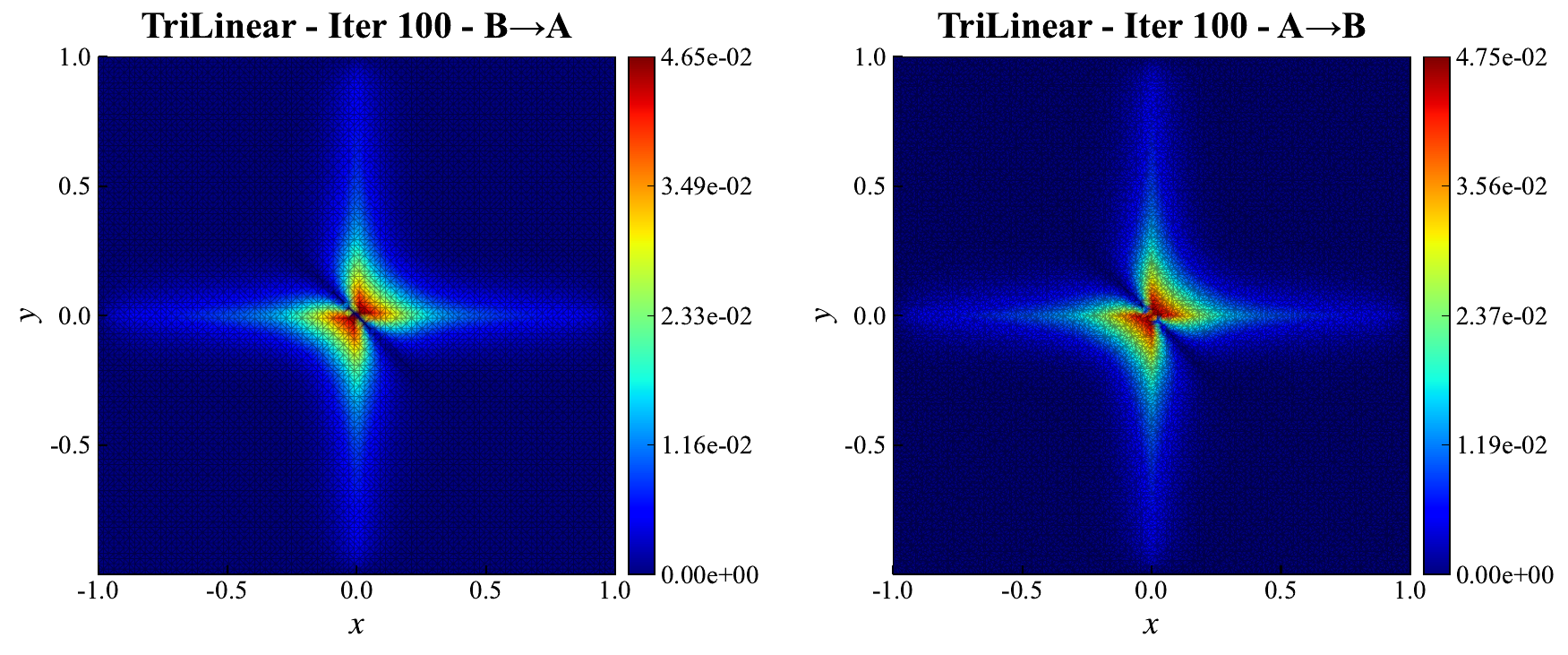}
			\caption{}
			\label{fig:APPbenchmark-C}
		\end{subfigure}
		\begin{subfigure}[b]{0.49\textwidth}
			\includegraphics[width=\textwidth]{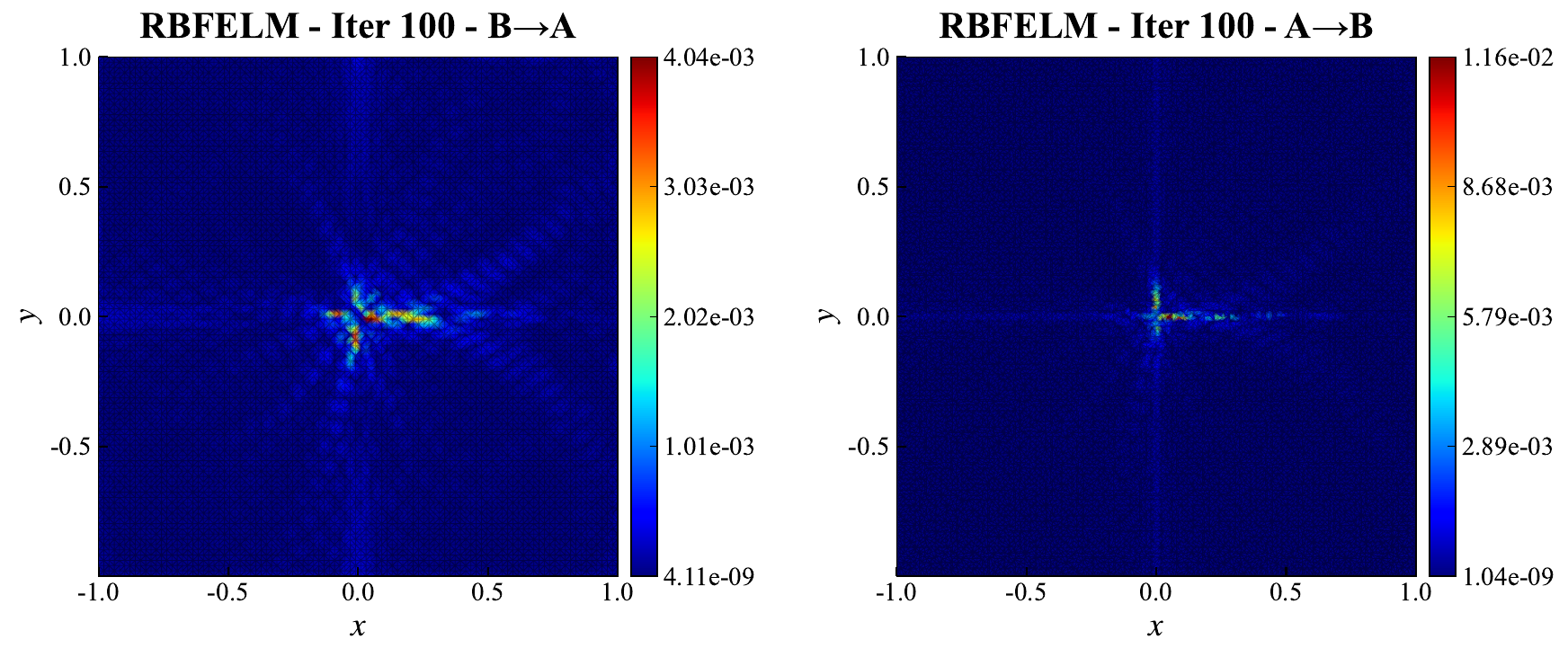}
			\caption{}
			\label{fig:APPbenchmark-D}
		\end{subfigure}	
		\caption{
				spatial distributions of interpolation errors for a singular 2D problem.
                (A) Piecewise linear interpolated error after 100 iterations;
                (B) RBF-ELM interpolated error after 100 iterations.}
        \label{fig:APPbenchmark2}
	\end{figure}

    Our final example focuses on a two-dimensional field with reduced regularity. 
    The reference solution, which is taken from a benchmark interface problem in finite element analysis, is used to assess the performance of repeated transfer in the presence of a weak singularity at the origin. 
    \begin{example}\label{exm2d3}
        In polar coordinates, define the reference function with a fractional exponent (\Cref{fig:APPbenchmark-A})
        \[
        u(r,\theta)=r^{\gamma}\mu(\theta),\qquad \gamma=0.1,
        \]
        where
        \[
        r=\sqrt{x^2+y^2},\qquad \theta=\operatorname{atan2}(y,x),
        \]
        and the angular component \(\mu(\theta)\) is given by
        \[
        \mu(\theta)=
        \begin{cases}
        \cos\!\big((\pi/2-\sigma)\gamma\big)\cos\!\big((\theta-\pi/2+\rho)\gamma\big), & 0\le\theta<\pi/2, \\[4pt]
        \cos(\rho\gamma)\cos\!\big((\theta-\pi+\sigma)\gamma\big), & \pi/2\le\theta<\pi, \\[4pt]
        \cos(\sigma\gamma)\cos\!\big((\theta-\pi-\rho)\gamma\big), & \pi\le\theta<3\pi/2, \\[4pt]
        \cos\!\big((\pi/2-\rho)\gamma\big)\cos\!\big((\theta-3\pi/2-\sigma)\gamma\big), & 3\pi/2\le\theta<2\pi,
        \end{cases}
        \]
        with \(\sigma=-14.92256510455152\) and \(\rho=\pi/4\). 
        The mesh configuration shown in \Cref{fig:APP2Dmesh} is used here as well, but the input features are transformed from Cartesian to polar coordinates before being fed into the RBF-ELM model. 
        The parameters are taken as \(N_c=1500\) and \(\varepsilon^2=1/25\). 
    \end{example}
    The error history in \Cref{fig:APPbenchmark-B} shows that the error of piecewise linear finite element interpolation grows rapidly throughout the repeated-transfer process, whereas the RBF-ELM error increases only in the initial stage and then remains essentially stable. 
    This difference once again indicates that RBF-ELM is more effective in controlling error accumulation near the singularity. 
    The error map at a representative transfer iteration is shown in \Cref{fig:APPbenchmark2}. 
    For piecewise linear finite element interpolation, the error spreads outward from the singular point and becomes pronounced in the first and third quadrants. 
    In the RBF-ELM result, most of the error remains confined to a small neighborhood of the singularity. 
    Further improvement is expected with adaptive mesh refinement and more careful parameter tuning. 

\section{Conclusions}\label{sec:CONCLU}
This work presents a systematic investigation of three neural network frameworks for transferring function data across distinct unstructured and non-nested meshes, along with a comparative assessment of their performance. The MLP-based approach delivers stable interpolation and moderate generalization across all test cases; however, its accuracy is limited to approximately \(10^{-3}\), and it incurs considerable computational cost due to its iterative training process. The ELM improves upon this performance by achieving significantly higher accuracy with minimal training effort, though its performance is sensitive to network architecture, parameter initialization, and the distribution of training samples. Augmenting the training set with Gaussian points in sparse regions can partially mitigate this sensitivity, but it also introduces additional approximation errors and lacks a well-defined stopping criterion.
By contrast, the RBF-ELM approach offers a favorable balance between accuracy, computational efficiency, and implementation simplicity. It exhibits reduced sensitivity to hyperparameters and maintains stable performance under irregular sampling conditions. Furthermore, in successive multi-step transfer experiments, the RBF-ELM demonstrates a strong capability to sustain long-term accuracy and stability.

Despite these promising results, several open issues remain. The mechanism underlying the observed accuracy improvement from the ELM to the RBF-ELM requires further investigation. In particular, while the conditioning of the system matrix remains largely unchanged, the prediction accuracy is significantly enhanced. One potential explanation is that the localized support of RBF basis functions reduces interactions between distant centers, which in turn weakens global coupling in the approximation and thereby improves accuracy. This phenomenon will be investigated in future work to gain a deeper understanding of the underlying approximation behavior.
    
\section*{Acknowledgments}
	This research was supported by the National Key R$\&$D Program of China (2024YFA1012600) and NSFC Project (12431014).
	
\bibliographystyle{plain}
\bibliography{references}
	
\end{document}